\documentclass[twoside]{article} 

\usepackage{jwn-article}
\usepackage{fullpage}
\usepackage{amsfonts}
\usepackage{mathptmx}
\usepackage{hyperref}

\usepackage{graphicx}
\usepackage{moreverb}
\usepackage{alltt}
\usepackage[all]{xy}

\usepackage{amssymb}
\renewcommand{\ge}{\geqslant}
\renewcommand{\le}{\leqslant}

\newcommand{\cmd}[1]{\mathrm{#1}}
\newcommand{\prob}[2][]{\ensuremath{\cmd{Pr}_{#1}\left(#2\right)}}

\newcommand{\condprob}[3][]{\prob[#1]{#2\,|\,#3}}
\newcommand{\condutil}[3][]{\util[#1]{#2\,|\,#3}}

\newcommand{\embed}[1]{\ensuremath{\left(#1\right)}}
\newcommand{\state}[1]{\ensuremath{\mathsf{#1}}}
\newenvironment{code}{\begin{quote}\raggedright\tt}{\end{quote}}

\usepackage{stmaryrd}
\newcommand{\punless}[1]{\ensuremath{\;\bbslash\; #1}}

\renewcommand{\condutil}[3][]{\ensuremath{#2}}

\newcommand{\model}[1]{\ensuremath{\mathcal{M}_{\mathrm{#1}}}}
\newcommand{\mpq}{\model{PQ}}

\title{A Tutorial Introduction to the Logic of Parametric Probability}

\author{Joseph W. Norman \\ University of Michigan \\ 
\texttt{jwnorman@umich.edu}}
\date{May, 2012}

\begin{document}

\maketitle

\begin{abstract}
  The computational method of parametric probability analysis is
  introduced.  It is demonstrated how to embed logical formulas from
  the propositional calculus into parametric probability networks,
  thereby enabling sound reasoning about the probabilities of logical
  propositions.  An alternative direct probability encoding scheme is
  presented, which allows statements of implication and quantification
  to be modeled directly as constraints on conditional probabilities.
  Several example problems are solved, from Johnson-Laird's aces to
  Smullyan's zombies.  Many apparently challenging problems in logic
  turn out to be simple problems in algebra and computer science:
  systems of polynomial equations or linear optimization problems.
  This work extends the mathematical logic and parametric probability
  methods invented by George Boole.
\end{abstract}



\section{Introduction}

This essay introduces \emph{parametric probability analysis}, a method
to compute useful symbolic and numeric results from probability models
that contain parameters treated algebraically as unknown variables.  A
convention is provided to embed formulas from the propositional
calculus into such parametric probability models, thereby enabling
sound reasoning about the probabilities of logical propositions.  An
alternative scheme of direct probability encoding is presented, which
allows statements of implication and consequence to be modeled
directly as constraints on conditional probabilities (without
intermediate formulas from the propositional calculus).  With direct
encoding, probabilities can be used to extend classical logical
quantifiers into more precise proportional statements of
quantification (or even arbitrary polynomial constraints involving
such proportions).  Several example problems are analyzed, using the
Probability Query Language (PQL) computer program developed by the
author.  It turns out that many apparently challenging problems in
logic and probability are in fact simple problems in algebra and
computer science: systems of polynomial equations and inequalities,
general search problems, polynomial fractional optimization problems,
and in some cases just linear optimization problems.

This work is a continuation of George Boole's pioneering formulation
of mathematical logic and probability, codified in his 1854 \emph{Laws
  of Thought} \cite{boole}.  It complements the formalization of
Hailperin \cite{hailperin} by parsing Boole's notation in a
substantially different way that respects Boole's overloaded use of
operator signs and numerals.  Recognizing Boole's operator
overloading, as did Venn many years ago \cite{venn}, obviates the need
to invoke unusual arithmetic or heaps that are not quite sets in order
to explain Boole's calculations.  There were several innovations in
Boole's methodology: the representation of logical formulas and axioms
as polynomial formulas and equations; a means to embed logical
formulas within probability models; a database-and-query model of
interaction; and two-phase inference, with a primary phase of symbolic
probability inference followed by a secondary phase of more general
algebraic and numerical analysis.  The computational method introduced
here adds several features to extend Boole's original work: explicit
probability-network models; structured probability queries; clearer
semantics for embedding propositional-calculus formulas versus
directly encoding implication with conditional probabilities; and
broadened secondary analysis that includes search and general algebra
as well as optimization.

\section{Parametric Probability Networks}

\subsection{The Problem with Penguins}

Let us begin with a problem that has vexed quite a few philosophers
and computer scientists.  How can logic be used to reason about an
implication that is true sometimes but not always?  Following
artificial-intelligence tradition \cite{reiter} we contemplate the
problem that most birds can fly but some cannot.  To model this
problem we shall build a \emph{parametric probability network} denoted
\mpq.  This modeling formalism is built from the symbolic algebra used
by de~Moivre and Bernoulli in their foundational 18th-century
treatises on probability \cite{demoivre,bernoulli}; from the
parametric treatment of probability developed by Boole in the 19th
century \cite{boole}; from the axioms of probability theory provided
by Kolmogorov in the early 20th century \cite{kolmogorov}; from the
relational databases developed by Codd in the 1970s \cite{codd}; and
from the Bayesian belief networks and influence diagrams developed in
the 1980s by Pearl, Howard, and others \cite{pearl,howard-book}.  A
parametric probability network has four parts: a set of variables, a
set of constraints, a network graph, and a set of component
probability tables.  These parts can be described in a formal,
structured language that is suitable for processing by computers as
well as by humans.

\subsection{Variables and Constraints}

First let us introduce two \emph{primary variables} to represent
logical propositions about a hypothetical creature: $P$ that it is a
bird, and $Q$ that it can fly.  For this example each of these primary
variables may be either true or false, abbreviated \state{T} and
\state{F}.  We would like to consider the truth value of the logical
statement `$P$ implies $Q$' and its relation to various probabilities
involving $P$ and $Q$.  To facilitate this we add a third primary
variable $R$ defined as the value of the formula $P \rightarrow Q$,
where the arrow denotes the usual `if/then' material-implication
operator of propositional logic \cite{delong}.  This definition is
denoted $R:=(P \rightarrow Q)$, with the custom that the expression
after the definition sign uses another mathematical system that is
\emph{embedded} within the probability model (in this case the
propositional calculus).  Following Nilsson, one way to describe the
truth value of a logical proposition in the context of probability is
to use the probability that the proposition is true \cite{nilsson};
let us call this `fractional truth value'.  Thus $\prob{P=\state{T}}$
describes the fractional truth value of the atomic formula $P$;
$\prob{Q=\state{T}}$ describes the fractional truth value of the
atomic formula $Q$; and $\prob{R=\state{T}}$ describes the fractional
truth value of the compound formula $P \rightarrow Q$.

\begin{table}
\begin{tabular}[c]{llll}\hline
\multicolumn{1}{l}{\scshape {Variable}} & 
\multicolumn{1}{l}{\scshape {Role}} & 
\multicolumn{1}{l}{\scshape {Description}} & 
\multicolumn{1}{l}{\scshape {Domain}} \\ \hline\hline
$P$ & 
Primary & 
Proposition: It is a bird & 
$\{ {\state{T}, \state{F}} \}$ \\ \hline
$Q$ & 
Primary & 
Proposition: It can fly & 
$\{ {\state{T}, \state{F}} \}$ \\ \hline
$R$ & 
Primary & 
Value of $(P \rightarrow Q)$ & 
$\{ {\state{T}, \state{F}} \}$ \\ \hline
$A$ & 
Primary & 
Value of $(3)$ & 
$\{ {\state{3}} \}$ \\ \hline
$B$ & 
Primary & 
Number true in $\{P,Q,R\}$ & 
$\{ {\state{0}, \state{1}, \state{2}, \state{3}} \}$ \\ \hline
$C$ & 
Primary & 
Value of $(A - B)$ & 
$\{ {\state{0}, \state{1}, \state{2}, \state{3}} \}$ \\ \hline
$x$ & 
Parameter & 
Probability that $P$ is true & 
$[0,1]$ \\ \hline
$y$ & 
Parameter & 
Probability that $Q$ is true if $P$ is true & 
$[0,1]$ \\ \hline
$z$ & 
Parameter & 
Probability that $Q$ is true if $P$ is false & 
$[0,1]$ \\ \hline
\end{tabular}
\caption{Variables in the probability network \mpq{} describing a
  creature that might be a bird and might be able to fly.}
\label{tbl:m1-variables}
\end{table}

Next let us introduce three more variables, $x$, $y$, and $z$, to be
used as \emph{parameters} for specifying the probabilities of $P$ and
$Q$.  Since parameters are treated differently from primary variables
during analysis, we maintain a distinction between these two
\emph{roles} that a variable may play.  In contrast to $P$, $Q$, and
$R$, which share the domain $\{\state{T},\state{F}\}$ of two possible
values, the parameters $x$, $y$, and $z$ may take real-number values
between zero and one: thus the domain of each is the interval $[0,1]$.
Finally let us add three more primary variables $A$, $B$, and $C$,
with the declarations that $B$ and $C$ share the domain $\{0,1,2,3\}$
of four possible integer values and that $A$ has the set $\{3\}$ of a
solitary possible value.  We define $A$ as the value of the number
$3$; $B$ as the number of true propositions in the set $\{P,Q,R\}$;
and $C$ as the value of the difference $A-B$.  Thus $A:=3$ and
$C:=(A-B)$, now using integer arithmetic instead of propositional
logic as the embedded mathematical system.
Table~\ref{tbl:m1-variables} lists all the variables in the parametric
probability network \mpq{}.  Using $m$ for the number of primary
variables and $n$ for the number of parameters, this model has $m=6$
and $n=3$.

This document uses the typographical conventions that primary
variables are rendered as uppercase italic letters ($A$, $P$, $Q$,
etc.) and parameters as lowercase italic letters ($x$, $y$, etc.);
variable names are case-sensitive, so for example $x$ and $X$ are
considered different variables.  States of primary variables are
rendered in sans-serif type (such as \state{True}, \state{False},
\state{T}, \state{F}, \state{0}, \state{1}, \state{2}, \state{3}).
These practices are intended to distinguish numbers and formulas in
embedded mathematical systems (such as `$P \rightarrow Q$' which is a
formula in the propositional calculus) from numbers and formulas in
the host probability model (such as `$1-p+pq$' which is a formula in
the algebra of polynomials with rational coefficients).  Additionally,
different symbols are used to disambiguate different meanings of the
equal sign: the colon and equal sign $:=$ for definition or
assignment; the double right arrow $\Rightarrow$ for evaluation (as in
$2+2 \Rightarrow 4$); and the standard equal sign $=$ for the test or
assertion of equality.  For this presentation, all primary variables
are discrete (with finite sets of possible values) and all parameters
are continuous (taking rational or real-number values).  But in
general, parametric probability networks are allowed to have
continuous primary variables and discrete parameters too.

\subsection{Network Graph and Component Probability Tables}
\label{sec:graph}

\begin{figure}
\includegraphics{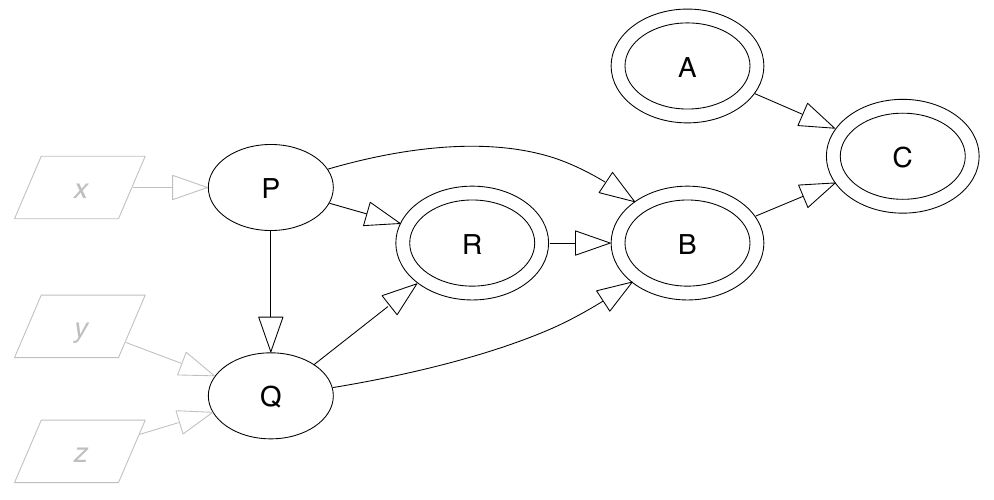}
\caption{Probability network graph of model \mpq, with
  idiosyncratic graphical notation explained in
  Section~\ref{sec:graph}.}
\label{fig:m1-id}
\end{figure}

The graph in Figure~\ref{fig:m1-id} shows how the variables in the
parametric probability network \mpq{} are related to one another; it
also dictates which component probability tables must be specified in
order to complete the model.  Following standard probability-network
notation, each oval node in this directed acyclic graph represents a
primary variable and each directed edge shows a correlation or
functional dependence relationship.  The absence of an edge is an
assertion of independence.  In the author's idiosyncratic notation
parameters are included in the graph and drawn with parallelogram
nodes; furthermore a clique (fully-connected subset of nodes) will be
indicated by a small diamond and undirected edges (as illustrated in
the examples in Section~\ref{sec:examples}).  Double borders on a node
indicate that the corresponding variable is \emph{deterministic}: for
a primary variable this means that every component probability must be
either $0$ or $1$; for a parameter this means that its value must be
fixed at some constant.

\begin{table}
(a) \quad
\begin{tabular}[c]{l|l}\hline
\multicolumn{1}{l|}{$P$} & 
\multicolumn{1}{l}{$\prob[0]{{P}}$} \\ \hline\hline
$\state{T}$ & 
$x$ \\ \hline
$\state{F}$ & 
$1 - x$ \\ \hline
\end{tabular}
\qquad (b) \quad
\begin{tabular}[c]{l|ll}\hline
\multicolumn{3}{l}{$\condprob[0]{{Q}}{{P}}$} \\ \hline\hline
\multicolumn{1}{l|}{$P$} & 
\multicolumn{1}{l}{${Q=\state{T}}$} & 
\multicolumn{1}{l}{${Q=\state{F}}$} \\ \hline\hline
$\state{T}$ & 
$y$ & 
$1 - y$ \\ \hline
$\state{F}$ & 
$z$ & 
$1 - z$ \\ \hline
\end{tabular}
\qquad (c) \quad
\begin{tabular}[c]{ll|ll}\hline
\multicolumn{4}{l}{$\condprob[0]{{R}}{{P, Q}}$} \\ \hline\hline
\multicolumn{1}{l}{$P$} & 
\multicolumn{1}{l|}{$Q$} & 
\multicolumn{1}{l}{${R=\state{T}}$} & 
\multicolumn{1}{l}{${R=\state{F}}$} \\ \hline\hline
$\state{T}$ & 
$\state{T}$ & 
$1$ & 
$0$ \\ \hline
$\state{T}$ & 
$\state{F}$ & 
$0$ & 
$1$ \\ \hline
$\state{F}$ & 
$\state{T}$ & 
$1$ & 
$0$ \\ \hline
$\state{F}$ & 
$\state{F}$ & 
$1$ & 
$0$ \\ \hline
\end{tabular}
\\\\
(d) \quad
\begin{tabular}[c]{lll|llll}\hline
\multicolumn{7}{l}{$\condprob[0]{{B}}{{P, Q, R}}$} \\ \hline\hline
\multicolumn{1}{l}{$P$} & 
\multicolumn{1}{l}{$Q$} & 
\multicolumn{1}{l|}{$R$} & 
\multicolumn{1}{l}{${B=\state{0}}$} & 
\multicolumn{1}{l}{${B=\state{1}}$} & 
\multicolumn{1}{l}{${B=\state{2}}$} & 
\multicolumn{1}{l}{${B=\state{3}}$} \\ \hline\hline
$\state{T}$ & 
$\state{T}$ & 
$\state{T}$ & 
$0$ & 
$0$ & 
$0$ & 
$1$ \\ \hline
$\state{T}$ & 
$\state{T}$ & 
$\state{F}$ & 
$0$ & 
$0$ & 
$1$ & 
$0$ \\ \hline
$\state{T}$ & 
$\state{F}$ & 
$\state{T}$ & 
$0$ & 
$0$ & 
$1$ & 
$0$ \\ \hline
$\state{T}$ & 
$\state{F}$ & 
$\state{F}$ & 
$0$ & 
$1$ & 
$0$ & 
$0$ \\ \hline
$\state{F}$ & 
$\state{T}$ & 
$\state{T}$ & 
$0$ & 
$0$ & 
$1$ & 
$0$ \\ \hline
$\state{F}$ & 
$\state{T}$ & 
$\state{F}$ & 
$0$ & 
$1$ & 
$0$ & 
$0$ \\ \hline
$\state{F}$ & 
$\state{F}$ & 
$\state{T}$ & 
$0$ & 
$1$ & 
$0$ & 
$0$ \\ \hline
$\state{F}$ & 
$\state{F}$ & 
$\state{F}$ & 
$1$ & 
$0$ & 
$0$ & 
$0$ \\ \hline
\end{tabular}
\qquad
\begin{tabular}{@{}l@{}}
(e) \quad
\begin{tabular}[c]{l|l}\hline
\multicolumn{1}{l|}{$A$} & 
\multicolumn{1}{l}{$\prob[0]{{A}}$} \\ \hline\hline
$\state{3}$ & 
$1$ \\ \hline
\end{tabular}
\\\\
(f) \quad
\begin{tabular}[c]{ll|llll}\hline
\multicolumn{6}{l}{$\condprob[0]{{C}}{{A, B}}$} \\ \hline\hline
\multicolumn{1}{l}{$A$} & 
\multicolumn{1}{l|}{$B$} & 
\multicolumn{1}{l}{${C=\state{0}}$} & 
\multicolumn{1}{l}{${C=\state{1}}$} & 
\multicolumn{1}{l}{${C=\state{2}}$} & 
\multicolumn{1}{l}{${C=\state{3}}$} \\ \hline\hline
$\state{3}$ & 
$\state{0}$ & 
$0$ & 
$0$ & 
$0$ & 
$1$ \\ \hline
$\state{3}$ & 
$\state{1}$ & 
$0$ & 
$0$ & 
$1$ & 
$0$ \\ \hline
$\state{3}$ & 
$\state{2}$ & 
$0$ & 
$1$ & 
$0$ & 
$0$ \\ \hline
$\state{3}$ & 
$\state{3}$ & 
$1$ & 
$0$ & 
$0$ & 
$0$ \\ \hline
\end{tabular}
\end{tabular}
\caption{Component probability tables for the model \mpq.  The
  subscript in \prob[0]{\cdots} identifies these as user input.}
\label{tbl:m1-cpt}
\end{table}

The probability of each primary variable must be specified as a
function of its parents in the network graph.  Thus for the model
\mpq{} we must provide several component probabilities based on the
graph in Figure~\ref{fig:m1-id}.  We must specify the probability of
the primary variable $P$ as a function of its parent, the parameter
$x$; and the conditional probability of $Q$ given its primary-variable
parent $P$ as a function of its parameter parents $y$ and $z$.  These
component probability distributions appear in Table~\ref{tbl:m1-cpt}
parts (a) and (b).  We must also specify the conditional probability
of $R$ given its parents $P$ and $Q$.  For this we transcribe the
truth table of the logical formula $P \rightarrow Q$ into the
conditional probability table shown as Table~\ref{tbl:m1-cpt}
part~(c).  The conditional probabilities for $B$ given $P$, $Q$, and
$R$, shown as Table~\ref{tbl:m1-cpt} part~(d), encode the number of
true parent variables.  The primary variable $A$ has one possible
state which is assigned probability one as shown in
Table~\ref{tbl:m1-cpt} part~(e).  Finally, for the primary variable
$C$ we specify as Table~\ref{tbl:m1-cpt} part~(f) a transcription of
the table that gives the value of the arithmetical formula $A-B$ using
$A=3$ and integer values of $B$ between $0$ and $3$.  Component
probabilities specified by the user are designated \prob[0]{\ldots},
with the subscript $0$ used to distinguish these input values from the
output probabilities \prob{\ldots} later computed from them.

In general the component probability table for a primary variable must
contain an element for each of its possible states, given every unique
combination of states of its primary-variable parents (an empty set of
parents is considered to have one combination of states).  Each of
these component probabilities must be a polynomial function of the
model parameters (with real or rational coefficients).  The user may
specify arbitrary polynomial equality and inequality constraints on
the model parameters; the system adds constraints as needed to enforce
the laws of probability (that the feasible values of each component
probability must lie between zero and one, and that at every feasible
point the probabilities of mutually exclusive and collectively
exhaustive events must add up to one).  Parameters do not get their
own probability distributions; that is precisely how they differ from
primary variables.  Hence in \mpq{} there are no component probability
tables for the parameters $x$, $y$, and $z$.  However parameters are
always subject to algebraic constraints, in this case the zero-one
bounds given in Table~\ref{tbl:m1-variables}: thus $0 \le x \le 1$, $0
\le y \le 1$, and $0 \le z \le 1$.

\subsection{Computable Specification}

Here is the computable specification of the parametric probability
model \mpq.  The corresponding file was processed by the author's
computer program
to generate the figures and tables in this section.
\begin{code}
\begin{verbatim}
// basic1.pql: with propositional calculus and integer arithmetic embedded

parameter x { label = "Probability that $P$ is true"; range = (0,1); }
parameter y { label = "Probability that $Q$ is true if $P$ is true"; range = (0,1); }
parameter z { label = "Probability that $Q$ is true if $P$ is false"; range = (0,1); }

primary P { label = "Proposition: It is a bird"; states = binary; }
probability ( P ) { data = (x, 1-x); noverify; }

primary Q { label = "Proposition: It can fly"; states = binary; }
probability ( Q | P ) { data = (y, 1-y, z, 1-z); noverify; }

primary R { label = "Value of $(P \rightarrow Q)$"; states = binary; }
probability ( R | P Q ) { function = "R <-> P -> Q ? 1 : 0"; }

primary A { label = "Value of $(3)$"; states = range( 3, 3 ); }
probability ( A ) { data = (1); }

primary B { label = "Number true in $\{P,Q,R\}$"; states = range( 0, 3 ); }
probability ( B | P Q R ) { function = "B == P + Q + R ? 1 : 0"; }

primary C { label = "Value of $(A - B)$"; states = range( 0, 3 ); }
probability ( C | A B ) { function = "C == A - B ? 1 : 0"; }

net { graph = 'subgraph { rank=same; "P"; "Q"; }'; }  // hint for graph drawing
\end{verbatim}
\end{code}
Like the Structured Query Language (SQL) for relational databases
\cite{date}, the author's Probability Query Language includes a data
definition language for specifying models and a data manipulation
language for asking queries.  The PQL data definition language, which
you see illustrated above, has syntax like the ubiquitous C
programming language \cite{c} and was also inspired by the
\textsc{net} modeling language from the Hugin system for
Bayesian-network inference \cite{hugin-api}.  The PQL data
manipulation language, which is demonstrated below, uses keywords that
the user may type as commands to an interactive shell or alternatively
include in source files for a batch processor.  The command-line shell
(called \texttt{pqlsh}) and the batch processor (called
\texttt{pqlpp}) are built atop the Tcl scripting language \cite{tcl}.
There is implicitly a structured language for the \emph{results} of
PQL queries as well; these results are essentially relational-database
tables whose entries include natural numbers, character strings, and
fractional polynomials with rational coefficients.

\section{Embedded Mathematical Systems}
\label{sec:embedding}

Probability can be used to reason about formulas that are governed by
other mathematical systems such as the propositional calculus or
integer arithmetic.  Such mathematical formulas can be \emph{embedded}
within parametric probability networks using the method described
here.  There are three main aspects to embedding mathematical formulas
in probability networks: assigning prior probabilities to the
probability-network variables copied from the mathematical variables;
assigning conditional probabilities to the probability-network
variables created to represent the non-atomic mathematical formulas;
and modeling unknown mathematical functions.

\subsection{Prior Probabilities for Propositional Variables}
\label{sec:prior}

To begin the embedding process we copy the variables used in the
mathematical formulas of interest into primary variables in the
parametric probability network; then we assign these primary variables
a prior probability distribution.  Let us take the simple view that a
mathematical formula $\phi$ is a finite-length string that may contain
only constants, variables, operator signs, and parentheses; a set of
grammatical rules dictates which strings constitute well-formed
formulas.  The constants must be members of some set $K$ of elementary
values, and each variable $P_1,P_2,\ldots,P_v$ ranges over values in
this set $K$.  When each variable in a formula is assigned a constant
value from $K$, the formula itself must also have a value in $K$; thus
we imagine the set $K$ to be the domain of a mathematical structure
that is closed under the allowed operations.  For this discussion, we
assume that each formula $\phi$ is finite in length and that the
number of formulas under consideration is finite; hence the number $v$
of primary variables is also finite.  Let us further limit our
attention to the case that the set $K$ of elementary values is finite,
with some size $d$; this is adequate for modeling `logical' systems
with limited numbers of elementary truth values (for example the usual
two).

\paragraph{Uninformative Priors, Parametrically}

The simplest way to specify a prior probability distribution on the
variables $P_1$ through $P_v$ is to avoid independence assertions and
to specify directly the joint probability distribution
$\prob[0]{P_1,P_2,\ldots,P_v}$.  Recall that the subscript $0$ in
\prob[0]{\cdots} designates component probabilities input by the user,
as opposed to computed probabilities \prob{\cdots} output by the
system.  Anyway in this joint-prior case the corresponding
probability-network graph has the variables $P_1$ through $P_v$ joined
as a clique.  With $d$ elementary values in the set $K$ it is
necessary to provide $d^v$ individual probabilities in order to
specify completely the distribution $\prob[0]{P_1,P_2,\ldots,P_v}$.
The laws of probability require that each probability must lie between
zero and one, and that the sum of the probabilities in this joint
distribution must be one (since they describe mutually exclusive and
collectively exhaustive events).  To provide a truly
\emph{uninformative} prior probability distribution, we should use
parameters to state exactly these properties and nothing more.  For
example, we might specify the respective probabilities as the
variables $x_1,x_2,\ldots,x_{d^v}$ subject to the constraints that
each $x_i \ge 0$, each $x_i \le 1$, and the sum $\sum_i x_i = 1$.

It is also possible to specify the joint distribution
$\prob[0]{P_1,P_2,\ldots,P_v}$ indirectly by using a
probability-network subgraph with directed edges; if this subgraph is
fully-connected then it still does not introduce independence
assertions.  For example we might specify the unconditioned
probability distribution $\prob[0]{P_1}$ and then the conditional
probability distributions $\condprob[0]{P_2}{P_1}$,
$\condprob[0]{P_3}{P_1,P_2}$, and so on until
$\condprob[0]{P_v}{P_1,P_2,\ldots,P_{v-1}}$.  De~Moivre used exactly
this construction centuries ago to describe the joint probability of
several dependent events \cite{demoivre}.  In general many different
fully-connected network graphs are possible; each requires a
particular set of probabilities to be supplied by the user.

\paragraph{Just Add Information}

There are several ways to add information beyond uninformative
parametric prior probabilities, should the user wish to do so: through
the choice of values in the component probability tables; through
explicit algebraic constraints on the parameters used to specify
component probabilities; and through independence assertions (modeled
as the absence of arcs in the probability-network graph) which
indirectly provide algebraic constraints.  For example, the user may
desire to specify that the primary variables $P_1,P_2,\ldots,P_v$ are
probabilistically independent: in this case they are not directly
connected in the network graph, and the joint probability
\prob{X_1=k_1,X_2=k_2,\ldots,X_v=k_v} that each variable $X_i$ takes
the elementary value $k_i$ is constrained to equal the product
of the concordant individual component
probabilities:
\begin{math}
  \prob[0]{X_1=k_1} \times \prob[0]{X_2=k_2} \times \cdots \times
  \prob[0]{X_v=k_v}
\end{math}.
As another example, it may be desirable to constrain each prior
probability to be strictly greater than zero, in order to specify that
no combination of primary-variable values is considered impossible
\emph{a priori}.  The tools of parametric probability, including
graphical models and algebraic constraints, allow the user to say
exactly what he or she means about prior probabilities.

\paragraph{For the Birds}

To illustrate copied mathematical variables and prior probabilities,
in the \mpq{} model there are two embedded formulas $P \rightarrow Q$
and $A-B$.  The first formula uses the propositional calculus, which
includes the set $\{\state{T},\state{F}\}$ of elementary values
representing truth and falsity as well as the operator $\rightarrow$
for material implication (along with the usual operators $\wedge$,
$\vee$, $\neg$, etc.).  We copy the propositional variables $P$ and
$Q$ into the probability network as primary variables, each with the
domain $\{\state{T},\state{F}\}$.  In lieu of the joint distribution
$\prob[0]{P,Q}$, for this model we specify the component probabilities
indicated by the fully-connected subgraph with a directed arc from $P$
to $Q$; hence the required component probabilities are $\prob[0]{P}$
and $\condprob[0]{Q}{P}$ which appear in Table~\ref{tbl:m1-cpt}.  The
constraints on the real-valued parameters $x$, $y$, and $z$ used in
these component probability tables provide no more information than
the laws of probability require.  The value of each parameter must lie
between $0$ and $1$.  Here algebra takes care of the sum-to-one
constraints; for example in \prob[0]{P} it is tautological that
$x+(1-x)=1$.  The second embedded formula $A-B$ uses integer
arithmetic, which includes set $\{0,1,2,3\}$ of elementary values (for
convenience we focus on this finite subset of $\mathbb{Z}$) and the
operator $-$ for subtraction (along with the usual operators $+$,
$\times$, etc.).  The variables $A$ and $B$ are copied into the
probability network as primary variables.  For the definition $A:=3$
we add the prior information that the only possible state of $A$ is
\state{3}; hence the simple component probability table in
Table~\ref{tbl:m1-cpt} part~(e) assigns probability one to this event.
The component probability table specified in Table~\ref{tbl:m1-cpt}
part~(d) encodes that the state of $B$ expresses the number of primary
variables among $P$, $Q$, and $R$ that have the state \state{T}.

\subsection{Conditional Probabilities for Operations}
\label{sec:conditional}

To continue the embedding process we introduce additional primary
variables for the compound formulas of interest.  For each non-atomic
formula $\phi_i$ we introduce a new primary variable $S_i$.  The
input component probability table for $S_i$ must specify the
conditional probability of $S_i$ given the variables
$P_{j_1},P_{j_2},\ldots,P_{j_a}$ used in the formula $\phi_i$.  This
conditional probability table
$\condprob[0]{S_i}{P_{j_1},P_{j_2},\ldots,P_{j_a}}$ must assign
probability one to the appropriate value of the formula given each
combination of values of its arguments.

Returning to the \mpq{} model, we introduce the primary variable $R$
to represent the value of the compound logical formula $P \rightarrow
Q$ and we add the primary variable $C$ for the compound arithmetical
formula $A-B$.  To complete the definition $R:=(P \rightarrow Q)$ the
input component probability table $\condprob[0]{R}{P,Q}$ assigns
probability one to the appropriate elementary truth value of the
statement of material implication, given each combination of truth
values of its arguments.  The conditional probability table is derived
in the obvious way from the related logical truth table:
\begin{eqnarray}
  \label{eq:p-impl-q}
  \begin{tabular}{cc|c}
    \hline
    $P$ & $Q$ & $P \rightarrow Q$ \\ \hline\hline
    \state{T} & \state{T} & \state{T} \\ \hline
    \state{T} & \state{F} & \state{F} \\ \hline
    \state{F} & \state{T} & \state{T} \\ \hline
    \state{F} & \state{F} & \state{T} \\ \hline
  \end{tabular}
  & \leadsto &
\begin{tabular}[c]{ll|ll}\hline
\multicolumn{4}{l}{$\condprob[0]{{\mbox{\embed{P \rightarrow Q}}}}{{\mbox{$P$}, \mbox{$Q$}}}$} \\ \hline\hline
\multicolumn{1}{l}{$\mbox{$P$}$} & 
\multicolumn{1}{l|}{$\mbox{$Q$}$} & 
\multicolumn{1}{l}{${\mbox{\embed{P \rightarrow Q}}=\state{T}}$} & 
\multicolumn{1}{l}{${\mbox{\embed{P \rightarrow Q}}=\state{F}}$} \\ \hline\hline
$\state{T}$ & 
$\state{T}$ & 
$1$ & 
$0$ \\ \hline
$\state{T}$ & 
$\state{F}$ & 
$0$ & 
$1$ \\ \hline
$\state{F}$ & 
$\state{T}$ & 
$1$ & 
$0$ \\ \hline
$\state{F}$ & 
$\state{F}$ & 
$1$ & 
$0$ \\ \hline
\end{tabular}
\end{eqnarray}
For example because the logical formula
$\state{T}\rightarrow\state{F}$ has the value $\state{F}$ we have
assigned the corresponding conditional probability $\condprob[0]{(P
  \rightarrow Q)=\state{F}}{P=\state{T},Q=\state{F}}$ the value $1$
and its complement $\condprob[0]{(P \rightarrow
  Q)=\state{T}}{P=\state{T},Q=\state{F}}$ the value $0$.  In
Equation~\ref{eq:p-impl-q} the complete formula $P \rightarrow Q$
instead of its abbreviation $R$ appears in the heading of the table
\condprob[0]{(P \rightarrow Q)}{P,Q}; otherwise this component
probability table is the same as \condprob[0]{R}{P,Q} which appears as
Table~\ref{tbl:m1-cpt} part~(c).  The other definition $C:=(A-B)$ in
the \mpq{} model is handled essentially the same way.  We consider the
portion of the mathematical function table for the arithmetical
subtraction operator when its first argument is $3$ and its second
argument is a member of the set $\{0,1,2,3\}$.  The derived component
probability table follows:
\begin{eqnarray}
  \begin{tabular}{cc|c}
    \hline
    $A$ & $B$ & $A-B$ \\ \hline\hline
    \state{3} & \state{0} & \state{3} \\ \hline
    \state{3} & \state{1} & \state{2} \\ \hline
    \state{3} & \state{2} & \state{1} \\ \hline
    \state{3} & \state{3} & \state{0} \\ \hline
  \end{tabular}
  & \leadsto &
\begin{tabular}[c]{ll|llll}\hline
\multicolumn{6}{l}{$\condprob[0]{{\mbox{$(A-B)$}}}{{\mbox{$A$}, \mbox{$B$}}}$} \\ \hline\hline
\multicolumn{1}{l}{$\mbox{$A$}$} & 
\multicolumn{1}{l|}{$\mbox{$B$}$} & 
\multicolumn{1}{l}{${\mbox{$(A-B)$}=\state{0}}$} & 
\multicolumn{1}{l}{${\mbox{$(A-B)$}=\state{1}}$} & 
\multicolumn{1}{l}{${\mbox{$(A-B)$}=\state{2}}$} & 
\multicolumn{1}{l}{${\mbox{$(A-B)$}=\state{3}}$} \\ \hline\hline
$\state{3}$ & 
$\state{0}$ & 
$0$ & 
$0$ & 
$0$ & 
$1$ \\ \hline
$\state{3}$ & 
$\state{1}$ & 
$0$ & 
$0$ & 
$1$ & 
$0$ \\ \hline
$\state{3}$ & 
$\state{2}$ & 
$0$ & 
$1$ & 
$0$ & 
$0$ \\ \hline
$\state{3}$ & 
$\state{3}$ & 
$1$ & 
$0$ & 
$0$ & 
$0$ \\ \hline
\end{tabular}
\end{eqnarray}
The same table appears as Table~\ref{tbl:m1-cpt} part~(f) labeled with
the abbreviation $C$ instead of the full embedded formula $A-B$.

\subsection{Known and Unknown Unknowns}
\label{sec:unknown}

The number of elementary truth values in a logical system is
independent of the idea of assigning a probability to each possible
value (or the idea of considering \emph{sets} of elementary values).
For example in ordinary algebra, in order to express the idea that an
integer is unknown, we do not imagine that there is an integer called
`unknown' in the set $\mathbb{Z}$.  Instead we introduce a symbolic
variable and declare that it is constrained to take values in the set
of integers, for example $Y$ with $Y \in \mathbb{Z}$.  We can expand
this set-based declaration $Y \in \mathbb{Z}$ into a parametric
probability distribution by further specifying that there is some
probability $p_k$ that the variable $Y$ takes each integer value $k$
in the set $\mathbb{Z}$.  In this context, to express perfect
ignorance about the value of the integer-valued variable $Y$, we
should admit only that each probability $p_k$ takes a real value
between zero and one and that all of the probabilities add up to one:
hence we constrain each $p_k \in \mathbb{R}$ with $0 \le p_k \le 1$
and the infinite sum $\sum_k p_k = 1$.

Returning to logic, let us consider an embedded system similar to the
propositional calculus but with the set
$\{\state{True},\state{False},\state{Unknown}\}$ of three elementary
truth values instead of the usual two.  Now, in this $3$-valued logic
it would be a different thing to say that the value of some variable
$V$ is unknown than to say that the value of $V$ is certainly the
elementary value called \state{Unknown}.  For the former assertion
(the value of $V$ is unknown) we should start with the set-based
declaration $V \in \{\state{True},\state{False},\state{Unknown}\}$ and
optionally expand this declaration into a parametric probability
distribution with only the essential probability constraints, such as:
\begin{equation}
\begin{tabular}[c]{l|l}\hline
\multicolumn{1}{l|}{$V$} & 
\multicolumn{1}{l}{$\prob[0]{{V}}$} \\ \hline\hline
$\state{True}$ & 
$x_{1}$ \\ \hline
$\state{False}$ & 
$x_{2}$ \\ \hline
$\state{Unknown}$ & 
$x_{3}$ \\ \hline
\end{tabular}
  \qquad x_i \in \mathbb{R}, \quad 0 \le x_i \le 1, 
  \quad x_1 + x_2 + x_3 = 1
\end{equation}
But for the latter assertion (\state{Unknown} is the value of $V$) we
should assign probability one to the value named \state{Unknown}:
\begin{equation}
\begin{tabular}[c]{l|l}\hline
\multicolumn{1}{l|}{$V$} & 
\multicolumn{1}{l}{$\prob[0]{{V}}$} \\ \hline\hline
$\state{True}$ & 
$0$ \\ \hline
$\state{False}$ & 
$0$ \\ \hline
$\state{Unknown}$ & 
$1$ \\ \hline
\end{tabular}
\end{equation}
Many multivalued logics confuse these two ideas (considering sets of
possible elementary values---perhaps with probabilities
attached---versus adding new elementary values).

We can extend the courtesy of parametric representation to unknown
\emph{functions} as well as to unknown \emph{variables}.  For example,
returning to $2$-valued logic, the following conditional probability
table and constraints describe an unknown binary operation $R^*$ whose
arguments $P$ and $Q$ take values in $\{\state{T},\state{F}\}$:
\begin{equation}
  \label{eq:sp}
\begin{tabular}[c]{ll|ll}\hline
\multicolumn{4}{l}{$\condprob[0]{{R^*}}{{P, Q}}$} \\ \hline\hline
\multicolumn{1}{l}{$P$} & 
\multicolumn{1}{l|}{$Q$} & 
\multicolumn{1}{l}{${R^*=\state{T}}$} & 
\multicolumn{1}{l}{${R^*=\state{F}}$} \\ \hline\hline
$\state{T}$ & 
$\state{T}$ & 
$t_{1}$ & 
$1 - t_{1}$ \\ \hline
$\state{T}$ & 
$\state{F}$ & 
$t_{2}$ & 
$1 - t_{2}$ \\ \hline
$\state{F}$ & 
$\state{T}$ & 
$t_{3}$ & 
$1 - t_{3}$ \\ \hline
$\state{F}$ & 
$\state{F}$ & 
$t_{4}$ & 
$1 - t_{4}$ \\ \hline
\end{tabular}
  \qquad
  t_i \in \mathbb{R}, \quad
  t_i \in \{0,1\}
\end{equation}
This parametric probability table encodes $2^4 \Rightarrow 16$
possible logical functions.  For example, the material-implication
function $R^*_{1011}:=(P \rightarrow Q)$ corresponds to the vector
$(1,0,1,1)$ of values for the parameters $(t_1,t_2,t_3,t_4)$; the
biconditional $R^*_{1001}:=(P \leftrightarrow Q)$ corresponds to the
vector $(1,0,0,1)$; and the always-false function
$R^*_{0000}:=\state{F}$ corresponds to the vector $(0,0,0,0)$.

\section{Primary Analysis: Symbolic Probability Inference}
\label{sec:primary}

Having defined the parametric probability network \mpq{} with embedded
formulas from the propositional calculus and from integer arithmetic,
let us now proceed with some analysis.  We follow the ingenious
framework laid out by Boole in his \emph{Laws of Thought}
\cite{boole}, which included a database-and-query model of interaction
between the user and the analytic system, as well as a two-phase model
of inference.  Boole considered a parametric probability model (the
\emph{data} in his terminology), to which a probability query could be
posed (his \emph{qu{\ae}situm}).  In the first phase of analysis, a
polynomial formula (Boole's \emph{final logical equation}) was
calculated to answer the query; the variables in this polynomial were
the unknown parameters in the probability model.  In the second phase
of analysis, the minimum and maximum feasible values of this
polynomial formula were computed (Boole's \emph{limits}), subject to
constraints that expressed the laws of probability.  Parametric
probability analysis follows Boole's two-phase model of inference: in
the primary phase of analysis, parametric probability networks and
structured probability queries are processed by \emph{symbolic
  probability-network inference} to compute polynomial solutions; and
in the secondary phase of analysis, these polynomial solutions are
used for additional algebraic and numerical analysis.

\subsection{Structured Probability Queries}
\label{sec:queries}

Each probability-table query asks for the probability distribution of
a \emph{principal} set of primary variables, given some
\emph{conditioning} set of primary variables; the remaining primary
variables are considered the \emph{marginal} set.  Any of these three
sets may be empty, but every primary variable must be assigned to one
of these positions to make a valid query.  Therefore for a parametric
probability network with $m$ primary variables, there are $3^m$
possible probability-table queries (thus $3^6 \Rightarrow 729$
possibilities for the \mpq{} model).  The result of a
probability-table query is a probability table with one or more
elements.  Each element in this result table is a polynomial function
of the component probabilities or a quotient of such polynomials; as a
special case an element can be a plain rational number.

For example, to ask for the fractional truth value of the proposition
$P \rightarrow Q$, which is represented in the model \mpq{} as the
variable $R$, we ask the probability-table query $\prob{R}$.  This
query has the principal set $\{R\}$; its conditioning set is empty;
its marginal set $\{P,Q,A,B,C\}$ contains the remaining primary
variables.  The following session shows how to use the \texttt{pqlsh}
command-line interface to load the probability model and to issue this
probability-table query; here \texttt{pqlsh>} is the system prompt,
user input is shown in italic type, and input follows the syntax of
the Tcl scripting language \cite{tcl}:
\begin{code}
\begin{alltt}
pqlsh> \emph{set m [pql::load basic1.pql t]; return;}
\end{alltt}
\begin{alltt}
pqlsh> \emph{set tr [\$m table R]; \$tr infer; \$tr print -index}
\end{alltt}
\begin{verbatimtab}
Index	| R	| Pr( {R} )	
-------	-------	-------
1	| T	| 1 - x + x*y	
2	| F	| x - x*y	
\end{verbatimtab}
\end{code}
Such results can also be generated in typeset form: PQL commands
included as special comments in \LaTeX{} documents are replaced by the
\texttt{pqlpp} preprocessor with the appropriate query results (that
is how this document was generated).  Anyway, inspecting the results
above and recalling the parameter definitions in
Table~\ref{tbl:m1-variables}, you can see that the fractional truth
value of the formula $P \rightarrow Q$ is a simple polynomial function
of the parameters $x$ and $y$, where $x$ is the input probability that
$P$ is true and $y$ is the input conditional probability that $Q$ is
true given that $P$ is true.  The bound constraints $x \in [0,1]$ and
$y \in [0,1]$ given in Table~\ref{tbl:m1-variables} are considered to
be part of the result (as would be additional parameter constraints if
there were any): thus the result
\begin{math}
\prob{{R}}
\end{math}
is the above probability table along with the constraints that the
parameters $x$ and $y$ take real values between zero and one.

Inspecting the first element
\begin{math}
1 - x + x y
\end{math}
in the result of the table query above, note that $\prob{R=\state{T}}
= 0$ if and only if $x=1$ and $y=0$ (considering also the preexisting
constraints $0 \le x \le 1$ and $0 \le y \le 1$).  In terms of the
example, the statement $P \rightarrow Q$ that bird implies flight is
certainly false exactly when the creature is certainly a bird ($x=1$)
and it is certain that no bird can fly ($y=0$).  Conversely, from
either constraint $x=0$ or $(x,y)=(1,1)$ it follows by simple algebra
that $\prob{R=\state{T}}=1$.  In other words, the implication $P
\rightarrow Q$ that bird implies flight is certainly true if the
creature is certainly not a bird ($x=0$); or if the creature is
certainly a bird and it is certain that a bird can fly ($x=1$ and
$y=1$).  Other values of $x$ and $y$ give intermediate values of
$\prob{R=\state{T}}$ and its complement $\prob{R=\state{F}}$.

Moving on, in order to ask the conditional probability of $Q$ given
$P$ we use the probability-table query $\condprob{Q}{P}$:
\begin{code}
\begin{alltt}
pqlsh> \emph{set tqp [[\$m table Q | P] infer]; \$tqp print -index}
\end{alltt}
\begin{verbatimtab}
Index	| P	Q	| Pr( {Q} | {P} )	
-------	-------	-------	-------
1	| T	T	| (x*y) / (x)	
2	| T	F	| (x - x*y) / (x)	
3	| F	T	| (z - x*z) / (1 - x)	
4	| F	F	| (1 - x - z + x*z) / (1 - x)	
\end{verbatimtab}
\end{code}
In this case the principal set is $\{Q\}$, the conditioning set is
$\{P\}$, and the marginal set is $\{R,A,B,C\}$.  Again the constraints
that the values of $x$, $y$, and $z$ must lie between $0$ and $1$ are
considered part of the result.

\subsection{Simple Table-Based Probability Inference}
\label{sec:symbolic}

Symbolic probability-network inference follows principles that were
already well-described by de~Moivre in the early 18th century (notably
before the famous paper from Bayes; in fact Bayes referred explicitly
to de~Moivre's work) \cite{demoivre,bayes}.  Modern
probability-network inference methods focus on efficiency through
clever factoring strategies and sophisticated graph manipulations
\cite{lauritzen,dambrosio}.  But for our purposes, inelegant
brute-force inference will suffice; this simple approach also helps to
elucidate the polynomial form of the probabilities that are computed.
Moreover this simple inference method offers its own routes to
improved performance: it happens that the necessary arithmetical
operations can be performed in parallel (taking advantage of computers
with multiple processors), and that the collected calculations are
essentially relational-database operations (for which well-optimized
software has been developed).

Simple table-based probability inference requires three steps:
joining, aggregation, and division.  In the first step the component
probability tables are \emph{joined} into the full-joint probability
table, by multiplying the component probabilities in a certain
fashion.  Each full-joint probability is the product of one element
from every component probability table; the overall calculation is
equivalent to a relational database join operation with a small amount
of post-processing.  In the second step the full-joint probabilities
are \emph{aggregated} to compute the marginal probabilities of the
queried events, by adding suitable joint probabilities together.  This
aggregation step corresponds to the relational database operation of
the same name.  For a conditional query, one aggregation must be
performed for the numerator events and a second aggregation for the
denominator events in the result table.  The numerator events use the
set union of primary variables in the query's principal and
conditioning sets; the denominator events use only the variables in
the conditioning set.  In the third step the aggregate probability of
each numerator event is \emph{divided} by the aggregate probability of
its corresponding denominator event; this is also a modified
relational join operation.  It emerges that each element of the result
table is the \emph{quotient} of \emph{sums} of \emph{products} of the
original component probabilities: in other words a fractional
polynomial in the model parameters, if each input component
probability was itself a polynomial.  When a query's conditioning set
is empty, no denominator table is constructed and no division is
performed; hence unconditional probability-table queries yield sums of
products of component probabilities (thus simple polynomial outputs
from simple polynomial inputs).

To illustrate symbolic probability-network inference, consider the
conditional probability query \condprob{Q}{P,R} for the model \mpq{}.
The principal set is $\{Q\}$ and the conditioning set is $\{P,R\}$.
In order to evaluate this query using simple table-based inference we
must first compute the full-joint probability table
\prob{P,Q,R,A,B,C}; then compute the numerator table \prob{P,R,Q}
(which uses the set union $\{Q\}\cup\{P,R\}$ of variables from the
principal and conditioning sets) and the denominator table \prob{P,R};
and finally divide the corresponding numerator and denominator
elements by one another.  The number of elements in the full-joint
table \prob{P,Q,R,A,B,C} is given by the product of the number of
states for each variable: in this case $2 \times 2 \times 2 \times 1
\times 4 \times 4$ which is $128$.  Every element in the full-joint
probability distribution is the product of several component
probabilities.  For example, the full-joint probability:
\begin{equation}
  \prob{P=\state{F},Q=\state{F},R=\state{T},
    A=\state{3},B=\state{1},C=\state{2}}
\end{equation}
is given by the product of the corresponding component probabilities:
\begin{equation}
  \begin{array}{l}
  \prob[0]{P=\state{F}} \;\times\;
  \condprob[0]{Q=\state{F}}{P=\state{F}} \;\times\;
  \condprob[0]{R=\state{T}}{P=\state{F},Q=\state{F}} \;\times\; \\
  \prob[0]{A=\state{3}} \;\times\;
  \condprob[0]{B=\state{1}}{P=\state{F},Q=\state{F},R=\state{T}} \;\times\;
  \condprob[0]{C=\state{2}}{A=\state{3},B=\state{1}}
  \end{array}
\end{equation}
Substituting the values from Table~\ref{tbl:m1-cpt}, this product
becomes
\begin{math}
  (1-x) (1-z) (1) (1) (1) (1)
\end{math}
which simplifies to $1-x-z+xz$.  For the \mpq{} model it happens that
only $4$ of the $128$ full-joint probabilities are not zero:
\begin{equation}
\label{eq:m1-joint}
\begin{tabular}[c]{r|llllll|l}\hline
\multicolumn{1}{l|}{\scshape {Index}} & 
\multicolumn{1}{l}{$P$} & 
\multicolumn{1}{l}{$Q$} & 
\multicolumn{1}{l}{$R$} & 
\multicolumn{1}{l}{$A$} & 
\multicolumn{1}{l}{$B$} & 
\multicolumn{1}{l|}{$C$} & 
\multicolumn{1}{l}{$\prob{{P, Q, R, A, B, C}}$} \\ \hline\hline
13 & 
$\state{T}$ & 
$\state{T}$ & 
$\state{T}$ & 
$\state{3}$ & 
$\state{3}$ & 
$\state{0}$ & 
$x y$ \\ \hline
55 & 
$\state{T}$ & 
$\state{F}$ & 
$\state{F}$ & 
$\state{3}$ & 
$\state{1}$ & 
$\state{2}$ & 
$x - x y$ \\ \hline
74 & 
$\state{F}$ & 
$\state{T}$ & 
$\state{T}$ & 
$\state{3}$ & 
$\state{2}$ & 
$\state{1}$ & 
$z - x z$ \\ \hline
103 & 
$\state{F}$ & 
$\state{F}$ & 
$\state{T}$ & 
$\state{3}$ & 
$\state{1}$ & 
$\state{2}$ & 
$1 - x - z + x z$ \\ \hline
\end{tabular}
\end{equation}
Here the index numbers are relative to all $128$ full-joint
probabilities, arranged in a particular lexicographic order.

The numerator table \prob{P,R,Q} and the denominator table \prob{P,R}
are computed by adding appropriate elements of this full-joint
probability table.  For example the marginal probability
\prob{P=\state{F},R=\state{T}} is given by the sum of the
probabilities of the corresponding nonzero elements of the full-joint
probability table (in rows 74 and 103):
\begin{equation}
  \prob{P=\state{F},Q=\state{T},R=\state{T},
    A=\state{3},B=\state{2},C=\state{1}}
  \;+\;
  \prob{P=\state{F},Q=\state{F},R=\state{T},
    A=\state{3},B=\state{1},C=\state{2}}
\end{equation}
Substituting the polynomial values from Equation~\ref{eq:m1-joint}
yields:
\begin{equation}
  \prob{P=\state{F},R=\state{T}} \quad \Rightarrow \quad
(z - x z) \;+\; 
(1 - x - z + x z)
\quad \Rightarrow \quad
1 - x
\end{equation}
The complete tables for \prob{P,R,Q} and \prob{P,R} are shown here:
\begin{equation}
\label{eq:m1-tpsqps}
\begin{tabular}[c]{r|lll|l}\hline
\multicolumn{1}{l|}{\scshape {Index}} & 
\multicolumn{1}{l}{$P$} & 
\multicolumn{1}{l}{$R$} & 
\multicolumn{1}{l|}{$Q$} & 
\multicolumn{1}{l}{$\prob{{P, R, Q}}$} \\ \hline\hline
1 & 
$\state{T}$ & 
$\state{T}$ & 
$\state{T}$ & 
$x y$ \\ \hline
2 & 
$\state{T}$ & 
$\state{T}$ & 
$\state{F}$ & 
$0$ \\ \hline
3 & 
$\state{T}$ & 
$\state{F}$ & 
$\state{T}$ & 
$0$ \\ \hline
4 & 
$\state{T}$ & 
$\state{F}$ & 
$\state{F}$ & 
$x - x y$ \\ \hline
5 & 
$\state{F}$ & 
$\state{T}$ & 
$\state{T}$ & 
$z - x z$ \\ \hline
6 & 
$\state{F}$ & 
$\state{T}$ & 
$\state{F}$ & 
$1 - x - z + x z$ \\ \hline
7 & 
$\state{F}$ & 
$\state{F}$ & 
$\state{T}$ & 
$0$ \\ \hline
8 & 
$\state{F}$ & 
$\state{F}$ & 
$\state{F}$ & 
$0$ \\ \hline
\end{tabular}
\qquad
\begin{tabular}[c]{r|ll|l}\hline
\multicolumn{1}{l|}{\scshape {Index}} & 
\multicolumn{1}{l}{$P$} & 
\multicolumn{1}{l|}{$R$} & 
\multicolumn{1}{l}{$\prob{{P, R}}$} \\ \hline\hline
1 & 
$\state{T}$ & 
$\state{T}$ & 
$x y$ \\ \hline
2 & 
$\state{T}$ & 
$\state{F}$ & 
$x - x y$ \\ \hline
3 & 
$\state{F}$ & 
$\state{T}$ & 
$1 - x$ \\ \hline
4 & 
$\state{F}$ & 
$\state{F}$ & 
$0$ \\ \hline
\end{tabular}
\end{equation}
Dividing each element of the numerator \prob{P,R,Q} result table by
the matching element of the denominator \prob{P,R} result table yields
the queried conditional probability table \condprob{Q}{P,R}:
\begin{equation}
\label{eq:m1-tqps}
\begin{tabular}[c]{r|lll|l}\hline
\multicolumn{1}{l|}{\scshape {Index}} & 
\multicolumn{1}{l}{$P$} & 
\multicolumn{1}{l}{$R$} & 
\multicolumn{1}{l|}{$Q$} & 
\multicolumn{1}{l}{$\condprob{{Q}}{{P, R}}$} \\ \hline\hline
1 & 
$\state{T}$ & 
$\state{T}$ & 
$\state{T}$ & 
$\left(x y\right) / \left(x y\right)$ \\ \hline
2 & 
$\state{T}$ & 
$\state{T}$ & 
$\state{F}$ & 
$\left(0\right) / \left(x y\right)$ \\ \hline
3 & 
$\state{T}$ & 
$\state{F}$ & 
$\state{T}$ & 
$\left(0\right) / \left(x - x y\right)$ \\ \hline
4 & 
$\state{T}$ & 
$\state{F}$ & 
$\state{F}$ & 
$\left(x - x y\right) / \left(x - x y\right)$ \\ \hline
5 & 
$\state{F}$ & 
$\state{T}$ & 
$\state{T}$ & 
$\left(z - x z\right) / \left(1 - x\right)$ \\ \hline
6 & 
$\state{F}$ & 
$\state{T}$ & 
$\state{F}$ & 
$\left(1 - x - z + x z\right) / \left(1 - x\right)$ \\ \hline
7 & 
$\state{F}$ & 
$\state{F}$ & 
$\state{T}$ & 
$\left(0\right) / \left(0\right)$ \\ \hline
8 & 
$\state{F}$ & 
$\state{F}$ & 
$\state{F}$ & 
$\left(0\right) / \left(0\right)$ \\ \hline
\end{tabular}
\end{equation}
Note that the condition \prob{P=\state{F},R=\state{F}} is impossible;
in terms of the embedded logical formulas the material implication
$R:=(P \rightarrow Q)$ cannot be false if its premise $P$ is false.
Hence both probabilities
\condprob{Q=\state{T}}{P=\state{F},R=\state{F}} and
\condprob{Q=\state{F}}{P=\state{F},R=\state{F}} conditioned on this
impossible event involve division by zero.  By default such
exceptional elements are not displayed; they are included above as the
quotient-expression $0/0$ to clarify the calculations that have
occurred.  Omitting these indeterminate elements and pivoting the
table to show the probabilities given each possible condition in the
same row, the same result for the query \condprob{Q}{P,R} is displayed
as:
\begin{equation}
\begin{tabular}[c]{r|ll|ll}\hline
\multicolumn{5}{l}{$\condprob{{Q}}{{P, R}}$} \\ \hline\hline
\multicolumn{1}{l|}{\scshape {Index}} & 
\multicolumn{1}{l}{$P$} & 
\multicolumn{1}{l|}{$R$} & 
\multicolumn{1}{l}{${Q=\state{T}}$} & 
\multicolumn{1}{l}{${Q=\state{F}}$} \\ \hline\hline
1, 2 & 
$\state{T}$ & 
$\state{T}$ & 
$\left(x y\right) / \left(x y\right)$ & 
$\left(0\right) / \left(x y\right)$ \\ \hline
3, 4 & 
$\state{T}$ & 
$\state{F}$ & 
$\left(0\right) / \left(x - x y\right)$ & 
$\left(x - x y\right) / \left(x - x y\right)$ \\ \hline
5, 6 & 
$\state{F}$ & 
$\state{T}$ & 
$\left(z - x z\right) / \left(1 - x\right)$ & 
$\left(1 - x - z + x z\right) / \left(1 - x\right)$ \\ \hline
\end{tabular}
\end{equation}

\subsection{Linear Functions Follow Form}  
\label{sec:linear}

If the primary variables $P_1,P_2,\ldots,P_v$ representing embedded
variables are modeled as a clique with a single parameter specifying
each probability in the joint distribution
$\prob[0]{P_1,P_2,\ldots,P_v}$, and furthermore if there are no other
parameters in the probability network, then all inferred probabilities
must be either linear functions of the parameters or quotients of such
linear functions.  This special case turns out to be quite useful, as
it is the natural way to model many problems combining logic and
probability.  For example, if we were to modify the \mpq{} model so
that the prior probabilities on $P$ and $Q$ were specified as this
joint distribution \prob[0]{P,Q} instead of as the separate
probabilities \prob[0]{P} and \condprob[0]{Q}{P}:
\begin{equation}
\begin{tabular}[c]{ll|l}\hline
\multicolumn{1}{l}{$P$} & 
\multicolumn{1}{l|}{$Q$} & 
\multicolumn{1}{l}{$\prob[0]{{P, Q}}$} \\ \hline\hline
$\state{T}$ & 
$\state{T}$ & 
$x_{1}$ \\ \hline
$\state{T}$ & 
$\state{F}$ & 
$x_{2}$ \\ \hline
$\state{F}$ & 
$\state{T}$ & 
$x_{3}$ \\ \hline
$\state{F}$ & 
$\state{F}$ & 
$x_{4}$ \\ \hline
\end{tabular}
\end{equation}
then all probabilities inferred from this revised model would be
linear functions of its parameters or quotients of such linear
functions.  For example:
\begin{equation}
\begin{tabular}[c]{lll|l}\hline
\multicolumn{1}{l}{$P$} & 
\multicolumn{1}{l}{$R$} & 
\multicolumn{1}{l|}{$Q$} & 
\multicolumn{1}{l}{$\prob{{P, R, Q}}$} \\ \hline\hline
$\state{T}$ & 
$\state{T}$ & 
$\state{T}$ & 
$x_{1}$ \\ \hline
$\state{T}$ & 
$\state{T}$ & 
$\state{F}$ & 
$0$ \\ \hline
$\state{T}$ & 
$\state{F}$ & 
$\state{T}$ & 
$0$ \\ \hline
$\state{T}$ & 
$\state{F}$ & 
$\state{F}$ & 
$x_{2}$ \\ \hline
$\state{F}$ & 
$\state{T}$ & 
$\state{T}$ & 
$x_{3}$ \\ \hline
$\state{F}$ & 
$\state{T}$ & 
$\state{F}$ & 
$x_{4}$ \\ \hline
$\state{F}$ & 
$\state{F}$ & 
$\state{T}$ & 
$0$ \\ \hline
$\state{F}$ & 
$\state{F}$ & 
$\state{F}$ & 
$0$ \\ \hline
\end{tabular}
\qquad
\begin{tabular}[c]{ll|l}\hline
\multicolumn{1}{l}{$P$} & 
\multicolumn{1}{l|}{$R$} & 
\multicolumn{1}{l}{$\prob{{P, R}}$} \\ \hline\hline
$\state{T}$ & 
$\state{T}$ & 
$x_{1}$ \\ \hline
$\state{T}$ & 
$\state{F}$ & 
$x_{2}$ \\ \hline
$\state{F}$ & 
$\state{T}$ & 
$x_{3} + x_{4}$ \\ \hline
$\state{F}$ & 
$\state{F}$ & 
$0$ \\ \hline
\end{tabular}
\qquad
\begin{tabular}[c]{ll|ll}\hline
\multicolumn{4}{l}{$\condprob{{Q}}{{P, R}}$} \\ \hline\hline
\multicolumn{1}{l}{$P$} & 
\multicolumn{1}{l|}{$R$} & 
\multicolumn{1}{l}{${Q=\state{T}}$} & 
\multicolumn{1}{l}{${Q=\state{F}}$} \\ \hline\hline
$\state{T}$ & 
$\state{T}$ & 
$\left(x_{1}\right) / \left(x_{1}\right)$ & 
$\left(0\right) / \left(x_{1}\right)$ \\ \hline
$\state{T}$ & 
$\state{F}$ & 
$\left(0\right) / \left(x_{2}\right)$ & 
$\left(x_{2}\right) / \left(x_{2}\right)$ \\ \hline
$\state{F}$ & 
$\state{T}$ & 
$\left(x_{3}\right) / \left(x_{3} + x_{4}\right)$ & 
$\left(x_{4}\right) / \left(x_{3} + x_{4}\right)$ \\ \hline
\end{tabular}
\end{equation}

\subsection{Handling Division by Zero}
\label{sec:divzero}

You may have noticed that quotients such as
\begin{math}
\left(x y\right) / \left(x y\right)
\end{math}
and
\begin{math}
\left(z - x z\right) / \left(1 - x\right)
\end{math}
were not simplified in the results displayed above.  That is because
it is important to recognize when division by zero is possible: this
corresponds to the imposition of an impossible condition, in which
case conditional probabilities are appropriately undefined.  (Having
established that there are no angels dancing on the point of a needle,
there is no unique answer to what proportion are boy-angels versus
girl-angels.)  Let us consider two additional options for handling
division by zero within the framework of parametric probability.

\paragraph{Double-Backslash Notation}

It is useful to have more compact notation to describe the value of a
quotient whose denominator might be zero.  For this the following
convention is proposed: let us say that the value of an expression
$\gamma \punless{\alpha=\beta}$ is usually $\gamma$, except that if
the condition $\alpha=\beta$ holds then the value of the expression is
undefined.  The double backslash $\punless$ may be read `unless' or
`except that it is undefined if'; you may think of it as a distant
relative of the set difference operator $\setminus$.  This new
construction is similar to the ternary conditional operator \verb|?:|
in the C programming language: in C the expression
\verb|a == b ? c : d| has the value \verb|c| if the condition
\verb|a == b| holds and \verb|d| otherwise \cite{c}.

Using this double-backslash notation the quotient $xy/x$ would be
rendered $y \punless{x=0}$, the quotient $x/x$ would be rendered $1
\punless{x=0}$, and the quotient $0/x$ would be rendered $0
\punless{x=0}$.  It is best to avoid double-backslash notation when
the undefining condition is tautological: thus $0/0$ should be
displayed as such (or some other designation for a value that is
\emph{always} undefined) instead of as the less intuitive $1
\punless{0=0}$.  By this double-backslash convention we can display
the result for the query \condprob{Q}{P,S} in the following way
(compare with the original table in Equation~\ref{eq:m1-tqps}):
\begin{equation}
\begin{tabular}[c]{r|lll|l}\hline
\multicolumn{1}{l|}{\scshape {Index}} & 
\multicolumn{1}{l}{$P$} & 
\multicolumn{1}{l}{$R$} & 
\multicolumn{1}{l|}{$Q$} & 
\multicolumn{1}{l}{$\condprob{{Q}}{{P, R}}$} \\ \hline\hline
1 & 
$\state{T}$ & 
$\state{T}$ & 
$\state{T}$ & 
$1 \punless{x y = 0}$ \\ \hline
2 & 
$\state{T}$ & 
$\state{T}$ & 
$\state{F}$ & 
$0 \punless{x y = 0}$ \\ \hline
3 & 
$\state{T}$ & 
$\state{F}$ & 
$\state{T}$ & 
$0 \punless{x y = x}$ \\ \hline
4 & 
$\state{T}$ & 
$\state{F}$ & 
$\state{F}$ & 
$1 \punless{x y = x}$ \\ \hline
5 & 
$\state{F}$ & 
$\state{T}$ & 
$\state{T}$ & 
$z \punless{x = 1}$ \\ \hline
6 & 
$\state{F}$ & 
$\state{T}$ & 
$\state{F}$ & 
$\left(1 - x - z + x z\right) / \left(1 - x\right)$ \\ \hline
7 & 
$\state{F}$ & 
$\state{F}$ & 
$\state{T}$ & 
$0/0$ \\ \hline
8 & 
$\state{F}$ & 
$\state{F}$ & 
$\state{F}$ & 
$0/0$ \\ \hline
\end{tabular}
\label{eq:m1-tqps-punless}
\end{equation}
Like most of the tables and figures in this document, this result
table was generated automatically by the author's computer program in
response to a structured query.  The current version of the program is
not very good at factoring polynomials: here it has not figured out
that the quotient
\begin{math}
\left(1 - x - z + x z\right) / \left(1 - x\right)
\end{math}
simplifies to $1-z \punless{x=1}$.  It would be good for a future
implementation of parametric probability analysis to be integrated
with a general-purpose computer algebra system; for temporal reasons
this was not done in the current version of the computer program.

\paragraph{Alternative Parametric Indeterminacy}

It may be desirable to handle division by zero in a different way, by
introducing additional parameters to encode indeterminate values while
preserving the semantics that probabilities are proportions that add
up to one.  For this example we might report the value of
$\condprob{Q=\state{T}}{P=\state{F},R=\state{F}}$ as $\theta$ and the
value of $\condprob{Q=\state{F}}{P=\state{F},R=\state{F}}$ as
$1-\theta$, where $\theta$ is a new parameter subject to the
constraint $0 \le \theta \le 1$ about which no other constraints are
allowed.  In this way we would maintain the property that these two
mutually exclusive and collectively exhaustive conditional
probabilities have the sum one and that the value of each probability
must lie between zero and one, but we would leave the precise value of
each probability indeterminate.  For example by this
parametric-indeterminacy convention we would consider the sum
\begin{math}
\left(x y\right) / \left(x y\right)
+
\left(0\right) / \left(x y\right)
\end{math}
to have the definite value $1$ even if it feasible that $xy=0$ (in
particular when using elements of the result table for the query
\condprob{{Q}}{{P, R}},
since in this context these fractional polynomial values would
describe the probabilities of mutually exclusive and collectively
exhaustive events).

\subsection{Boolean Polynomials and Coincident Probabilities}

\begin{table}
  \begin{tabular}{l@{\quad}l@{\quad}l@{\quad}l} \hline
    \scshape Logical & 
    \scshape Polynomial in $\mathbb{R}[P,Q]$ & 
    \scshape Polynomial in $\mathbb{F}_2[P,Q]$ & 
    \scshape Description \\ 
    \hline\hline
    $\state{T}$ & $1$ & $1$ & Elementary truth \\ \hline
    $\state{F}$ & $0$ & $0$ & Elementary falsity \\ \hline
    $\neg P$ & $1-P$ & $1+P$ & Negation (\textsc{not}) \\ \hline
    $P \wedge Q$ & $PQ$ & $PQ$ & Conjunction (\textsc{and}) \\ \hline
    $P \oplus Q$ & $P+Q-2PQ$ & $P+Q$ & 
    Exclusive disjunction (\textsc{xor}) \\ \hline
    $P \vee Q$ & $P+Q-PQ$ & $P+Q+PQ$ & 
    Inclusive disjunction (\textsc{or}) \\ \hline
    $P \rightarrow Q$ & $1-P+PQ$ & $1+P+PQ$ & Material implication \\ \hline
    $P \leftrightarrow Q$ & $1-P-Q+2PQ$ & $1+P+Q$ & 
    Biconditional (\textsc{xnor}) \\ \hline
    $P \uparrow Q$ & $1-PQ$ & $1+PQ$ & Nonconjunction (\textsc{nand})
    \\ \hline
    $P \downarrow Q$ & $1-P-Q+PQ$ & $1+P+Q+PQ$ & 
    Nondisjunction (\textsc{nor}) \\ \hline
  \end{tabular}
  \caption{Boolean representation of logical formulas, illustrated for
    propositional variables $P$ and $Q$ and for polynomials with real
    or finite-field coefficients.  Using real-number coefficients
    Boole's `special law' constraints $P^2=P$ and $Q^2=Q$ are
    necessary to limit the possible values of each variable to
    $\{0,1\}$.}
  \label{tbl:translation}
\end{table}

Let us briefly review Boole's polynomial notation for logical
formulas.  The mappings between Boole's polynomials and what is now
standard logical notation (our mash-up from Hilbert, Peano, and
others) are shown in Table~\ref{tbl:translation}.  These are sometimes
called the `Stone isomorphisms' after \cite{stone}.  Despite the
convenient phrase `Boolean translation' it should be noted that Boole
did not use polynomials to translate from some other conventional
system of symbolic notation for logic---for in his time there was no
such convention. (Frege's \emph{Begriffsschrift} was published in
1879, many years after Boole's death in 1864; likewise Boole's
lifetime predated the works of Peano and Hilbert in which much of
modern logical notation was developed \cite{cajori}.)  Nonetheless,
when viewed as translation, what Boole described turns out to be a
special case of Lagrange polynomial interpolation
\cite{carnielli-polynomizing}.

Contrary to a very common misrepresentation, Boole's polynomials used
ordinary integer coefficients for which $1+1=2$.  However there are
some advantages to using instead coefficients in the finite field
$\mathbb{F}_2$ of order $2$, which uses integer arithmetic modulo $2$
(hence $1+1=0$, addition and subtraction are the same operation, and
each value is its own additive inverse).  With coefficients in
$\mathbb{F}_2$ the polynomials that represent logical formulas are
simpler in form and the number of distinct polynomials is finite
(given a finite set of propositional variables).
Table~\ref{tbl:translation} includes mappings from conventional
logical formulas to polynomials with coefficients in the binary finite
field as well as to polynomials with coefficients in the real numbers.

Boole's polynomial notation for logical formulas is often understood
in a monolithic way but it was actually the expression of two
different ideas: first, the idea that classical logical operations
(conjunction, disjunction, negation, and so on) are equivalent to
certain combinations of ordinary arithmetical operations, when
elementary truth values are represented as ordinary numbers; and
second, the idea that the probability that a logical formula is true
is a polynomial function of the probabilities that its constituent
propositional variables are true.  In a special case the polynomial
that denotes a logical formula coincides with the polynomial that
expresses the probability that the formula itself is true.  It is
important to recognize the independence property required for this
coincidence, and to generalize a means to compute appropriate
probabilities in the case that this independence property does not
hold (Boole did both).

\paragraph{Coincidence from Independent Propositional Variables}

When the propositional variables in use are modeled as
probabilistically independent of one another (as described in
Section~\ref{sec:prior}), then the probability that any compound
formula is true coincides with its Boolean polynomial representation
(using real coefficients).  For example, as shown in
Table~\ref{tbl:translation}, the Boolean representation of the logical
formula $X \leftrightarrow Y$ is the polynomial $1-X-Y+2XY$.  This
Boolean coincidence principle provides that if $X$ and $Y$ are
independent in the probability model, with $\prob[0]{X=\state{T}}:=x$
and $\prob[0]{Y=\state{T}}:=y$, then the probability $\prob{(X
  \leftrightarrow Y)=\state{T}}$ that the compound formula $X
\leftrightarrow Y$ is true has the value $1-x-y+2xy$ that mirrors the
Boolean polynomial representation $1-X-Y+2XY$ of this compound
formula.  This coincidence is not generally present when the
propositional variables in use are probabilistically correlated.  For
example if $\prob[0]{X=\state{T}}:=x$,
$\condprob[0]{Y=\state{T}}{X=\state{T}}:=y$, and
$\condprob[0]{Y=\state{T}}{X=\state{F}}:=z$, then the probability
$\prob{(X \leftrightarrow Y)=\state{T}}$ has the value $1-x-z+xy+xz$.
This value is different from $1-x-y+2xy$ exactly when $y \neq z$, in
other words when the probabilities of $X$ and $Y$ are nontrivially
correlated.  The general symbolic probability-network inference method
discussed in this section computes correct answers with or without the
assertion that the propositional variables are independent.

\section{Secondary Analysis: Algebra, Optimization, and Search}

The results generated by symbolic probability-network inference are
always algebraic functions of the model parameters; more specifically
they must be polynomials or quotients of polynomials, when the domain
of each primary variable is finite and when each input component
probability is itself a polynomial.  These computed polynomials are
ordinary mathematical objects that can be manipulated by ordinary
mathematical methods.  Boole focused on \emph{optimization} as the
secondary analysis that followed his symbolic probability inference.
In addition to this very useful technique, we can broaden the scope of
secondary analysis to include more general applications of symbolic
and numerical analysis to the polynomials generated by symbolic
probability-network inference.  Here we consider algebra,
optimization, and search.

\subsection{Algebra with Polynomials}
\label{sec:algebra}

Perhaps the simplest kind of secondary analysis for the fractional
polynomials computed by symbolic probability-network inference is
elementary algebra: these formulas can be added, subtracted,
multiplied, and divided to form new fractional polynomials.  For
example, consider the difference between the probability that the
logical formula $P \rightarrow Q$ holds and the conditional
probability that $Q$ is true given that $P$ is true.  Since the model
\mpq{} uses the definition $R:=(P \rightarrow Q)$, the requisite
difference is
$\prob{R=\state{T}}-\condprob{Q=\state{T}}{P=\state{T}}$.  To compute
this we first select the appropriate elements from the result tables
presented in Section~\ref{sec:queries}:
\begin{code}
\begin{alltt}
pqlsh> \emph{\$tr item 1}
\end{alltt}
\begin{verbatimtab}
1 - x + x*y
\end{verbatimtab}
\begin{alltt}
pqlsh> \emph{\$tqp item 1}
\end{alltt}
\begin{verbatimtab}
(x*y) / (x)
\end{verbatimtab}
\end{code}
We then ask \texttt{pqlsh} to compute the difference (as a human or
computer algebra system could easily do):
\begin{code}
\begin{alltt}
pqlsh> \emph{pql::expr "[\$tr item 1] - ([\$tqp item 1])"}
\end{alltt}
\begin{verbatimtab}
(x - x^2 - x*y + x^2*y) / (x)
\end{verbatimtab}
\end{code}
Rewriting this difference with \TeX{} formatting and the
double-backslash notation from Section~\ref{sec:divzero} we have:
\begin{eqnarray}
\label{eq:sqp}
\prob{R=\state{T}}-\condprob{Q=\state{T}}{P=\state{T}}
& \Rightarrow &
1 - x - y + x y \punless{x = 0}
\end{eqnarray}
It follows from Equation~\ref{eq:sqp} that the two quantities
$\prob{R=\state{T}}$ and $\condprob{Q=\state{T}}{P=\state{T}}$ are
different unless $x=1$ or $y=1$ or both (with the caveat that the
difference is undefined when $x=0$).  Perhaps this is clearer when the
polynomial difference in Equation~\ref{eq:sqp} is factored as
$(1-x)(1-y)(x/x)$, which plainly has roots $x=1$ and $y=1$.  In terms
of the example, the probability that the statement `bird implies
flight' is true differs from the conditional probability of flight
given bird, unless the creature is certainly a bird ($x=1$) or it is
certain that all birds can fly ($y=1$) or both.  However the
difference is undefined if the creature could not possibly be a bird
($x=0$), because in that case the requested condition $P=\state{T}$
that the creature is a bird is impossible.

As another example of secondary analysis by algebra, let us calculate
the expectation $\mathbb{E}(B)$: the mean value of the number $B$ of
\state{True} propositions among $\{P,Q,R\}$.  For this we begin with
the result for the probability-table query $\prob{B}$, computed as
described in Section~\ref{sec:symbolic}:
\begin{equation}
\begin{tabular}[c]{r|l|l}\hline
\multicolumn{1}{l|}{\scshape {Index}} & 
\multicolumn{1}{l|}{$B$} & 
\multicolumn{1}{l}{$\prob{{B}}$} \\ \hline\hline
1 & 
$\state{0}$ & 
$0$ \\ \hline
2 & 
$\state{1}$ & 
$1 - z - x y + x z$ \\ \hline
3 & 
$\state{2}$ & 
$z - x z$ \\ \hline
4 & 
$\state{3}$ & 
$x y$ \\ \hline
\end{tabular}
\end{equation}
The expected value $\mathbb{E}(B)$ is then computed in the usual way
as the sum of each possible value $b_i \in \{0,1,2,3\}$ of $B$
weighted by the probability that it is attained:
\begin{math}
  \sum_{i} \left[ b_i \times \prob{B=b_i} \right]
\end{math}.  
Thus we calculate:
\begin{eqnarray}
\label{eq:m1-b}
0 \times (
0
  ) + 1 \times (
1 - z - x y + x z
  ) + 2 \times (
z - x z
  ) + 3 \times (
x y
  )
& \Rightarrow &
1 + z + 2 x y - x z
\end{eqnarray}
which establishes using simple algebra that the expected value
$\mathbb{E}(B)$ is the polynomial
\begin{math}
1 + z + 2 x y - x z
\end{math}.

\subsection{Polynomial Optimization}
\label{sec:optimization}

Next we consider secondary analysis using optimization (mathematical
programming).  The polynomials computed by symbolic
probability-network inference, and the additional polynomials derived
by algebraic calculation, can be used as constraints and objective
functions in optimization problems.  In \emph{Laws of Thought}, Boole
routinely sought to calculate the minimum and maximum feasible values
of some polynomial objective computed by symbolic probability-network
inference, subject to constraints reflecting the laws of probability.
However, Boole's 19th-century optimization methods were not very
robust.  Polynomial optimization problems are still challenging to
solve; they are generally nonlinear and nonconvex, and therefore they
can have local solutions that are not global solutions.  But there are
much better global optimization algorithms now, including
reformulation-linearization and semidefinite programming techniques
\cite{sherali,lasserre}.  Specifically for the application of
parametric probability analysis, the author has developed a new
reformulation and linearization algorithm that computes interval
bounds on the global solutions to multivariate polynomial optimization
problems, using mixed integer-linear program approximations; the user
controls the tightness of the bounds and the time required to compute
them by setting the number of reformulation variables
\cite{norman-nlp}.

To illustrate optimization as secondary analysis, let us consider the
minimum and maximum feasible values of the expectation $\mathbb{E}(B)$
given in Equation~\ref{eq:m1-b}, with the added constraint that there
a $75\%$ chance or less that the statement that `bird implies flight'
is true: $\prob{R=\state{T}} \le 0.75$.  Recalling the specification
of the \mpq{} probability model given in Tables \ref{tbl:m1-variables}
and \ref{tbl:m1-cpt}, we include the constraints $0 \le x \le 1$, $0
\le y \le 1$, and $0 \le z \le 1$.  For the additional constraint
regarding $\prob{R=\state{T}}$ we recall from
Section~\ref{sec:queries} that \prob{R=\state{T}} evaluates to
\begin{math}
1 - x + x y
\end{math}.
The resulting inequality
\begin{math}
1 - x + x y
\le 0.75
\end{math}
simplifies to
\begin{math}
0.25 + x y \le x
\end{math}.
Hence to find minimum and maximum bounds on $\mathbb{E}(B)$ we must
solve the paired polynomial optimization problems:
\begin{equation}
\begin{array}{r@{\quad}l}
\mbox{minimize}: & 1 + z + 2 x y - x z \\
\mbox{subject to}:
& 0.25 + x y \le x \\
\mbox{and}:
& 0 \le x \le 1 \\
& 0 \le y \le 1 \\
& 0 \le z \le 1
\end{array}
\qquad
\begin{array}{r@{\quad}l}
\mbox{maximize}: & 1 + z + 2 x y - x z \\
\mbox{subject to}:
& 0.25 + x y \le x \\
\mbox{and}:
& 0 \le x \le 1 \\
& 0 \le y \le 1 \\
& 0 \le z \le 1
\end{array}
\end{equation}
Using a small number of reformulation variables, the author's bounded
global polynomial optimization solver computes that the global minimum
lies in the interval
\begin{math}
[
0.938
,
1.000
]
\end{math}
and that the global maximum lies in the interval
\begin{math}
[
2.500
,
2.594
]
\end{math}.
With more reformulation variables the solver generates tighter
intervals
\begin{math}
[
1.000
,
1.000
]
\end{math}
and
\begin{math}
[
2.500
,
2.500
]
\end{math}.
The solver reports that the global minimum
\begin{math}
1.000
\end{math}
is achieved at the point
\begin{math}
(
x = 0.984, y = 0.000, z = 0.000
)
\end{math}
and that the global maximum
\begin{math}
2.500
\end{math}
is achieved at the point
\begin{math}
(
x = 1.000, y = 0.750, z = 0.266
)
\end{math}.
You may notice by inspection of Equation~\ref{eq:m1-b} that, absent
any constraints besides the bounds on the parameters $x$, $y$, and
$z$, the global minimum of the expected value $\mathbb{E}(B)$ is
exactly $1$ (when $x=y=z=0$) and its global maximum is exactly $3$
(when $x=y=z=1$).  The constraint $\prob{R=\state{T}} \le 0.75$ added
for this example has rendered some of that range infeasible.

\subsection{General Search}
\label{sec:search}

Providing yet another mode of secondary analysis, polynomials
generated by symbolic probability-network inference can be used in
general computer-science search problems that might be awkward to
formulate in terms of algebraic equations or numeric optimization.  To
illustrate, note that in the model \mpq{} it is impossible for all
three propositions $P$, $Q$, and $R$ to be false; as
Equation~\ref{eq:m1-tpsqps} shows, the joint probability
\prob{P=\state{F},R=\state{F},Q=\state{F}} evaluates to zero.  This
makes semantic sense based on the definition $R:=(P \rightarrow Q)$,
for if the premise $P$ is false then the material implication $P
\rightarrow Q$ must be true; hence both propositions cannot be false
simultaneously.  Let us search for a logical proposition that has a
different property: a logical function $R^*$ of $P$ and $Q$ such that
the number of true propositions among $\{P,Q,R^*\}$ must be odd (either
$1$ or $3$).

To set up this search let us replace the original component
probability table \condprob{R}{P,Q} shown in Table~\ref{tbl:m1-cpt}
part~(c) with the table for \condprob{R^*}{P,Q} given in
Equation~\ref{eq:sp} in which the probabilities of each value of $R^*$
given each combination of values for $P$ and $Q$ are encoded by the
parameters $t_1$, $t_2$, $t_3$, and $t_4$ as described in
Section~\ref{sec:unknown}.  Using this replacement parametric table
and symbolic probability-network inference as in
Section~\ref{sec:symbolic} we compute the probability distribution on
the number $B^*$ of true propositions among $\{P,Q,R^*\}$:
\begin{equation}
\label{eq:m1s-b}
\begin{tabular}[c]{l|l}\hline
\multicolumn{1}{l|}{$B^*$} & 
\multicolumn{1}{l}{$\prob{{B^*}}$} \\ \hline\hline
$\state{0}$ & 
$1 - x - z - t_{4} + x z + x t_{4} + z t_{4} - x z t_{4}$ \\ \hline
$\state{1}$ & 
$x + z + t_{4} - x y - x z - x t_{2} - x t_{4} - z t_{3} - z t_{4} + x y t_{2} + x z t_{3} + x z t_{4}$ \\ \hline
$\state{2}$ & 
$x y + x t_{2} + z t_{3} - x y t_{1} - x y t_{2} - x z t_{3}$ \\ \hline
$\state{3}$ & 
$x y t_{1}$ \\ \hline
\end{tabular}
\end{equation}
Now we desire to find the values of $(t_1,t_2,t_3,t_4)$ for which the
only possible values of $B^*$ are $1$ and $3$.  In other words, we
require that the polynomials \prob{B^*=\state{0}} and
\prob{B^*=\state{2}} are both identically zero after substituting the
selected values of $(t_1,t_2,t_3,t_4)$.  Considering each value $t_i
\in \{0,1\}$ there are $2^4 \Rightarrow 16$ possible values of the
vector $(t_1,t_2,t_3,t_4)$.  For this small problem we simply
enumerate every possibility and substitute these values into the
polynomials in the result table for \prob{B^*} given in
Equation~\ref{eq:m1s-b}. (Of course the point of most search
algorithms is to avoid exhaustive enumeration of the search space; we
eschew such luxury for now.)  In Table~\ref{tbl:util-b} the table on
the left gives the instantiations of \prob{B^*=\state{0}} and the table
on the right gives the instantiations of \prob{B^*=\state{2}} at all
$16$ possible values of $(t_1,t_2,t_3,t_4)$.

\begin{table}
\begin{tabular}[c]{r|llll|l}\hline
\multicolumn{1}{l|}{\scshape {Index}} & 
\multicolumn{1}{l}{$t_{1}$} & 
\multicolumn{1}{l}{$t_{2}$} & 
\multicolumn{1}{l}{$t_{3}$} & 
\multicolumn{1}{l|}{$t_{4}$} & 
\multicolumn{1}{l}{$\condutil[0]{{\prob{B^*=\state{0}}}}{{t_{1}, t_{2}, t_{3}, t_{4}}}$} \\ \hline\hline
1 & 
$\state{0}$ & 
$\state{0}$ & 
$\state{0}$ & 
$\state{0}$ & 
$1 - x - z + x z$ \\ \hline
2 & 
$\state{0}$ & 
$\state{0}$ & 
$\state{0}$ & 
$\state{1}$ & 
$0$ \\ \hline
3 & 
$\state{0}$ & 
$\state{0}$ & 
$\state{1}$ & 
$\state{0}$ & 
$1 - x - z + x z$ \\ \hline
4 & 
$\state{0}$ & 
$\state{0}$ & 
$\state{1}$ & 
$\state{1}$ & 
$0$ \\ \hline
5 & 
$\state{0}$ & 
$\state{1}$ & 
$\state{0}$ & 
$\state{0}$ & 
$1 - x - z + x z$ \\ \hline
6 & 
$\state{0}$ & 
$\state{1}$ & 
$\state{0}$ & 
$\state{1}$ & 
$0$ \\ \hline
7 & 
$\state{0}$ & 
$\state{1}$ & 
$\state{1}$ & 
$\state{0}$ & 
$1 - x - z + x z$ \\ \hline
8 & 
$\state{0}$ & 
$\state{1}$ & 
$\state{1}$ & 
$\state{1}$ & 
$0$ \\ \hline
9 & 
$\state{1}$ & 
$\state{0}$ & 
$\state{0}$ & 
$\state{0}$ & 
$1 - x - z + x z$ \\ \hline
10 & 
$\state{1}$ & 
$\state{0}$ & 
$\state{0}$ & 
$\state{1}$ & 
$0$ \\ \hline
11 & 
$\state{1}$ & 
$\state{0}$ & 
$\state{1}$ & 
$\state{0}$ & 
$1 - x - z + x z$ \\ \hline
12 & 
$\state{1}$ & 
$\state{0}$ & 
$\state{1}$ & 
$\state{1}$ & 
$0$ \\ \hline
13 & 
$\state{1}$ & 
$\state{1}$ & 
$\state{0}$ & 
$\state{0}$ & 
$1 - x - z + x z$ \\ \hline
14 & 
$\state{1}$ & 
$\state{1}$ & 
$\state{0}$ & 
$\state{1}$ & 
$0$ \\ \hline
15 & 
$\state{1}$ & 
$\state{1}$ & 
$\state{1}$ & 
$\state{0}$ & 
$1 - x - z + x z$ \\ \hline
16 & 
$\state{1}$ & 
$\state{1}$ & 
$\state{1}$ & 
$\state{1}$ & 
$0$ \\ \hline
\end{tabular}
\qquad
\begin{tabular}[c]{r|llll|l}\hline
\multicolumn{1}{l|}{\scshape {Index}} & 
\multicolumn{1}{l}{$t_{1}$} & 
\multicolumn{1}{l}{$t_{2}$} & 
\multicolumn{1}{l}{$t_{3}$} & 
\multicolumn{1}{l|}{$t_{4}$} & 
\multicolumn{1}{l}{$\condutil[0]{{\prob{B^*=\state{2}}}}{{t_{1}, t_{2}, t_{3}, t_{4}}}$} \\ \hline\hline
1 & 
$\state{0}$ & 
$\state{0}$ & 
$\state{0}$ & 
$\state{0}$ & 
$x y$ \\ \hline
2 & 
$\state{0}$ & 
$\state{0}$ & 
$\state{0}$ & 
$\state{1}$ & 
$x y$ \\ \hline
3 & 
$\state{0}$ & 
$\state{0}$ & 
$\state{1}$ & 
$\state{0}$ & 
$z + x y - x z$ \\ \hline
4 & 
$\state{0}$ & 
$\state{0}$ & 
$\state{1}$ & 
$\state{1}$ & 
$z + x y - x z$ \\ \hline
5 & 
$\state{0}$ & 
$\state{1}$ & 
$\state{0}$ & 
$\state{0}$ & 
$x$ \\ \hline
6 & 
$\state{0}$ & 
$\state{1}$ & 
$\state{0}$ & 
$\state{1}$ & 
$x$ \\ \hline
7 & 
$\state{0}$ & 
$\state{1}$ & 
$\state{1}$ & 
$\state{0}$ & 
$x + z - x z$ \\ \hline
8 & 
$\state{0}$ & 
$\state{1}$ & 
$\state{1}$ & 
$\state{1}$ & 
$x + z - x z$ \\ \hline
9 & 
$\state{1}$ & 
$\state{0}$ & 
$\state{0}$ & 
$\state{0}$ & 
$0$ \\ \hline
10 & 
$\state{1}$ & 
$\state{0}$ & 
$\state{0}$ & 
$\state{1}$ & 
$0$ \\ \hline
11 & 
$\state{1}$ & 
$\state{0}$ & 
$\state{1}$ & 
$\state{0}$ & 
$z - x z$ \\ \hline
12 & 
$\state{1}$ & 
$\state{0}$ & 
$\state{1}$ & 
$\state{1}$ & 
$z - x z$ \\ \hline
13 & 
$\state{1}$ & 
$\state{1}$ & 
$\state{0}$ & 
$\state{0}$ & 
$x - x y$ \\ \hline
14 & 
$\state{1}$ & 
$\state{1}$ & 
$\state{0}$ & 
$\state{1}$ & 
$x - x y$ \\ \hline
15 & 
$\state{1}$ & 
$\state{1}$ & 
$\state{1}$ & 
$\state{0}$ & 
$x + z - x y - x z$ \\ \hline
16 & 
$\state{1}$ & 
$\state{1}$ & 
$\state{1}$ & 
$\state{1}$ & 
$x + z - x y - x z$ \\ \hline
\end{tabular}
\caption{Instantiations of \prob{B^*=\state{0}} and \prob{B^*=\state{2}}
  calculated by substituting the listed values of $(t_1,t_2,t_3,t_4)$
  into the corresponding polynomials from Equation~\ref{eq:m1s-b}.  In
  row~10 both instantiated polynomials are identically zero.}
\label{tbl:util-b}
\end{table}

Comparing these tables you can see that only in row~10 are both
polynomials identically zero.  In every other case it is either
possible that \prob{B^*=\state{0}} is not zero, that
\prob{B^*=\state{2}} is not zero, or that both probabilities are not
zero.\footnote{For this problem it happens that every polynomial in
  Table~\ref{tbl:util-b} that is not identically zero has a value that
  is strictly greater than zero at some feasible value of the
  parameters $(x,y,z)$.  In general it would be necessary to solve an
  optimization problem to verify that a polynomial which is
  symbolically different from zero indeed has a feasible value greater
  than zero, taking into account all of the constraints provided.}
Row~10 corresponds to the vector $(t_1,t_2,t_3,t_4)=(1,0,0,1)$ and in
turn to the logical formula $R^*_{1001}:=(P \leftrightarrow Q)$ that
combines $P$ and $Q$ using the biconditional operator (see
Section~\ref{sec:unknown}).  Here is the probability distribution on
the number $B^*_{1001}$ of true propositions among the set
$\{P,Q,R^*_{1001}\}$, which is obtained from substituting the solution
values $(t_1,t_2,t_3,t_4)=(1,0,0,1)$ into Equation~\ref{eq:m1s-b}:
\begin{equation}
\begin{tabular}[c]{l|l}\hline
\multicolumn{1}{l|}{$B^*_{1001}$} & 
\multicolumn{1}{l}{$\prob{{B^*_{1001}}}$} \\ \hline\hline
$\state{0}$ & 
$0$ \\ \hline
$\state{1}$ & 
$1 - x y$ \\ \hline
$\state{2}$ & 
$0$ \\ \hline
$\state{3}$ & 
$x y$ \\ \hline
\end{tabular}
\end{equation}
Thus by applying general search as secondary analysis, we have
computed that the number of true propositions among the set $\{P,Q,P
\leftrightarrow Q\}$ must be odd; it cannot happen that exactly $0$ or
$2$ of these propositions are true, regardless of the prior
probabilities $x$, $y$, and $z$.  Moreover, search has demonstrated
that the only other formulas $R^*$ with this odd-number property must
have the same truth table as $P \leftrightarrow Q$.  For a formula
with any other truth table there would be feasible values of the
parameters $(x,y,z)$ for which $\prob{B^*=\state{0}}>0$,
$\prob{B^*=\state{2}}>0$, or both.

\section{Additional Modeling Issues}

Having discussed primary and secondary analysis of parametric
probability networks, there are two additional modeling issues to
consider.  First, there are two different techniques to model
conditions during parametric probability analysis: as denominator
events in conditional probability queries and as constraints in
optimization problems.  Second, it is possible to encode some
statements about implication directly as conditional probabilities,
without the intermediate device of embedded formulas from the
propositional calculus.  Using such direct probability encoding,
constraints on conditional probabilities can express quantification
without the need for classical logical quantifiers and with the option
to specify more precise fractional values than just `some'.

\subsection{Subjunctive Conditions, Imperative Constraints}
\label{sec:subjunctive}

Problems in logic and probability commonly involve \emph{conditions},
and there are two idioms for representing conditions during parametric
probability analysis.  Using the \emph{subjunctive} idiom, a condition
is modeled as the denominator event in a conditional probability-table
query.  Using the \emph{imperative} idiom, a condition is modeled as
an equality constraint in an optimization problem.  These alternative
formulations have slightly different semantics.  In the subjunctive
formulation we ask hypothetically what would be the probability of
some event, if the stated condition were to hold.  In the imperative
formulation we assert factually that the stated condition must hold,
and then ask what is the probability of some event under this
necessary condition.  These two idioms yield essentially the same
solutions, although they report the exception that the stated
condition is impossible in two different ways: in the subjunctive
formulation an impossible condition causes division by zero, but in
the imperative formulation an impossible condition produces an
unsatisfiable system of equations.  Both conditioning idioms can be
used for parametric probability networks with or without embedded
logical formulas.

\paragraph{Example: Two Modes of Modus Ponens}

There are two different ways to express the familiar phenomenon of
\emph{modus ponens} with a parametric probability network, and these
correspond to two slightly different questions.  First we might ask in
a \emph{subjunctive} mood: if $P$ and $P \rightarrow Q$ were to be
true, what would be the probability that $Q$ is also true?  Second we
might specify in an \emph{indicative} mood that $P$ and $P \rightarrow
Q$ must certainly be true, and then ask the probability that $Q$ is
also true.  Recall that in model \mpq{} the variable $R$ stands for
the proposition $P \rightarrow Q$.

For the subjunctive formulation we use a probability-table query to
compute \condprob{Q}{P,R}, the result of which is shown in
Equation~\ref{eq:m1-tqps}.  The first element of this result table
gives the desired conditional probability:
\begin{eqnarray}
\condprob{Q=\state{T}}{P=\state{T},R=\state{T}} & \Rightarrow &
\left(x y\right) / \left(x y\right)
\end{eqnarray}
The value of this quotient is one unless its denominator is zero (that
is, $1 \punless{xy=0}$ using the notation of
Section~\ref{sec:divzero}).  We interpret this to mean that that if
both propositions $P$ and $P \rightarrow Q$ were to be true, then $Q$
would also be true with probability $1$ (unless $x=0$ or $y=0$ or
both, in which exceptional cases the condition
$\{P=\state{T},R=\state{T}\}$ would be impossible and thus the
requested conditional probability would be the indeterminate
expression $0/0$).  In other words, if the creature happens to be a
bird, and if it happens to be true that bird implies flight, then it
would certainly be true that the creature can fly.  Except that if it
were already known \emph{a priori} that it is impossible for the
creature to be a bird, and/or that there is no chance that a bird can
fly, then the question of flight assuming the stated conditions would
not have a definite answer because the conditions would be impossible.

Alternatively, to use the indicative formulation, we build an
optimization query in which we ask the minimum feasible value of the
objective \prob{Q=\state{T}} subject to the following constraints that
$P$ and $R$ must certainly be true:
\begin{eqnarray}
\prob{P=\state{T}} & = & 1 \\
\prob{R=\state{T}} & = & 1
\end{eqnarray}
We specify the constraint on $P$ using the first element of its result
table (the output \prob{P} happens to be the same as the input
\prob[0]{P} shown in Table~\ref{tbl:m1-cpt}) and the constraint on $R$
from the first element of \prob{R} shown in Section~\ref{sec:queries}:
\begin{code}
\begin{alltt}
pqlsh> \emph{set cp "[[[\$m1 table P] infer] item 1] == 1"}
\end{alltt}
\begin{verbatimtab}
x == 1
\end{verbatimtab}
\begin{alltt}
pqlsh> \emph{set cr "[\$tr item 1] == 1"}
\end{alltt}
\begin{verbatimtab}
1 - x + x*y == 1
\end{verbatimtab}
\end{code}
The second constraint simplifies to
\begin{math}
x = x y
\end{math}.
Moving on, we obtain the objective function as the first element of
the result table for the query \prob{Q}:
\begin{code}
\begin{alltt}
pqlsh> \emph{set tq [[\$m1 table Q] infer]; \$tq print;}
\end{alltt}
\begin{verbatimtab}
Q	| Pr( {Q} )	
-------	-------
T	| z + x*y - x*z	
F	| 1 - z - x*y + x*z	
\end{verbatimtab}
\end{code}
Our optimization query is a request for a polynomial program using
this objective function and these constraints:
\begin{code}
\begin{alltt}
pqlsh> \emph{set pq [\$m1 pprog -min "[\$tq item 1]" \$cp \$cr]; return;}
\end{alltt}
\end{code}
This optimization query generates the following polynomial
optimization problem:
\begin{equation}
\label{eq:m1-qmin}
\begin{array}{r@{\quad}l}
\mbox{minimize}: & z + x y - x z \\
\mbox{subject to}:
& x = 1 \\
& x = x y \\
\mbox{and}:
& 0 \le x \le 1 \\
& 0 \le y \le 1 \\
& 0 \le z \le 1
\end{array}
\end{equation}
Solving this problem gives bounds on the minimum feasible value of
\prob{Q=\state{T}} under the constraints $\prob{P=\state{T}}=1$ and
$\prob{R=\state{T}}=1$:
\begin{code}
\begin{alltt}
pqlsh> \emph{\$pq solve; \$pq solution}
\end{alltt}
\begin{verbatimtab}
1.000 1.000
\end{verbatimtab}
\begin{alltt}
pqlsh> \emph{\$pq point}
\end{alltt}
\begin{verbatimtab}
{x = 1.000} {y = 1.000} {z = 0.000}
\end{verbatimtab}
\end{code}
This solution to the optimization problem in Equation~\ref{eq:m1-qmin}
demonstrates that proposition $Q$ must be true (that is, the variable
$Q$ attains the state \state{T} with minimum probability in the
interval
\begin{math}
[
1.000
,
1.000
]
\end{math}) 
if $P$ and $P \rightarrow Q$ are constrained to be true.  (You can see
by inspection that the global solution to Equation~\ref{eq:m1-qmin} is
exactly $1$.)  In other words, if it is certainly true that the
creature is a bird, and it is certainly true that bird implies flight,
then it is certainly true that the creature can fly.  The exceptional
cases $x=0$ and $y=0$ are handled differently in this formulation:
instead of causing division by zero, they identify infeasible points
relative to the constraints in Equation~\ref{eq:m1-qmin}.



\subsection{Direct Probability Encoding}
\label{sec:direct}

It is possible to model implication and quantification directly in
parametric probability networks, without using the classical logical
devices of material implication or universal and existential
quantifiers.  We have already seen the statement of material
implication $P \rightarrow Q$ used to model the idea that $Q$ is a
necessary consequence of $P$ (with subjunctive and indicative idioms
to impose the condition that this statement of material implication is
true).  We could instead constrain the conditional probability
$\condprob{Q=\state{T}}{P=\state{T}}$ to $1$ to express the idea that
that $Q$ is a necessary consequence of the premise $P$.  Similarly the
constraint $\condprob{Q=\state{T}}{P=\state{T}}=0$ is an alternative
way to model the assertion that $Q$ is never a consequence of $P$.

These equality constraints on conditional probabilities are
alternatives to universally-quantified statements such as $\forall
\alpha (P(\alpha) \rightarrow Q(\alpha))$ to say that all $P$ are $Q$;
or $\forall \alpha (P(\alpha) \rightarrow \neg Q(\alpha))$ to say that
no $P$ are $Q$.  In a similar fashion, to model the assertion that $Q$
sometimes follows $P$ we could constrain the relevant conditional
probability to be strictly greater than zero:
$\condprob{Q=\state{T}}{P=\state{T}} > 0$.  And to model the assertion
that $Q$ sometimes does not follow $P$ we could constrain the relevant
conditional probability to be strictly less than one:
$\condprob{Q=\state{T}}{P=\state{T}} < 1$.  These inequality
constraints on conditional probabilities are alternatives to
existentially-quantified formulas such as $\exists \alpha (P(\alpha)
\rightarrow Q(\alpha))$ or $\exists \alpha (P(\alpha) \wedge
Q(\alpha))$ to say that some $P$ are $Q$; or $\exists \alpha
(P(\alpha) \wedge \neg Q(\alpha))$ to say that some $P$ are not $Q$.
In the framework of parametric probability, quantified variables like
$\alpha$ are distinct both from primary variables and from parameters.

Recall from the secondary analysis in Section~\ref{sec:algebra} that
the conditional probability that $Q$ is true given that $P$ is true is
mathematically distinct from the probability that the
material-implication statement $P \rightarrow Q$ is true.  These
correspond to symbolically different polynomials; the arithmetical
difference between them depends on the prior probabilities assigned to
the events $P$ and $Q$ as reported in Equation~\ref{eq:sqp}.

\paragraph{The Option of Existential Import}

The constraints introduced to quantify propositions can be specified
such that they do or do not have existential import, as the user
desires.  In general, constraining \emph{input} component
probabilities specified by the user can have different effects from
constraining \emph{output} probabilities computed by the system;
existential import is one of those mutable effects.  Elementary
algebra helps to clarify the consequences of various polynomial
constraints.

To illustrate, here are the input component probability table
\condprob[0]{Q}{P} from Table~\ref{tbl:m1-cpt} part~(b); the computed
table \condprob{Q}{P} for this same conditional probability, created
as in Section~\ref{sec:symbolic}; and the computed probability table
\prob{P}:
\begin{equation}
\begin{tabular}[c]{l|ll}\hline
\multicolumn{3}{l}{$\condprob[0]{{Q}}{{P}}$} \\ \hline\hline
\multicolumn{1}{l|}{$P$} & 
\multicolumn{1}{l}{${Q=\state{T}}$} & 
\multicolumn{1}{l}{${Q=\state{F}}$} \\ \hline\hline
$\state{T}$ & 
$y$ & 
$1 - y$ \\ \hline
$\state{F}$ & 
$z$ & 
$1 - z$ \\ \hline
\end{tabular}
\qquad
\begin{tabular}[c]{l|ll}\hline
\multicolumn{3}{l}{$\condprob{{Q}}{{P}}$} \\ \hline\hline
\multicolumn{1}{l|}{$P$} & 
\multicolumn{1}{l}{${Q=\state{T}}$} & 
\multicolumn{1}{l}{${Q=\state{F}}$} \\ \hline\hline
$\state{T}$ & 
$\left(x y\right) / \left(x\right)$ & 
$\left(x - x y\right) / \left(x\right)$ \\ \hline
$\state{F}$ & 
$\left(z - x z\right) / \left(1 - x\right)$ & 
$\left(1 - x - z + x z\right) / \left(1 - x\right)$ \\ \hline
\end{tabular}
\qquad
\begin{tabular}[c]{l|l}\hline
\multicolumn{1}{l|}{$P$} & 
\multicolumn{1}{l}{$\prob{{P}}$} \\ \hline\hline
$\state{T}$ & 
$x$ \\ \hline
$\state{F}$ & 
$1 - x$ \\ \hline
\end{tabular}
\label{eq:qp-cpt}
\label{eq:qp-infer}
\end{equation}
Now, adding the constraint $y>0$ that the \emph{input} value
\condprob[0]{Q=\state{T}}{P=\state{T}} must be strictly greater than
zero would specify that $Q$ is sometimes a consequence of $P$ (when
$P$ happens to be true) without asserting that $P$ is ever true.  For
if $x=0$ then the computed probability \prob{P=\state{T}} would be
zero even if $y>0$.  Thus the input-value constraint
$\condprob[0]{Q=\state{T}}{P=\state{T}}>0$, meaning $y>0$, does not
affect the values in \prob{P}; it has no existential import.
Incidentally, in the case $x=0$ the inferred conditional probability
\condprob{Q=\state{T}}{P=\state{T}}, calculated by the laws of
probability as the quotient
$\prob{P=\state{T},Q=\state{T}}\;/\;\prob{P=\state{T}}$, would be
indeterminate due to division by zero---regardless of the value $y$
assigned to the corresponding input conditional probability
$\condprob[0]{Q=\state{T}}{P=\state{T}}$.  For conditional
probabilities, it is essential to distinguish between input and output
values; what is computed by the system may differ from what was input
by the user (particularly regarding denominators that may be zero).

On the other hand, adding the constraint
\begin{math}
\left(x y\right) / \left(x\right)
> 0
\end{math}
that the \emph{output} value \condprob{Q=\state{T}}{P=\state{T}} is
strictly greater than zero would indeed carry existential import.
Since both parameters $x$ and $y$ are constrained to be nonnegative,
by simple algebra this constraint would require the product $xy$ in
the numerator to be strictly greater than zero, which would in turn
require both $x>0$ and $y>0$.  Consequently the output value
$\prob{P=\state{T}}$, which has the polynomial value $x$, would also
be required to be strictly greater than zero.  So the output-value
constraint $\condprob{Q=\state{T}}{P=\state{T}}>0$, meaning
\begin{math}
\left(x y\right) / \left(x\right)
> 0
\end{math}, would also assert $\prob{P=\state{T}}>0$; it has
existential import.

It is always possible to constrain output probabilities; but which
input probabilities are available to constrain depends on the graph
structure of the probability network.  For example the presence of an
input table \condprob[0]{Q}{P} requires that $P$ is a parent of $Q$ in
the network graph.  If instead $P$ and $Q$ were joined a clique then
their joint probability distribution \prob[0]{P,Q} would be specified
as the input component probability table; in this case there would be
no separate table \condprob[0]{Q}{P} of input values to be
constrained.

As an implementation detail many optimization solvers to not
distinguish between strict and weak inequality constraints.  To work
around this limitation we can use some small positive constant
$\epsilon$ and weak inequality constraints $p \ge \epsilon$ and $p \le
1-\epsilon$ to approximate the strict inequalities $p>0$ and $p<1$ for
a polynomial of interest $p$.  For example using $\epsilon=0.1$ the
constraint $y \ge 0.1$ approximates the strict inequality
$\condprob[0]{Q=\state{T}}{P=\state{T}}>0$.  This approximation scheme
allows common optimization solvers to compute results that reliably
distinguish between probabilities that are exactly zero, those that
are exactly one, and those that have some intermediate value.

\paragraph{Fractional Quantification}

Using probability directly to model quantification offers an important
benefit: we are not limited to the classical quantifiers `all' and
`some'.  Instead it is possible to describe and to constrain the
precise proportion of cases for which some logical formula is true or
false.  In other words, besides the constraints $p=0$, $p=1$, $p>0$,
and $p<1$ that the probability $p$ encoding some statement of
quantification is equal to zero, equal to one, strictly greater than
zero, or strictly less than one, it is possible to specify arbitrary
polynomial constraints on $p$.  Thus in addition to statements like
`all $P$ are $Q$' or `some $P$ are not $Q$' we can model such
assertions as `exactly $c$ percent of $P$ are $Q$' or `between $a$ and
$b$ percent of $P$ are $Q$' or `if there are any $P$, then twice as
many $P$ are $Q$ as $P$ are $R$'.  In certain cases, the requisite
constraints are guaranteed to be linear in the model parameters.

As a philosophical aside, \emph{probability} is best understood as the
proportion of some underlying basic measure.  There is diversity in
what that basic measure can represent: number or cardinality (in which
case the proportion is frequency); the absolute weight of subjective
belief or of causal propensity (in which case the proportion is
subjective probability); monetary value; mass; or some other property.
It is essential for the property chosen as a basic measure to be
additive across set unions of measured events, which is the
quintessential mathematical property of a measure.  In the course of
contemplating proportional statements of quantification more precise
than `all' and `some', and for reasoning with Nilsson-style fractional
truth values, it may be worthwhile to clarify what the basic measure
is intended to mean.

\section{Analysis of Selected Problems}
\label{sec:examples}

Now let us apply parametric probability analysis to an assortment of
problems from the literature, each of which has been advertised as
being difficult or impossible to solve by formal mathematical methods.
We shall see that several well-known problems---about card games with
logical rules, axioms with uncertainty, counterfactual conditions, and
truthful knights and lying knaves---are nothing more than parametric
probability problems.  Such problems are easily solved by parametric
probability analysis; many of them turn out to be linear optimization
problems whose constraints and objectives are the solutions to
probability queries.

\subsection{Johnson-Laird's Winning Hand}

Let us begin with a problem from Johnson-Laird that was also discussed
by Bringsjord \cite{johnson-laird-illusory,bringsjord}:
\begin{quote}
  If one of the following assertions is true then so is the other:
  \begin{enumerate}
  \item[1.] There is a king in the hand if and only if there is an ace in
    the hand.
  \item[2.] There is a king in the hand.
  \end{enumerate}
  Which is more likely to be in the hand, if either: the king or the
  ace?  Prove that you are correct.
\end{quote}
It is worth repeating the challenge that Bringsjord issued with this
problem:
\begin{quote}
  I don't even think Bayesian systems can possibly solve logic
  problems that involve probability.  \ldots{} I would very much like to
  see a Bayesian system take this declarative information as input,
  and yield the correct answer, and a proof that this \emph{is} the
  answer.  I assure you that I will not hold my breath.
\end{quote}
Though the computational method presented here is more appropriately
called `Boolean' than `Bayesian', there is no difficulty in solving
logic problems that involve probability using parametric probability
analysis.  Proof, such as it is, is supplied by elementary algebra.
Let us consider three ways to analyze this ace-king problem: using
Boolean polynomials to simplify the logical formula involved; using
parametric probability analysis with subjunctive conditioning; and
using parametric probability analysis with indicative conditioning.

\paragraph{Polynomial Simplification of Logical Formulas}

Perhaps the easiest way to solve this ace-king problem is to simplify
the logical assertion in it.  Using $A$ to represent the proposition
that there is an ace in the hand and $K$ to represent the proposition
that there is a king in the hand, the problem specifies the assertion
\begin{math}
  (K \leftrightarrow A) \leftrightarrow K
\end{math}.  This compound logical formula simplifies to the atomic
formula $A$ after Boolean polynomial representation using the rules in
Table~\ref{tbl:translation}.  For example using coefficients in the
binary finite field $\mathbb{F}_2$ the inner biconditional translates
to the polynomial $1+K+A$.  Substituting this value, the entire
formula $(1+K+A) \leftrightarrow K$ translates to the polynomial
$1+(1+K+A)+K$.  Recall that using integer arithmetic modulo $2$ either
elementary value $0$ or $1$ is its own additive inverse; thus all
terms in this polynomial $1+1+K+A+K$ cancel out except $A$.  Using
real-number coefficients would generate the same answer (keeping in
mind that $A$ can be substituted for $A^2$ and $K$ for $K^2$ due to
Boole's `special law' constraints $A^2=A$ and $K^2=K$).  Taking
advantage of such polynomial representation and simplification, an
equivalent problem statement would be:
\begin{quote}
  There is an ace in the hand. \\
  Which is more likely to be in the hand, if either: the king or the ace?
\end{quote}
It is evident that the ace must be equally likely or more likely than
the king, since the ace is present with certainty.

\paragraph{Parametric Probability Model}

Next let us construct an explicit parametric probability network for
this ace-king problem, using the technique described in
Section~\ref{sec:embedding}.  We copy the propositional variables $A$
and $K$ into the probability network as primary variables with the set
$\{\state{T},\state{F}\}$ of possible values representing elementary
truth and falsity.  We add a third primary variable $P$ which is
defined as the value of the compound logical formula $(K
\leftrightarrow A) \leftrightarrow K$ asserted in the problem
description.  We introduce parameters $x_1$ through $x_4$ to specify
prior probabilities.  Here is the network graph:
\begin{equation}
  \includegraphics{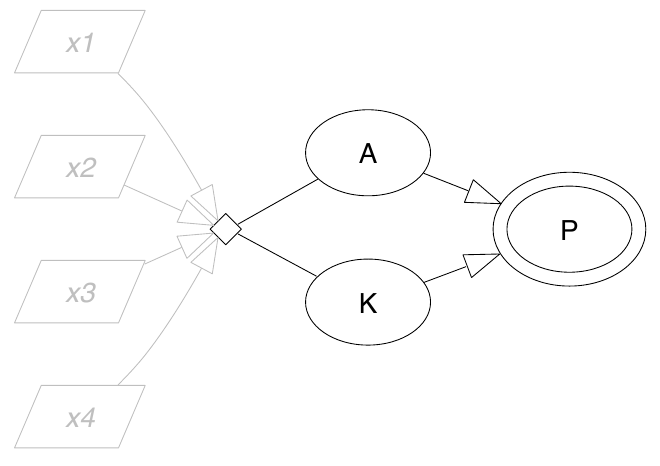}
\end{equation}
As in Section~\ref{sec:prior} we specify an uninformative prior
probability distribution \prob[0]{A,K} on the variables $A$ and $K$
using the $x_i$ parameters, with the constraints $0 \le x_i \le 1$ and
$x_1+x_2+x_3+x_4=1$ to enforce the laws of probability.  Following the
method of Section~\ref{sec:conditional} we construct a conditional
probability table \condprob[0]{P}{A,K} to express the definition
$P:=((K \leftrightarrow A) \leftrightarrow K)$.  These two component
probability tables complete the ace-king model:
\begin{equation}
\begin{tabular}[c]{ll|l}\hline
\multicolumn{1}{l}{$A$} & 
\multicolumn{1}{l|}{$K$} & 
\multicolumn{1}{l}{$\prob[0]{{A, K}}$} \\ \hline\hline
$\state{T}$ & 
$\state{T}$ & 
$x_{1}$ \\ \hline
$\state{T}$ & 
$\state{F}$ & 
$x_{2}$ \\ \hline
$\state{F}$ & 
$\state{T}$ & 
$x_{3}$ \\ \hline
$\state{F}$ & 
$\state{F}$ & 
$x_{4}$ \\ \hline
\end{tabular}
\label{eq:tacea}
\qquad
\begin{tabular}[c]{ll|ll}\hline
\multicolumn{4}{l}{$\condprob[0]{{P}}{{A, K}}$} \\ \hline\hline
\multicolumn{1}{l}{$A$} & 
\multicolumn{1}{l|}{$K$} & 
\multicolumn{1}{l}{${P=\state{T}}$} & 
\multicolumn{1}{l}{${P=\state{F}}$} \\ \hline\hline
$\state{T}$ & 
$\state{T}$ & 
$1$ & 
$0$ \\ \hline
$\state{T}$ & 
$\state{F}$ & 
$1$ & 
$0$ \\ \hline
$\state{F}$ & 
$\state{T}$ & 
$0$ & 
$1$ \\ \hline
$\state{F}$ & 
$\state{F}$ & 
$0$ & 
$1$ \\ \hline
\end{tabular}
\label{eq:tacep}
\end{equation}
Using this parametric probability network we can compare the
probabilities that $A$ and $K$ are true, given the condition that $P$
is true---using both subjunctive and imperative formulations to
express the condition.

\paragraph{Analysis in the Subjunctive Mood}

Using the subjunctive idiom for conditioning, we desire the difference
between the conditional probabilities
$\condprob{A=\state{T}}{P=\state{T}}$ and
$\condprob{K=\state{T}}{P=\state{T}}$.  The relevant fractional
polynomial values are included in the result tables for the queries
\begin{math}
\condprob{{A}}{{P}}
\end{math}
and
\begin{math}
\condprob{{K}}{{P}}
\end{math}, computed as described in Section~\ref{sec:symbolic}:
\begin{equation}
\begin{tabular}[c]{l|ll}\hline
\multicolumn{3}{l}{$\condprob{{A}}{{P}}$} \\ \hline\hline
\multicolumn{1}{l|}{$P$} & 
\multicolumn{1}{l}{${A=\state{T}}$} & 
\multicolumn{1}{l}{${A=\state{F}}$} \\ \hline\hline
$\state{T}$ & 
$\left(x_{1} + x_{2}\right) / \left(x_{1} + x_{2}\right)$ & 
$\left(0\right) / \left(x_{1} + x_{2}\right)$ \\ \hline
$\state{F}$ & 
$\left(0\right) / \left(x_{3} + x_{4}\right)$ & 
$\left(x_{3} + x_{4}\right) / \left(x_{3} + x_{4}\right)$ \\ \hline
\end{tabular}
\qquad
\begin{tabular}[c]{l|ll}\hline
\multicolumn{3}{l}{$\condprob{{K}}{{P}}$} \\ \hline\hline
\multicolumn{1}{l|}{$P$} & 
\multicolumn{1}{l}{${K=\state{T}}$} & 
\multicolumn{1}{l}{${K=\state{F}}$} \\ \hline\hline
$\state{T}$ & 
$\left(x_{1}\right) / \left(x_{1} + x_{2}\right)$ & 
$\left(x_{2}\right) / \left(x_{1} + x_{2}\right)$ \\ \hline
$\state{F}$ & 
$\left(x_{3}\right) / \left(x_{3} + x_{4}\right)$ & 
$\left(x_{4}\right) / \left(x_{3} + x_{4}\right)$ \\ \hline
\end{tabular}
\end{equation}
The desired difference uses the first element of each table, combined
by elementary algebra:
\begin{eqnarray}
\left(x_{1} + x_{2}\right) / \left(x_{1} + x_{2}\right)
-
\left(x_{1}\right) / \left(x_{1} + x_{2}\right)
& \Rightarrow &
\left(x_{2}\right) / \left(x_{1} + x_{2}\right)
\end{eqnarray}
To answer the question posed in the problem we must determine the
minimum and maximum feasible values of this difference, subject to the
constraints on the parameters involved.  These extreme values are the
solutions to the following pair of optimization problems, which share
a set of linear constraints and a fractional linear objective function:
\begin{equation}
\begin{array}{r@{\quad}l}
\mbox{minimize}: & \left(x_{2}\right) / \left(x_{1} + x_{2}\right) \\
\mbox{subject to}:
& x_{1} + x_{2} + x_{3} + x_{4} = 1 \\
\mbox{and}:
& 0 \le x_{1} \le 1 \\
& 0 \le x_{2} \le 1 \\
& 0 \le x_{3} \le 1 \\
& 0 \le x_{4} \le 1
\end{array}
\qquad
\begin{array}{r@{\quad}l}
\mbox{maximize}: & \left(x_{2}\right) / \left(x_{1} + x_{2}\right) \\
\mbox{subject to}:
& x_{1} + x_{2} + x_{3} + x_{4} = 1 \\
\mbox{and}:
& 0 \le x_{1} \le 1 \\
& 0 \le x_{2} \le 1 \\
& 0 \le x_{3} \le 1 \\
& 0 \le x_{4} \le 1
\end{array}
\end{equation}
One way to solve these fractional linear programs is through the
Charnes-Cooper transformation, which converts them into ordinary
linear programs \cite{charnes}.  After such reformulation, standard
linear optimization finds the minimum value
\begin{math}
0
\end{math}
which is achieved at the point
\begin{math}
(
x_{1} = 1, x_{2} = 0, x_{3} = 0, x_{4} = 0
)
\end{math}
and the maximum value
\begin{math}
1
\end{math}
which is achieved at the point
\begin{math}
(
x_{1} = 0, x_{2} = 1, x_{3} = 0, x_{4} = 0
)
\end{math}.  Thus the difference
\begin{math}
  \condprob{A=\state{T}}{P=\state{T}} -
  \condprob{K=\state{T}}{P=\state{T}}
\end{math}
is bounded by zero and one.  This implies the inequality:
\begin{eqnarray}
  \condprob{A=\state{T}}{P=\state{T}} & \ge &
  \condprob{K=\state{T}}{P=\state{T}}
\end{eqnarray}
In other words, given the condition $P:=((K \leftrightarrow A)
\leftrightarrow K)$ stated in the problem description, it is at least
as likely that there is an ace in the hand as a king.

\paragraph{Parametric Probability in the Indicative Mood}

We can find the equivalent solution using the same parametric
probability network model, but with an imperative rather than
subjunctive query formulation.  In this idiom our objective function
is the difference $\prob{A=\state{T}}-\prob{K=\state{T}}$ between the
unconditional probabilities that there is an ace versus a king in the
hand.  The relevant probabilities, computed as in
Section~\ref{sec:symbolic}, appear in the results for the
probability-table queries \prob{A} and \prob{K}:
\begin{equation}
\begin{tabular}[c]{l|l}\hline
\multicolumn{1}{l|}{$A$} & 
\multicolumn{1}{l}{$\prob{{A}}$} \\ \hline\hline
$\state{T}$ & 
$x_{1} + x_{2}$ \\ \hline
$\state{F}$ & 
$x_{3} + x_{4}$ \\ \hline
\end{tabular}
\qquad
\begin{tabular}[c]{l|l}\hline
\multicolumn{1}{l|}{$K$} & 
\multicolumn{1}{l}{$\prob{{K}}$} \\ \hline\hline
$\state{T}$ & 
$x_{1} + x_{3}$ \\ \hline
$\state{F}$ & 
$x_{2} + x_{4}$ \\ \hline
\end{tabular}
\end{equation}
The requisite difference
\begin{math}
(
x_{1} + x_{2}
) - (
x_{1} + x_{3}
) 
\end{math}
simplifies to
\begin{math}
x_{2} - x_{3}
\end{math}.
Continuing on, in this indicative formulation the condition in the
problem statement is now modeled as the additional constraint
$\prob{P=\state{T}}=1$, using this result for the query \prob{P}:
\begin{equation}
\begin{tabular}[c]{l|l}\hline
\multicolumn{1}{l|}{$P$} & 
\multicolumn{1}{l}{$\prob{{P}}$} \\ \hline\hline
$\state{T}$ & 
$x_{1} + x_{2}$ \\ \hline
$\state{F}$ & 
$x_{3} + x_{4}$ \\ \hline
\end{tabular}
\end{equation}
Therefore we formulate the following pair of optimization problems to
find the minimum and maximum feasible values of the difference between
the probability of the ace and the king, subject to the constraints
that the assertion in the problem description is true and that the
laws of probability are followed:
\begin{equation}
\begin{array}{r@{\quad}l}
\mbox{minimize}: & x_{2} - x_{3} \\
\mbox{subject to}:
& x_{1} + x_{2} + x_{3} + x_{4} = 1 \\
& x_{1} + x_{2} = 1 \\
\mbox{and}:
& 0 \le x_{1} \le 1 \\
& 0 \le x_{2} \le 1 \\
& 0 \le x_{3} \le 1 \\
& 0 \le x_{4} \le 1
\end{array}
\qquad
\begin{array}{r@{\quad}l}
\mbox{maximize}: & x_{2} - x_{3} \\
\mbox{subject to}:
& x_{1} + x_{2} + x_{3} + x_{4} = 1 \\
& x_{1} + x_{2} = 1 \\
\mbox{and}:
& 0 \le x_{1} \le 1 \\
& 0 \le x_{2} \le 1 \\
& 0 \le x_{3} \le 1 \\
& 0 \le x_{4} \le 1
\end{array}
\end{equation}
These are simple linear programs whose objective functions are not
fractional.  Solving them with standard methods yields the minimum
value
\begin{math}
0
\end{math}
which is achieved at the point
\begin{math}
(
x_{1} = 1, x_{2} = 0, x_{3} = 0, x_{4} = 0
)
\end{math}
and the maximum value
\begin{math}
1
\end{math}
which is achieved at the point
\begin{math}
(
x_{1} = 0, x_{2} = 1, x_{3} = 0, x_{4} = 0
)
\end{math}.
In other words, subject to the constraint $\prob{P=\state{T}}=1$ that
the formula $(K \leftrightarrow A) \leftrightarrow K$ asserted in the
problem description is true, the difference
$\prob{A=\state{T}}-\prob{K=\state{T}}$ between the probability of the
ace and the king is bounded by zero and one.  It follows that
$\prob{A=\state{T}} \ge \prob{K=\state{T}}$ and therefore that the ace
is equally likely or more likely than the king, subject to the
constraint $\prob{P=\state{T}}=1$ that the assertion in the problem
statement is true.

\subsection{The Monster from Paris}
\label{sec:paris}

Next we examine a problem from Paris, Mui\~no, and Rosefield
concerning inconsistent propositions \cite{paris-muino}.  We shall see
that the parameterized consequence relation $^\eta
\triangleright_\zeta$ and the related concepts of maximal consistency
and primary and secondary probability thresholds presented by these
authors are simply indirect ways of describing linear optimization
problems and their solutions.  For this example there is a
hypothetical creature about which there are three propositions: $P$
that it is a chicken killer; $Q$ that it is Japanese; and $R$ that it
is a salamander.  There are also three compound formulas used as
axiom-like assertions (here designated $S_1$, $S_2$ and $S_3$):
\begin{equation}
  S_1 := (P \wedge Q) \qquad
  S_2 := (\neg(Q \wedge R) \wedge P) \qquad
  S_3 := (R \wedge ( \neg P \rightarrow (R \wedge Q)))
  \label{eq:s3}
\end{equation}
Additionally, there are several more compound formulas which are used
as queries (here designated $S_4$ through $S_8$):
\begin{equation}
  S_4 := (P \wedge R) \qquad
  S_5 := (P \wedge (Q \vee R)) \qquad
  S_6 := R \qquad
  S_7 := \neg R \qquad
  S_8 := (R \wedge \neg R)
  \label{eq:s8}
\end{equation}
The probability-network graph, which also includes the parameters $x_1$
through $x_8$, is shown here:
\begin{equation}
  \includegraphics{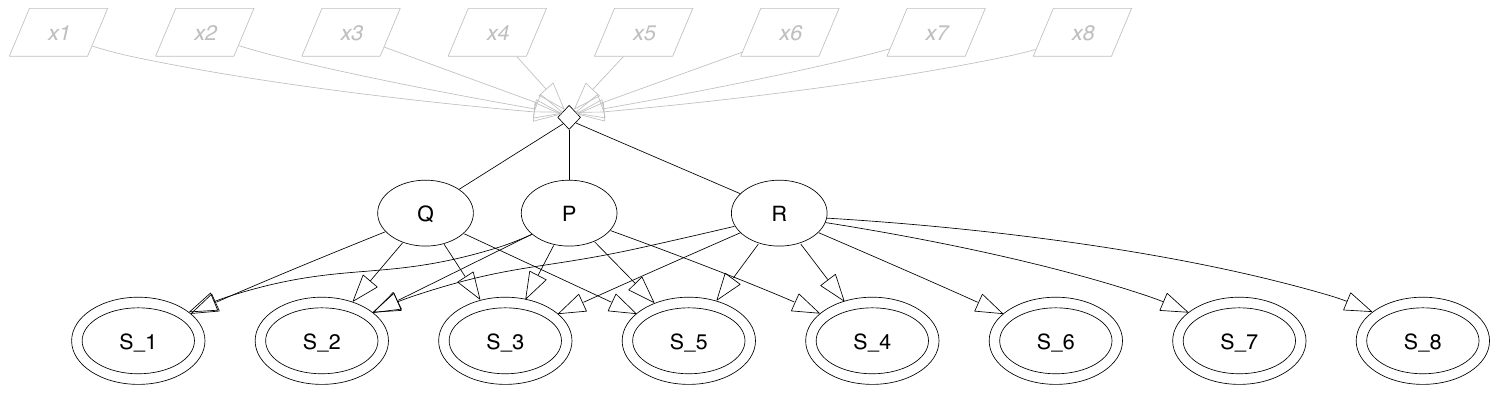}
\end{equation}
For this parametric probability network we specify an uninformative
probability distribution \prob[0]{P,Q,R} on the primary variables $P$,
$Q$, and $R$ copied from the propositional variables, using the $x_i$
parameters as in Section~\ref{sec:prior}:
\begin{equation}
\begin{tabular}[c]{lll|l}\hline
\multicolumn{1}{l}{$P$} & 
\multicolumn{1}{l}{$Q$} & 
\multicolumn{1}{l|}{$R$} & 
\multicolumn{1}{l}{$\prob[0]{{P, Q, R}}$} \\ \hline\hline
$\state{T}$ & 
$\state{T}$ & 
$\state{T}$ & 
$x_{1}$ \\ \hline
$\state{T}$ & 
$\state{T}$ & 
$\state{F}$ & 
$x_{2}$ \\ \hline
$\state{T}$ & 
$\state{F}$ & 
$\state{T}$ & 
$x_{3}$ \\ \hline
$\state{T}$ & 
$\state{F}$ & 
$\state{F}$ & 
$x_{4}$ \\ \hline
$\state{F}$ & 
$\state{T}$ & 
$\state{T}$ & 
$x_{5}$ \\ \hline
$\state{F}$ & 
$\state{T}$ & 
$\state{F}$ & 
$x_{6}$ \\ \hline
$\state{F}$ & 
$\state{F}$ & 
$\state{T}$ & 
$x_{7}$ \\ \hline
$\state{F}$ & 
$\state{F}$ & 
$\state{F}$ & 
$x_{8}$ \\ \hline
\end{tabular}
\end{equation}
Each parameter $x_i$ is subject to the constraint $0 \le x_i \le 1$
and they are collectively constrained by $\sum_i x_i = 1$.  Following
the method of Section~\ref{sec:conditional} we construct a conditional
probability table for each compound formula:
\begin{equation}
\begin{tabular}[c]{ll|ll}\hline
\multicolumn{4}{l}{$\condprob[0]{{S_1}}{{P, Q}}$} \\ \hline\hline
\multicolumn{1}{l}{$P$} & 
\multicolumn{1}{l|}{$Q$} & 
\multicolumn{1}{l}{${S_1=\state{T}}$} & 
\multicolumn{1}{l}{${S_1=\state{F}}$} \\ \hline\hline
$\state{T}$ & 
$\state{T}$ & 
$1$ & 
$0$ \\ \hline
$\state{T}$ & 
$\state{F}$ & 
$0$ & 
$1$ \\ \hline
$\state{F}$ & 
$\state{T}$ & 
$0$ & 
$1$ \\ \hline
$\state{F}$ & 
$\state{F}$ & 
$0$ & 
$1$ \\ \hline
\end{tabular}
\qquad
\begin{tabular}[c]{lll|ll}\hline
\multicolumn{5}{l}{$\condprob[0]{{S_2}}{{P, Q, R}}$} \\ \hline\hline
\multicolumn{1}{l}{$P$} & 
\multicolumn{1}{l}{$Q$} & 
\multicolumn{1}{l|}{$R$} & 
\multicolumn{1}{l}{${S_2=\state{T}}$} & 
\multicolumn{1}{l}{${S_2=\state{F}}$} \\ \hline\hline
$\state{T}$ & 
$\state{T}$ & 
$\state{T}$ & 
$0$ & 
$1$ \\ \hline
$\state{T}$ & 
$\state{T}$ & 
$\state{F}$ & 
$1$ & 
$0$ \\ \hline
$\state{T}$ & 
$\state{F}$ & 
$\state{T}$ & 
$1$ & 
$0$ \\ \hline
$\state{T}$ & 
$\state{F}$ & 
$\state{F}$ & 
$1$ & 
$0$ \\ \hline
$\state{F}$ & 
$\state{T}$ & 
$\state{T}$ & 
$0$ & 
$1$ \\ \hline
$\state{F}$ & 
$\state{T}$ & 
$\state{F}$ & 
$0$ & 
$1$ \\ \hline
$\state{F}$ & 
$\state{F}$ & 
$\state{T}$ & 
$0$ & 
$1$ \\ \hline
$\state{F}$ & 
$\state{F}$ & 
$\state{F}$ & 
$0$ & 
$1$ \\ \hline
\end{tabular}
\qquad
\begin{tabular}[c]{lll|ll}\hline
\multicolumn{5}{l}{$\condprob[0]{{S_3}}{{P, Q, R}}$} \\ \hline\hline
\multicolumn{1}{l}{$P$} & 
\multicolumn{1}{l}{$Q$} & 
\multicolumn{1}{l|}{$R$} & 
\multicolumn{1}{l}{${S_3=\state{T}}$} & 
\multicolumn{1}{l}{${S_3=\state{F}}$} \\ \hline\hline
$\state{T}$ & 
$\state{T}$ & 
$\state{T}$ & 
$1$ & 
$0$ \\ \hline
$\state{T}$ & 
$\state{T}$ & 
$\state{F}$ & 
$0$ & 
$1$ \\ \hline
$\state{T}$ & 
$\state{F}$ & 
$\state{T}$ & 
$1$ & 
$0$ \\ \hline
$\state{T}$ & 
$\state{F}$ & 
$\state{F}$ & 
$0$ & 
$1$ \\ \hline
$\state{F}$ & 
$\state{T}$ & 
$\state{T}$ & 
$1$ & 
$0$ \\ \hline
$\state{F}$ & 
$\state{T}$ & 
$\state{F}$ & 
$0$ & 
$1$ \\ \hline
$\state{F}$ & 
$\state{F}$ & 
$\state{T}$ & 
$0$ & 
$1$ \\ \hline
$\state{F}$ & 
$\state{F}$ & 
$\state{F}$ & 
$0$ & 
$1$ \\ \hline
\end{tabular}
\end{equation}
\begin{equation}
\begin{tabular}[c]{ll|ll}\hline
\multicolumn{4}{l}{$\condprob[0]{{S_4}}{{P, R}}$} \\ \hline\hline
\multicolumn{1}{l}{$P$} & 
\multicolumn{1}{l|}{$R$} & 
\multicolumn{1}{l}{${S_4=\state{T}}$} & 
\multicolumn{1}{l}{${S_4=\state{F}}$} \\ \hline\hline
$\state{T}$ & 
$\state{T}$ & 
$1$ & 
$0$ \\ \hline
$\state{T}$ & 
$\state{F}$ & 
$0$ & 
$1$ \\ \hline
$\state{F}$ & 
$\state{T}$ & 
$0$ & 
$1$ \\ \hline
$\state{F}$ & 
$\state{F}$ & 
$0$ & 
$1$ \\ \hline
\end{tabular}
\qquad
\begin{tabular}[c]{lll|ll}\hline
\multicolumn{5}{l}{$\condprob[0]{{S_5}}{{P, Q, R}}$} \\ \hline\hline
\multicolumn{1}{l}{$P$} & 
\multicolumn{1}{l}{$Q$} & 
\multicolumn{1}{l|}{$R$} & 
\multicolumn{1}{l}{${S_5=\state{T}}$} & 
\multicolumn{1}{l}{${S_5=\state{F}}$} \\ \hline\hline
$\state{T}$ & 
$\state{T}$ & 
$\state{T}$ & 
$1$ & 
$0$ \\ \hline
$\state{T}$ & 
$\state{T}$ & 
$\state{F}$ & 
$1$ & 
$0$ \\ \hline
$\state{T}$ & 
$\state{F}$ & 
$\state{T}$ & 
$1$ & 
$0$ \\ \hline
$\state{T}$ & 
$\state{F}$ & 
$\state{F}$ & 
$0$ & 
$1$ \\ \hline
$\state{F}$ & 
$\state{T}$ & 
$\state{T}$ & 
$0$ & 
$1$ \\ \hline
$\state{F}$ & 
$\state{T}$ & 
$\state{F}$ & 
$0$ & 
$1$ \\ \hline
$\state{F}$ & 
$\state{F}$ & 
$\state{T}$ & 
$0$ & 
$1$ \\ \hline
$\state{F}$ & 
$\state{F}$ & 
$\state{F}$ & 
$0$ & 
$1$ \\ \hline
\end{tabular}
\end{equation}
\begin{equation}
\begin{tabular}[c]{l|ll}\hline
\multicolumn{3}{l}{$\condprob[0]{{S_6}}{{R}}$} \\ \hline\hline
\multicolumn{1}{l|}{$R$} & 
\multicolumn{1}{l}{${S_6=\state{T}}$} & 
\multicolumn{1}{l}{${S_6=\state{F}}$} \\ \hline\hline
$\state{T}$ & 
$1$ & 
$0$ \\ \hline
$\state{F}$ & 
$0$ & 
$1$ \\ \hline
\end{tabular}
\qquad
\begin{tabular}[c]{l|ll}\hline
\multicolumn{3}{l}{$\condprob[0]{{S_7}}{{R}}$} \\ \hline\hline
\multicolumn{1}{l|}{$R$} & 
\multicolumn{1}{l}{${S_7=\state{T}}$} & 
\multicolumn{1}{l}{${S_7=\state{F}}$} \\ \hline\hline
$\state{T}$ & 
$0$ & 
$1$ \\ \hline
$\state{F}$ & 
$1$ & 
$0$ \\ \hline
\end{tabular}
\qquad
\begin{tabular}[c]{l|ll}\hline
\multicolumn{3}{l}{$\condprob[0]{{S_8}}{{R}}$} \\ \hline\hline
\multicolumn{1}{l|}{$R$} & 
\multicolumn{1}{l}{${S_8=\state{T}}$} & 
\multicolumn{1}{l}{${S_8=\state{F}}$} \\ \hline\hline
$\state{T}$ & 
$0$ & 
$1$ \\ \hline
$\state{F}$ & 
$0$ & 
$1$ \\ \hline
\end{tabular}
\end{equation}
Now let us ask, is it possible that the propositions $P \wedge Q$,
$\neg(Q \wedge R) \wedge P$, and $R \wedge ( \neg P \rightarrow (R
\wedge Q))$ hold simultaneously?  Since these are the compound
formulas defined as $S_1$, $S_2$, and $S_3$ in Equation~\ref{eq:s3},
the computed probability table \prob{S_1,S_2,S_3} provides the answer:
\begin{equation}
\begin{tabular}[c]{lll|l}\hline
\multicolumn{1}{l}{$S_1$} & 
\multicolumn{1}{l}{$S_2$} & 
\multicolumn{1}{l|}{$S_3$} & 
\multicolumn{1}{l}{$\prob{{S_1, S_2, S_3}}$} \\ \hline\hline
$\state{T}$ & 
$\state{T}$ & 
$\state{T}$ & 
$0$ \\ \hline
$\state{T}$ & 
$\state{T}$ & 
$\state{F}$ & 
$x_{2}$ \\ \hline
$\state{T}$ & 
$\state{F}$ & 
$\state{T}$ & 
$x_{1}$ \\ \hline
$\state{T}$ & 
$\state{F}$ & 
$\state{F}$ & 
$0$ \\ \hline
$\state{F}$ & 
$\state{T}$ & 
$\state{T}$ & 
$x_{3}$ \\ \hline
$\state{F}$ & 
$\state{T}$ & 
$\state{F}$ & 
$x_{4}$ \\ \hline
$\state{F}$ & 
$\state{F}$ & 
$\state{T}$ & 
$x_{5}$ \\ \hline
$\state{F}$ & 
$\state{F}$ & 
$\state{F}$ & 
$x_{6} + x_{7} + x_{8}$ \\ \hline
\end{tabular}
\end{equation}
There are two impossible cases: that $S_1$, $S_2$, and $S_3$ hold
simultaneously; and that $S_1$ holds but $S_2$ and $S_3$ do not.

Next let us address the matter of `maximal consistency', which is a
property defined by Paris \emph{et al} to describe how compatible a
set of logical formulas are with one another.  According to their
definition we seek the maximum threshold value such that the
probability of each proposition in some set $\Gamma$ attains at least
that threshold; this maximum value is designated $\eta$.  This
definition of maximal consistency describes a linear optimization
problem.  In particular, the maximal consistency of the set
$\{S_1,S_2,S_3\}$ of formulas from Equation~\ref{eq:s3} can be
formulated as the linear optimization problem in which we ask for the
maximum value of a new parameter $z$ subject to the constraints
$\prob{S_1=\state{T}} \ge z$, $\prob{S_2=\state{T}} \ge z$, and
$\prob{S_3=\state{T}} \ge z$.  This is straightforward to construct
using parametric probability analysis.  First we compute the marginal
probability distribution on each compound formula in the set
$\{S_1,S_2,S_3\}$ as in Section~\ref{sec:symbolic}:
\begin{equation}
\begin{tabular}[c]{l|l}\hline
\multicolumn{1}{l|}{$S_1$} & 
\multicolumn{1}{l}{$\prob{{S_1}}$} \\ \hline\hline
$\state{T}$ & 
$x_{1} + x_{2}$ \\ \hline
$\state{F}$ & 
$x_{3} + x_{4} + x_{5} + x_{6} + x_{7} + x_{8}$ \\ \hline
\end{tabular}
\qquad
\begin{tabular}[c]{l|l}\hline
\multicolumn{1}{l|}{$S_2$} & 
\multicolumn{1}{l}{$\prob{{S_2}}$} \\ \hline\hline
$\state{T}$ & 
$x_{2} + x_{3} + x_{4}$ \\ \hline
$\state{F}$ & 
$x_{1} + x_{5} + x_{6} + x_{7} + x_{8}$ \\ \hline
\end{tabular}
\qquad
\begin{tabular}[c]{l|l}\hline
\multicolumn{1}{l|}{$S_3$} & 
\multicolumn{1}{l}{$\prob{{S_3}}$} \\ \hline\hline
$\state{T}$ & 
$x_{1} + x_{3} + x_{5}$ \\ \hline
$\state{F}$ & 
$x_{2} + x_{4} + x_{6} + x_{7} + x_{8}$ \\ \hline
\end{tabular}
\end{equation}
Then using the first element of each table in its respective
constraint (along with the usual constraints to enforce the laws of
probability) we construct the following optimization problem:
\begin{equation}
\begin{array}{r@{\quad}l}
\mbox{maximize}: & z \\
\mbox{subject to}:
& x_{1} + x_{2} + x_{3} + x_{4} + x_{5} + x_{6} + x_{7} + x_{8} = 1 \\
& x_{1} + x_{2} \ge z \\
& x_{2} + x_{3} + x_{4} \ge z \\
& x_{1} + x_{3} + x_{5} \ge z \\
\mbox{and}:
& 0 \le x_{1} \le 1 \\
& 0 \le x_{2} \le 1 \\
& 0 \le x_{3} \le 1 \\
& 0 \le x_{4} \le 1 \\
& 0 \le x_{5} \le 1 \\
& 0 \le x_{6} \le 1 \\
& 0 \le x_{7} \le 1 \\
& 0 \le x_{8} \le 1 \\
& 0 \le z \le 1
\end{array}
\label{eq:pp-eta}
\end{equation}
Solving this linear program yields the maximum value
\begin{math}
\eta =
0.667
\end{math}
which is achieved at the following point:
\begin{equation}
\left(\;
x_{1} = 0.333,\; x_{2} = 0.333,\; x_{3} = 0.333,\; x_{4} = 0.000,\; x_{5} = 0.000,\; x_{6} = 0.000,\; x_{7} = 0.000,\; x_{8} = 0.000,\; z = 0.667
\;\right)
\end{equation}
This agrees with the result reported in \cite{paris-muino}.  Note
that, although the author's polynomial-optimization solver used
floating-point arithmetic to compute the optimization results
displayed here, there are exact rational solvers for linear programs
which would instead calculate the solution precisely as the fraction
$2/3$.

\begin{table}
  \begin{tabular}{cll}\hline
    \textsc{Variable} & \textsc{Formula} & 
    \textsc{Threshold $\zeta$} \\ \hline\hline
    $S_4$ & $P \wedge R$ & 
\begin{math}
0.667
\end{math}
    \\ \hline
    $S_5$ & $P \wedge (Q \vee R)$ & 
\begin{math}
1.000
\end{math}
    \\ \hline
    $S_6$ & $R$ & 
\begin{math}
0.667
\end{math}
    \\ \hline
    $S_7$ & $\neg R$ & 
\begin{math}
0.333
\end{math}
    \\ \hline
    $S_8$ & $R \wedge \neg R$ & 
\begin{math}
0.000
\end{math}
    \\ \hline
  \end{tabular}
  \caption{Secondary threshold probabilities $\zeta$ for the formulas
    from Equation~\ref{eq:s8}, using the maximal consistency
    $
    \eta =
0.667
    $
    of the
    formulas $\{S_1,S_2,S_3\}$ from Equation~\ref{eq:s3} as the
    primary threshold probability.}
  \label{tbl:zeta}
\end{table}

Moving on, we can formulate additional optimization problems to
compute for other logical formulas the `secondary threshold
probability' designated $\zeta$, when the solution value $\eta$ above
is used as the `primary threshold probability'.  By the definition in
\cite{paris-muino} we seek the maximum probability $\zeta$ of each
queried formula subject to the constraints that every formula in the
designated set $\Gamma$ must attain probability at least $\eta$.
Again this definition describes certain linear optimization problems.
For this example, to find secondary threshold probabilities of the
formulas defined as $S_4$ through $S_8$ in Equation~\ref{eq:s8} we
first calculate their symbolic polynomial probabilities:
\begin{equation}
\begin{tabular}[c]{l|l}\hline
\multicolumn{1}{l|}{$S_4$} & 
\multicolumn{1}{l}{$\prob{{S_4}}$} \\ \hline\hline
$\state{T}$ & 
$x_{1} + x_{3}$ \\ \hline
$\state{F}$ & 
$x_{2} + x_{4} + x_{5} + x_{6} + x_{7} + x_{8}$ \\ \hline
\end{tabular}
\qquad
\begin{tabular}[c]{l|l}\hline
\multicolumn{1}{l|}{$S_5$} & 
\multicolumn{1}{l}{$\prob{{S_5}}$} \\ \hline\hline
$\state{T}$ & 
$x_{1} + x_{2} + x_{3}$ \\ \hline
$\state{F}$ & 
$x_{4} + x_{5} + x_{6} + x_{7} + x_{8}$ \\ \hline
\end{tabular}
\end{equation}
\begin{equation}
\begin{tabular}[c]{l|l}\hline
\multicolumn{1}{l|}{$S_6$} & 
\multicolumn{1}{l}{$\prob{{S_6}}$} \\ \hline\hline
$\state{T}$ & 
$x_{1} + x_{3} + x_{5} + x_{7}$ \\ \hline
$\state{F}$ & 
$x_{2} + x_{4} + x_{6} + x_{8}$ \\ \hline
\end{tabular}
\qquad
\begin{tabular}[c]{l|l}\hline
\multicolumn{1}{l|}{$S_7$} & 
\multicolumn{1}{l}{$\prob{{S_7}}$} \\ \hline\hline
$\state{T}$ & 
$x_{2} + x_{4} + x_{6} + x_{8}$ \\ \hline
$\state{F}$ & 
$x_{1} + x_{3} + x_{5} + x_{7}$ \\ \hline
\end{tabular}
\qquad
\begin{tabular}[c]{l|l}\hline
\multicolumn{1}{l|}{$S_8$} & 
\multicolumn{1}{l}{$\prob{{S_8}}$} \\ \hline\hline
$\state{T}$ & 
$0$ \\ \hline
$\state{F}$ & 
$x_{1} + x_{2} + x_{3} + x_{4} + x_{5} + x_{6} + x_{7} + x_{8}$ \\ \hline
\end{tabular}
\qquad
\end{equation}
Now to compute the secondary threshold probability for the query
$S_4:=(P \wedge R)$ relative to the set $\{S_1,S_2,S_3\}$ using the
primary threshold probability $\eta$, we query the maximum value of
\prob{S_4=\state{T}} subject to the constraints $\prob{S_1=\state{T}}
\ge z$, $\prob{S_2=\state{T}} \ge z$, and $\prob{S_3=\state{T}} \ge z$
where the auxiliary parameter $z$ is now fixed at the desired primary
threshold $\eta$.  Using the solution
\begin{math}
\eta =
0.667
\end{math}
from the optimization problem in Equation~\ref{eq:pp-eta} as the
primary threshold probability, we construct this optimization problem
to compute the secondary threshold probability for the formula
$S_4:=(P \wedge R)$:
\begin{equation}
\begin{array}{r@{\quad}l}
\mbox{maximize}: & x_{1} + x_{3} \\
\mbox{subject to}:
& x_{1} + x_{2} + x_{3} + x_{4} + x_{5} + x_{6} + x_{7} + x_{8} = 1 \\
& x_{1} + x_{2} \ge z \\
& x_{2} + x_{3} + x_{4} \ge z \\
& x_{1} + x_{3} + x_{5} \ge z \\
\mbox{and}:
& 0 \le x_{1} \le 1 \\
& 0 \le x_{2} \le 1 \\
& 0 \le x_{3} \le 1 \\
& 0 \le x_{4} \le 1 \\
& 0 \le x_{5} \le 1 \\
& 0 \le x_{6} \le 1 \\
& 0 \le x_{7} \le 1 \\
& 0 \le x_{8} \le 1 \\
& z = 0.667
\end{array}
\label{eq:pp-zeta}
\end{equation}
Solving this linear program yields the maximum value
\begin{math}
\zeta =
0.667
\end{math}
for the queried formula $S_4:=(P \wedge R)$.  By similar analysis we
compute secondary threshold probabilities for the remaining formulas
in Equation~\ref{eq:s8}, with the results displayed in
Table~\ref{tbl:zeta}.  Notably for the query formula $R$ the secondary
threshold $\zeta=1/3$ was reported in \cite{paris-muino}; but the
maximal solution to the optimization problem suggested by the text is
twice this value.

\subsection{Goodman's Hot Buttered Conditionals}

Next let us visit Goodman's treatment of counterfactual conditional
statements using his opening example from \emph{Fact, Fiction, and
  Forecast} \cite{goodman}.  In his own words:
\begin{quote}
  What, then, is the \emph{problem} about counterfactual conditionals?
  Let us confine ourselves to those in which antecedent and consequent
  are inalterably false---as, for example, when I say of a piece of
  butter that was eaten yesterday, and that had never been heated,
  \begin{quote}
    If that piece of butter had been heated to 150$^\circ$ F., it would
    have melted.
  \end{quote}
  Considered as truth-functional compounds, all counterfactuals are of
  course true, since their antecedents are false.  Hence
  \begin{quote}
    If that piece of butter had been heated to 150$^\circ$ F., it would
    not have melted
  \end{quote}
  would also hold.  Obviously something different is intended, and the
  problem is to define the circumstances under which a given
  counterfactual holds while the opposing counterfactual with the
  contradictory consequent fails to hold.
\end{quote}
Let us use $H$ to represent the proposition that the butter had been
heated, and $M$ for the proposition that it melted.  Each variable has
the set $\{\state{T},\state{F}\}$ of possible elementary values
representing truth and falsity.  It is indeed the case that given the
axiom that $H$ is false, both statements of material implication $H
\rightarrow M$ and $H \rightarrow \neg M$ are true.  In fact, Boolean
polynomial representation shows that the conjunction $(H \rightarrow
M)\wedge(H \rightarrow \neg M)$ simplifies to the negation $\neg H$.
In other words, interpreting them as statements of material
implication, the logical conjunction of the two opposing conditional
statements above is equivalent to the single unconditional statement:
\begin{quote}
  That piece of butter had not been heated to 150$^\circ$ F.
\end{quote}
If it desired to meet Goodman's dichotomy criterion (that one of the
conditional statements holds but the opposing statement with the
contradictory consequent does not hold) then it is not wise to model
his conditional sentences as statements of material implication.

Parametric probability networks and direct probability encoding
provide the desired `something different' to model conditional and
counterfactual statements.  The resulting models and analysis meet
Goodman's dichotomy criterion and give otherwise reasonable results.
To illustrate, let us build a parametric probability network including
the primary variables $H$ and $M$ to represent the propositional
variables, along with the primary variables $C_1:=(H \rightarrow M)$
and $C_2:=(H \rightarrow \neg M)$ to represent the corresponding
statements of material implication; we add parameters $x$, $y$, and
$z$.  Here is the network graph:
\begin{equation}
  \includegraphics{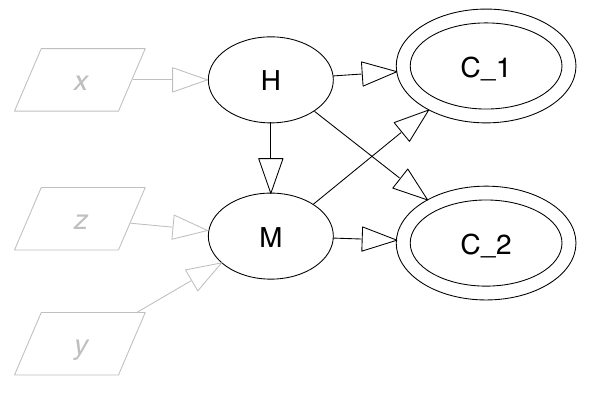}
\end{equation}
First we specify an uninformative prior distribution on $H$ and $M$
using parametric distributions \prob[0]{H} and \condprob[0]{M}{H},
with constraints $0 \le x \le 1$, $0 \le y \le 1$, and $0 \le z \le
1$:
\begin{equation}
\begin{tabular}[c]{l|l}\hline
\multicolumn{1}{l|}{$H$} & 
\multicolumn{1}{l}{$\prob[0]{{H}}$} \\ \hline\hline
$\state{T}$ & 
$x$ \\ \hline
$\state{F}$ & 
$1 - x$ \\ \hline
\end{tabular}
\qquad
\begin{tabular}[c]{l|ll}\hline
\multicolumn{3}{l}{$\condprob[0]{{M}}{{H}}$} \\ \hline\hline
\multicolumn{1}{l|}{$H$} & 
\multicolumn{1}{l}{${M=\state{T}}$} & 
\multicolumn{1}{l}{${M=\state{F}}$} \\ \hline\hline
$\state{T}$ & 
$y$ & 
$1 - y$ \\ \hline
$\state{F}$ & 
$z$ & 
$1 - z$ \\ \hline
\end{tabular}
\end{equation}
Let us interpret Goodman's first sentence as the constraint
$\condprob[0]{M=\state{T}}{H=\state{T}}=1$ that the input value of the
conditional probability that $M$ is true given that $H$ is true is
one; likewise let us interpret his second sentence as the constraint
$\condprob[0]{M=\state{T}}{H=\state{T}}=0$ or the equivalent
$\condprob[0]{M=\state{F}}{H=\state{T}}=1$.  Now it is a matter of
elementary algebra that these statements are inconsistent: the first
says $y=1$ and the others say $y=0$.  This algebraic inconsistency
does not involve the parameter $x$ that encodes the prior probability
$\prob[0]{H=\state{T}}$ that the butter had been heated.  In other
words, these opposing conditional sentences are inconsistent
specifically because their consequents disagree and not because they
are counterfactual.  Additionally, in this interpretation these
conditional sentences do not have existential import.  Neither
constraint $y=1$ nor $y=0$ requires that $H$ must certainly be true
($x=1$), nor even that $H$ must possibly be true ($x>0$); these
constraints on $y$ do not affect $x$ at all.

Next let us compute the conditional probability \condprob{M}{H} of
whether the butter melted given whether it had been heated, alongside
the marginal probability \prob{M} of whether the butter melted
(integrating the cases that it had or had not been heated) and the
marginal probability \prob{H} of whether the butter had been heated:
\begin{equation}
\begin{tabular}[c]{l|ll}\hline
\multicolumn{3}{l}{$\condprob{{M}}{{H}}$} \\ \hline\hline
\multicolumn{1}{l|}{$H$} & 
\multicolumn{1}{l}{${M=\state{T}}$} & 
\multicolumn{1}{l}{${M=\state{F}}$} \\ \hline\hline
$\state{T}$ & 
$\left(x y\right) / \left(x\right)$ & 
$\left(x - x y\right) / \left(x\right)$ \\ \hline
$\state{F}$ & 
$\left(z - x z\right) / \left(1 - x\right)$ & 
$\left(1 - x - z + x z\right) / \left(1 - x\right)$ \\ \hline
\end{tabular}
\qquad
\begin{tabular}[c]{l|l}\hline
\multicolumn{1}{l|}{$M$} & 
\multicolumn{1}{l}{$\prob{{M}}$} \\ \hline\hline
$\state{T}$ & 
$z + x y - x z$ \\ \hline
$\state{F}$ & 
$1 - z - x y + x z$ \\ \hline
\end{tabular}
\qquad
\begin{tabular}[c]{l|l}\hline
\multicolumn{1}{l|}{$H$} & 
\multicolumn{1}{l}{$\prob{{H}}$} \\ \hline\hline
$\state{T}$ & 
$x$ \\ \hline
$\state{F}$ & 
$1 - x$ \\ \hline
\end{tabular}
\end{equation}
These computed probability tables already tell an interesting story
about Goodman's conditional sentences.  First, the output probability
\prob{H=\state{T}} that the butter had been heated has exactly the
same value as the input probability \prob[0]{H=\state{T}}, namely the
parameter $x$.  However the output probability
\condprob{M=\state{T}}{H=\state{T}} differs in a subtle way from the
input probability \condprob[0]{M=\state{T}}{H=\state{T}}; the input
value is the parameter $y$ but the output value is the quotient
$xy/x$.  Therefore if it is certain that the butter had not been
heated ($x=0$) then the computed conditional probability $xy/x$ that
the butter melted given this now-impossible condition is algebraically
indeterminate: it is the quotient $0/0$.

Notwithstanding this exceptional \emph{conditional} probability, the
overall \emph{marginal} probability that the butter melted subject to
the constraint $x=0$ that the butter certainly had not been heated
does not involve division by zero; in fact it does not involve
division at all.  The computed probability $\prob{M=\state{T}}$, whose
value is the polynomial
\begin{math}
z + x y - x z
\end{math},
simplifies to $z$ when $x=0$.  In other words, subject to the
constraint that the butter certainly had not been heated, the
probability that the butter melted is exactly the value $z$ specified
as the input component probability
\condprob[0]{M=\state{T}}{H=\state{F}}.  If the user has not provided
any more information about $z$ then its precise value remains
indeterminate; it is only constrained by $0 \le z \le 1$ to satisfy
the laws of probability.  In the terminology of
Section~\ref{sec:subjunctive}, we have just considered subjunctive and
imperative modes of asking the same question: What is the probability
that $M$ is true, given the condition that $H$ is false?  In either
formulation the answer is the same: The queried probability is
indeterminate.

We have just analyzed Goodman's counterfactual conditional sentences
using parametric probability directly, without intermediate formulas
from the propositional calculus.  However it is instructive to embed
the formulas $H \rightarrow M$ and $H \rightarrow \neg M$ into our
probability network and to compute the probabilities associated with
them.  Using the definitions $C_1:=(H \rightarrow M)$ and $C_2:=(H
\rightarrow \neg M)$ and the method of Section~\ref{sec:conditional}
we construct the following conditional probability tables (here
labeled with the embedded formulas instead of with the
primary-variable names $C_1$ and $C_2$):
\begin{equation}
\begin{tabular}[c]{ll|ll}\hline
\multicolumn{4}{l}{$\condprob[0]{{(H \rightarrow M)}}{{H, M}}$} \\ \hline\hline
\multicolumn{1}{l}{$H$} & 
\multicolumn{1}{l|}{$M$} & 
\multicolumn{1}{l}{${(H \rightarrow M)=\state{T}}$} & 
\multicolumn{1}{l}{${(H \rightarrow M)=\state{F}}$} \\ \hline\hline
$\state{T}$ & 
$\state{T}$ & 
$1$ & 
$0$ \\ \hline
$\state{T}$ & 
$\state{F}$ & 
$0$ & 
$1$ \\ \hline
$\state{F}$ & 
$\state{T}$ & 
$1$ & 
$0$ \\ \hline
$\state{F}$ & 
$\state{F}$ & 
$1$ & 
$0$ \\ \hline
\end{tabular}
\qquad
\begin{tabular}[c]{ll|ll}\hline
\multicolumn{4}{l}{$\condprob[0]{{(H \rightarrow \neg M)}}{{H, M}}$} \\ \hline\hline
\multicolumn{1}{l}{$H$} & 
\multicolumn{1}{l|}{$M$} & 
\multicolumn{1}{l}{${(H \rightarrow \neg M)=\state{T}}$} & 
\multicolumn{1}{l}{${(H \rightarrow \neg M)=\state{F}}$} \\ \hline\hline
$\state{T}$ & 
$\state{T}$ & 
$0$ & 
$1$ \\ \hline
$\state{T}$ & 
$\state{F}$ & 
$1$ & 
$0$ \\ \hline
$\state{F}$ & 
$\state{T}$ & 
$1$ & 
$0$ \\ \hline
$\state{F}$ & 
$\state{F}$ & 
$1$ & 
$0$ \\ \hline
\end{tabular}
\end{equation}
We compute the marginal probabilities that each embedded formula is
true:
\begin{equation}
\begin{tabular}[c]{l|l}\hline
\multicolumn{1}{l|}{$(H \rightarrow M)$} & 
\multicolumn{1}{l}{$\prob{{(H \rightarrow M)}}$} \\ \hline\hline
$\state{T}$ & 
$1 - x + x y$ \\ \hline
$\state{F}$ & 
$x - x y$ \\ \hline
\end{tabular}
\qquad
\begin{tabular}[c]{l|l}\hline
\multicolumn{1}{l|}{$(H \rightarrow \neg M)$} & 
\multicolumn{1}{l}{$\prob{{(H \rightarrow \neg M)}}$} \\ \hline\hline
$\state{T}$ & 
$1 - x y$ \\ \hline
$\state{F}$ & 
$x y$ \\ \hline
\end{tabular}
\label{eq:tc1}
\end{equation}
Solving either equation
\begin{math}
1 - x + x y
=
1 - x y
\end{math}
or
\begin{math}
x - x y
=
x y
\end{math}
reveals that there are exactly two cases in which the opposing
statements of material implication have the same probability of truth:
when $x=0$ (in which case each statement is true with certainty) and
when $y=\frac{1}{2}$ (in which case each statement is true with
probability $1-\frac{1}{2}x$).  In other words, if it is certain
\emph{a priori} that the butter had not been heated ($x=0$) then both
formulas $H \rightarrow M$ and $H \rightarrow \neg M$ must certainly
be true.  However if there is exactly a 50\% chance that the butter
melted after it had been heated ($y=\frac{1}{2}$) then it is also
equally likely that the formulas $H \rightarrow M$ and $H \rightarrow
\neg M$ are true; but now their mutual probability ($1-\frac{1}{2}x$)
is one minus half the prior probability that the butter had been
heated.

The relation that both formulas $H \rightarrow M$ and $H \rightarrow
\neg M$ are \emph{equally likely} to be true is different from the
relation that both formulas are \emph{simultaneously} true; both
relations occur when $x=0$ but only the former occurs when
$y=\frac{1}{2}$ and $x \neq 0$.  Evaluating the joint probability
\prob{H,C_1,C_2} of $H$ and both statements of material implication
provides additional detail that is hidden in the marginal
probabilities above:
\begin{equation}
\begin{tabular}[c]{r|lll|l}\hline
\multicolumn{1}{l|}{\scshape {Index}} & 
\multicolumn{1}{l}{$H$} & 
\multicolumn{1}{l}{$(H \rightarrow M)$} & 
\multicolumn{1}{l|}{$(H \rightarrow \neg M)$} & 
\multicolumn{1}{l}{$\prob{{H, (H \rightarrow M), (H \rightarrow \neg M)}}$} \\ \hline\hline
1 & 
$\state{T}$ & 
$\state{T}$ & 
$\state{T}$ & 
$0$ \\ \hline
2 & 
$\state{T}$ & 
$\state{T}$ & 
$\state{F}$ & 
$x y$ \\ \hline
3 & 
$\state{T}$ & 
$\state{F}$ & 
$\state{T}$ & 
$x - x y$ \\ \hline
4 & 
$\state{T}$ & 
$\state{F}$ & 
$\state{F}$ & 
$0$ \\ \hline
5 & 
$\state{F}$ & 
$\state{T}$ & 
$\state{T}$ & 
$1 - x$ \\ \hline
6 & 
$\state{F}$ & 
$\state{T}$ & 
$\state{F}$ & 
$0$ \\ \hline
7 & 
$\state{F}$ & 
$\state{F}$ & 
$\state{T}$ & 
$0$ \\ \hline
8 & 
$\state{F}$ & 
$\state{F}$ & 
$\state{F}$ & 
$0$ \\ \hline
\end{tabular}
\label{eq:tc12}
\end{equation}
This probability table shows that it is impossible that the antecedent
$H$ and both opposing statements of material implication are
simultaneously true (the probability of this event, which appears in
row~1 of the table, is identically zero).  Moreover when $H$ is true
exactly one of the opposing statements of material implication must
hold (see rows 1 through 4).  But when $H$ is false both statements of
material implication must be true (see rows 5 through 8).  To recover
the probabilities in Equation~\ref{eq:tc1} from the probabilities in
Equation~\ref{eq:tc12} it is necessary to add appropriate elements of
the latter probability table.  For example $\prob{(H \rightarrow
  M)=\state{T}}$ is given by the sum
\begin{math}
(
0
) + (
x y
) + (
1 - x
) + (
0
)
\end{math}
of the polynomials from rows 1, 2, 5, and 6, which yields
\begin{math}
1 - x + x y
\end{math}.

\subsection{Smullyan's Knights, Knaves, and Zombies}

Finally we consider two of Smullyan's problems from \emph{What Is the
  Name of This Book?} which were also analyzed by Kolany using a
rather different technique \cite{smullyan-name,kolany}.  First a basic
knights and knaves problem for which we must simply answer a
probability-table query; and second a zombie problem in which we must
first answer a probability-table query and then perform a search to
find certain values of the parameters in the resulting probability
table.

\paragraph{We Are the Knights Who Say\ldots}

The background for the first problem (number~36 in
\cite{smullyan-name}) is that on the imagined island, knights always
tell the truth and knaves always lie.
\begin{quote}
  Once when I visited the island of knights and knaves, I came across
  two of the inhabitants resting under a tree.  I asked one of them,
  ``Is either of you a knight?''  He responded, and I knew the answer
  to my question.

  What is the person to whom I addressed the question---is he a knight
  or knave; And what is the other one?
\end{quote}
In the parametric probability model let us use the variable $A$ to
represent the proposition that the respondent is a knight, and $B$ for
the proposition that the other inhabitant is a knight.  We define
$Q:=(A \vee B)$ to represent the true answer to the question of
whether either inhabitant is a knight.  We introduce $R$ to represent
the response that is given: $R$ is true if the response is `yes' and
false if the response is `no'.  By Smullyan's rules for the island, if
$Q$ were true then the inhabitant would respond `yes' if and only if
he were a knight; and if $Q$ were false then the responses would be
opposite.  Hence the definition $R:=(A \leftrightarrow Q)$ using the
biconditional.  Here is the network graph, which also includes the
parameters $x_1$ through $x_4$:
\begin{equation}
  \includegraphics{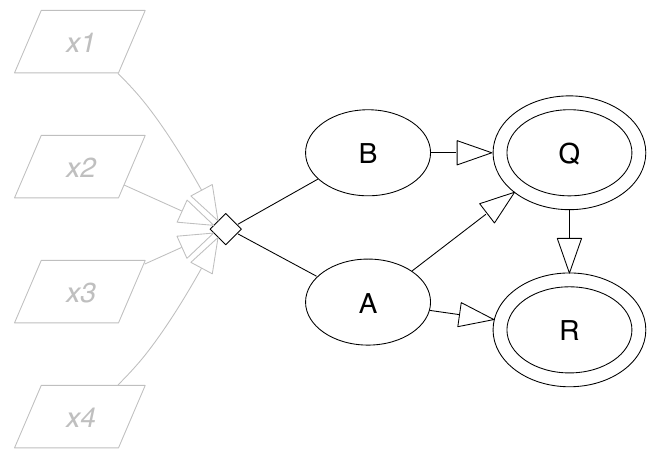}
\end{equation}
We specify an uninformative parametric probability distribution on $A$
and $B$ using the $x_i$ parameters with the usual constraints $0 \le
x_i \le 1$ and $\sum_i x_i = 1$ as in Section~\ref{sec:prior}; we
construct the appropriate component probability tables for $Q$ and $R$
according to Section~\ref{sec:conditional}:
\begin{equation}
\begin{tabular}[c]{ll|l}\hline
\multicolumn{1}{l}{$A$} & 
\multicolumn{1}{l|}{$B$} & 
\multicolumn{1}{l}{$\prob[0]{{A, B}}$} \\ \hline\hline
$\state{T}$ & 
$\state{T}$ & 
$x_{1}$ \\ \hline
$\state{T}$ & 
$\state{F}$ & 
$x_{2}$ \\ \hline
$\state{F}$ & 
$\state{T}$ & 
$x_{3}$ \\ \hline
$\state{F}$ & 
$\state{F}$ & 
$x_{4}$ \\ \hline
\end{tabular}
\qquad
\begin{tabular}[c]{ll|ll}\hline
\multicolumn{4}{l}{$\condprob[0]{{Q}}{{A, B}}$} \\ \hline\hline
\multicolumn{1}{l}{$A$} & 
\multicolumn{1}{l|}{$B$} & 
\multicolumn{1}{l}{${Q=\state{T}}$} & 
\multicolumn{1}{l}{${Q=\state{F}}$} \\ \hline\hline
$\state{T}$ & 
$\state{T}$ & 
$1$ & 
$0$ \\ \hline
$\state{T}$ & 
$\state{F}$ & 
$1$ & 
$0$ \\ \hline
$\state{F}$ & 
$\state{T}$ & 
$1$ & 
$0$ \\ \hline
$\state{F}$ & 
$\state{F}$ & 
$0$ & 
$1$ \\ \hline
\end{tabular}
\qquad
\begin{tabular}[c]{ll|ll}\hline
\multicolumn{4}{l}{$\condprob[0]{{R}}{{Q, A}}$} \\ \hline\hline
\multicolumn{1}{l}{$Q$} & 
\multicolumn{1}{l|}{$A$} & 
\multicolumn{1}{l}{${R=\state{T}}$} & 
\multicolumn{1}{l}{${R=\state{F}}$} \\ \hline\hline
$\state{T}$ & 
$\state{T}$ & 
$1$ & 
$0$ \\ \hline
$\state{T}$ & 
$\state{F}$ & 
$0$ & 
$1$ \\ \hline
$\state{F}$ & 
$\state{T}$ & 
$0$ & 
$1$ \\ \hline
$\state{F}$ & 
$\state{F}$ & 
$1$ & 
$0$ \\ \hline
\end{tabular}
\end{equation}
We compute the joint probability of the identities $A$ and $B$
conditioned on the response $R$:
\begin{equation}
\begin{tabular}[c]{l|llll}\hline
\multicolumn{5}{l}{$\condprob{{A, B}}{{R}}$} \\ \hline\hline
\multicolumn{1}{l|}{$R$} & 
\multicolumn{1}{l}{${A=\state{T}, B=\state{T}}$} & 
\multicolumn{1}{l}{${A=\state{T}, B=\state{F}}$} & 
\multicolumn{1}{l}{${A=\state{F}, B=\state{T}}$} & 
\multicolumn{1}{l}{${A=\state{F}, B=\state{F}}$} \\ \hline\hline
$\state{T}$ & 
$\left(x_{1}\right) / \left(x_{1} + x_{2} + x_{4}\right)$ & 
$\left(x_{2}\right) / \left(x_{1} + x_{2} + x_{4}\right)$ & 
$\left(0\right) / \left(x_{1} + x_{2} + x_{4}\right)$ & 
$\left(x_{4}\right) / \left(x_{1} + x_{2} + x_{4}\right)$ \\ \hline
$\state{F}$ & 
$\left(0\right) / \left(x_{3}\right)$ & 
$\left(0\right) / \left(x_{3}\right)$ & 
$\left(x_{3}\right) / \left(x_{3}\right)$ & 
$\left(0\right) / \left(x_{3}\right)$ \\ \hline
\end{tabular}
\end{equation}
This result table for \condprob{A,B}{R} gives the solution.  An answer
of `yes' ($R=\state{T}$) leaves three possible configurations of
identities, but an answer of `no' ($R=\state{F}$) leaves only one
possible configuration of identities: that the responder is a knave
and the other inhabitant is a knight ($A=\state{F},B=\state{T}$).
Therefore if the identities of the inhabitants are known with
certainty after the response, then the response must have been `no',
the responder must be a knave, and his fellow inhabitant must be a
knight.

\paragraph{Zombieland}

For the second problem (number~160 in \cite{smullyan-name}) we visit
Smullyan's island of zombies.  The custom here is that zombies always
lie and humans always tell the truth.  However instead of `yes' and
`no' the inhabitants answer questions with `Bal' or `Da'; one means
`yes' and the other means `no' but we do not know which is which.  The
problem asks the following:
\begin{quote}
Suppose you are not interested in what ``Bal'' means, but only in
whether the speaker is a zombie.  How can you find this out in only
one question? (Again, he will answer ``Bal'' or ``Da.'')
\end{quote}
For the parametric probability model of this problem we use the
variable $H$ to represent whether the speaker is human ($H=\state{T}$)
or zombie; $B$ to represent whether `Bal' means `yes' ($B=\state{T}$)
or `no'; $Q$ for the true answer to the unknown question that is
sought; and $R$ for whether the speaker gives the response `Bal'
($R=\state{T}$) or `Da'.  Here is the network graph:
\begin{equation}
  \includegraphics{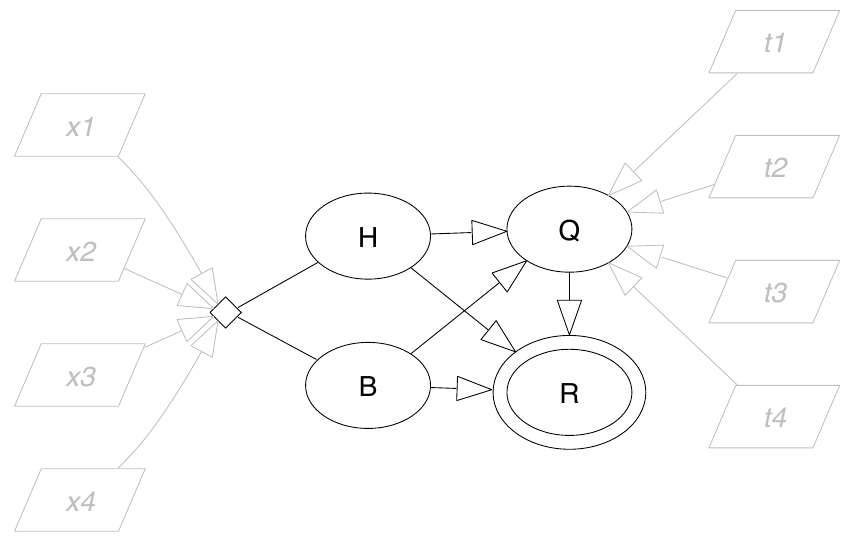}
\end{equation}
The prior distribution \prob[0]{H,B} on $H$ and $B$ is specified using
parameters $x_1$ through $x_4$ with $0 \le x_i \le 1$ and $\sum_i x_i
= 1$.  Additional parameters $t_1$ through $t_4$ with each $t_i \in
\{0,1\}$ are introduced as described in Section~\ref{sec:unknown} to
specify the component probability table \condprob[0]{Q}{H,B} of the
unknown question $Q$.  Here are both component probability tables:
\begin{equation}
\begin{tabular}[c]{ll|l}\hline
\multicolumn{1}{l}{$H$} & 
\multicolumn{1}{l|}{$B$} & 
\multicolumn{1}{l}{$\prob[0]{{H, B}}$} \\ \hline\hline
$\state{T}$ & 
$\state{T}$ & 
$x_{1}$ \\ \hline
$\state{T}$ & 
$\state{F}$ & 
$x_{2}$ \\ \hline
$\state{F}$ & 
$\state{T}$ & 
$x_{3}$ \\ \hline
$\state{F}$ & 
$\state{F}$ & 
$x_{4}$ \\ \hline
\end{tabular}
\qquad
\begin{tabular}[c]{ll|ll}\hline
\multicolumn{4}{l}{$\condprob[0]{{Q}}{{H, B}}$} \\ \hline\hline
\multicolumn{1}{l}{$H$} & 
\multicolumn{1}{l|}{$B$} & 
\multicolumn{1}{l}{${Q=\state{T}}$} & 
\multicolumn{1}{l}{${Q=\state{F}}$} \\ \hline\hline
$\state{T}$ & 
$\state{T}$ & 
$t_{1}$ & 
$1 - t_{1}$ \\ \hline
$\state{T}$ & 
$\state{F}$ & 
$t_{2}$ & 
$1 - t_{2}$ \\ \hline
$\state{F}$ & 
$\state{T}$ & 
$t_{3}$ & 
$1 - t_{3}$ \\ \hline
$\state{F}$ & 
$\state{F}$ & 
$t_{4}$ & 
$1 - t_{4}$ \\ \hline
\end{tabular}
\label{eq:cpt-q}
\end{equation}
The speaker's response $R$ can be modeled in the propositional
calculus using nested biconditionals to represent the rules of the
island.  The formula $H \leftrightarrow Q$ reveals whether the the
speaker will answer in the affirmative.  For example a zombie
($H=\state{F}$) will dishonestly provide an affirmative answer to a
question that is actually false ($Q=\state{F}$) but a human will
honestly provide a negative answer to a question that is actually
false.  Relating this inner biconditional to whether `Bal' means `yes'
using another biconditional $(H \leftrightarrow Q) \leftrightarrow B$
then reveals whether the speaker will answer `Bal' (if this outer
biconditional is true) or `Da'.  For example, if the speaker will
answer in the affirmative, then he will respond with `Bal' if and only
if that means `yes'.  The following component probability table
\condprob[0]{R}{H,Q,B} constructed as in Section~\ref{sec:conditional}
implements the definition $R:=((H \leftrightarrow Q) \leftrightarrow
B)$:
\begin{equation}
\begin{tabular}[c]{lll|ll}\hline
\multicolumn{5}{l}{$\condprob[0]{{R}}{{H, Q, B}}$} \\ \hline\hline
\multicolumn{1}{l}{$H$} & 
\multicolumn{1}{l}{$Q$} & 
\multicolumn{1}{l|}{$B$} & 
\multicolumn{1}{l}{${R=\state{T}}$} & 
\multicolumn{1}{l}{${R=\state{F}}$} \\ \hline\hline
$\state{T}$ & 
$\state{T}$ & 
$\state{T}$ & 
$1$ & 
$0$ \\ \hline
$\state{T}$ & 
$\state{T}$ & 
$\state{F}$ & 
$0$ & 
$1$ \\ \hline
$\state{T}$ & 
$\state{F}$ & 
$\state{T}$ & 
$0$ & 
$1$ \\ \hline
$\state{T}$ & 
$\state{F}$ & 
$\state{F}$ & 
$1$ & 
$0$ \\ \hline
$\state{F}$ & 
$\state{T}$ & 
$\state{T}$ & 
$0$ & 
$1$ \\ \hline
$\state{F}$ & 
$\state{T}$ & 
$\state{F}$ & 
$1$ & 
$0$ \\ \hline
$\state{F}$ & 
$\state{F}$ & 
$\state{T}$ & 
$1$ & 
$0$ \\ \hline
$\state{F}$ & 
$\state{F}$ & 
$\state{F}$ & 
$0$ & 
$1$ \\ \hline
\end{tabular}
\end{equation}

\begin{table}
\begin{tabular}{@{}l@{\qquad}l@{}}
(a) \quad
\begin{tabular}[c]{r|llll|l}\hline
\multicolumn{1}{l|}{\scshape {Index}} & 
\multicolumn{1}{l}{$t_{1}$} & 
\multicolumn{1}{l}{$t_{2}$} & 
\multicolumn{1}{l}{$t_{3}$} & 
\multicolumn{1}{l|}{$t_{4}$} & 
\multicolumn{1}{l}{$\condutil[0]{{\prob{R=\state{T},H=\state{T}}}}{{t_{1}, t_{2}, t_{3}, t_{4}}}$} \\ \hline\hline
1 & 
$\state{0}$ & 
$\state{0}$ & 
$\state{0}$ & 
$\state{0}$ & 
$x_{2}$ \\ \hline
2 & 
$\state{0}$ & 
$\state{0}$ & 
$\state{0}$ & 
$\state{1}$ & 
$x_{2}$ \\ \hline
3 & 
$\state{0}$ & 
$\state{0}$ & 
$\state{1}$ & 
$\state{0}$ & 
$x_{2}$ \\ \hline
4 & 
$\state{0}$ & 
$\state{0}$ & 
$\state{1}$ & 
$\state{1}$ & 
$x_{2}$ \\ \hline
5 & 
$\state{0}$ & 
$\state{1}$ & 
$\state{0}$ & 
$\state{0}$ & 
$0$ \\ \hline
6 & 
$\state{0}$ & 
$\state{1}$ & 
$\state{0}$ & 
$\state{1}$ & 
$0$ \\ \hline
7 & 
$\state{0}$ & 
$\state{1}$ & 
$\state{1}$ & 
$\state{0}$ & 
$0$ \\ \hline
8 & 
$\state{0}$ & 
$\state{1}$ & 
$\state{1}$ & 
$\state{1}$ & 
$0$ \\ \hline
9 & 
$\state{1}$ & 
$\state{0}$ & 
$\state{0}$ & 
$\state{0}$ & 
$x_{1} + x_{2}$ \\ \hline
10 & 
$\state{1}$ & 
$\state{0}$ & 
$\state{0}$ & 
$\state{1}$ & 
$x_{1} + x_{2}$ \\ \hline
11 & 
$\state{1}$ & 
$\state{0}$ & 
$\state{1}$ & 
$\state{0}$ & 
$x_{1} + x_{2}$ \\ \hline
12 & 
$\state{1}$ & 
$\state{0}$ & 
$\state{1}$ & 
$\state{1}$ & 
$x_{1} + x_{2}$ \\ \hline
13 & 
$\state{1}$ & 
$\state{1}$ & 
$\state{0}$ & 
$\state{0}$ & 
$x_{1}$ \\ \hline
14 & 
$\state{1}$ & 
$\state{1}$ & 
$\state{0}$ & 
$\state{1}$ & 
$x_{1}$ \\ \hline
15 & 
$\state{1}$ & 
$\state{1}$ & 
$\state{1}$ & 
$\state{0}$ & 
$x_{1}$ \\ \hline
16 & 
$\state{1}$ & 
$\state{1}$ & 
$\state{1}$ & 
$\state{1}$ & 
$x_{1}$ \\ \hline
\end{tabular}
&
(b) \quad
\begin{tabular}[c]{r|llll|l}\hline
\multicolumn{1}{l|}{\scshape {Index}} & 
\multicolumn{1}{l}{$t_{1}$} & 
\multicolumn{1}{l}{$t_{2}$} & 
\multicolumn{1}{l}{$t_{3}$} & 
\multicolumn{1}{l|}{$t_{4}$} & 
\multicolumn{1}{l}{$\condutil[0]{{\prob{R=\state{T},H=\state{F}}}}{{t_{1}, t_{2}, t_{3}, t_{4}}}$} \\ \hline\hline
1 & 
$\state{0}$ & 
$\state{0}$ & 
$\state{0}$ & 
$\state{0}$ & 
$x_{3}$ \\ \hline
2 & 
$\state{0}$ & 
$\state{0}$ & 
$\state{0}$ & 
$\state{1}$ & 
$x_{3} + x_{4}$ \\ \hline
3 & 
$\state{0}$ & 
$\state{0}$ & 
$\state{1}$ & 
$\state{0}$ & 
$0$ \\ \hline
4 & 
$\state{0}$ & 
$\state{0}$ & 
$\state{1}$ & 
$\state{1}$ & 
$x_{4}$ \\ \hline
5 & 
$\state{0}$ & 
$\state{1}$ & 
$\state{0}$ & 
$\state{0}$ & 
$x_{3}$ \\ \hline
6 & 
$\state{0}$ & 
$\state{1}$ & 
$\state{0}$ & 
$\state{1}$ & 
$x_{3} + x_{4}$ \\ \hline
7 & 
$\state{0}$ & 
$\state{1}$ & 
$\state{1}$ & 
$\state{0}$ & 
$0$ \\ \hline
8 & 
$\state{0}$ & 
$\state{1}$ & 
$\state{1}$ & 
$\state{1}$ & 
$x_{4}$ \\ \hline
9 & 
$\state{1}$ & 
$\state{0}$ & 
$\state{0}$ & 
$\state{0}$ & 
$x_{3}$ \\ \hline
10 & 
$\state{1}$ & 
$\state{0}$ & 
$\state{0}$ & 
$\state{1}$ & 
$x_{3} + x_{4}$ \\ \hline
11 & 
$\state{1}$ & 
$\state{0}$ & 
$\state{1}$ & 
$\state{0}$ & 
$0$ \\ \hline
12 & 
$\state{1}$ & 
$\state{0}$ & 
$\state{1}$ & 
$\state{1}$ & 
$x_{4}$ \\ \hline
13 & 
$\state{1}$ & 
$\state{1}$ & 
$\state{0}$ & 
$\state{0}$ & 
$x_{3}$ \\ \hline
14 & 
$\state{1}$ & 
$\state{1}$ & 
$\state{0}$ & 
$\state{1}$ & 
$x_{3} + x_{4}$ \\ \hline
15 & 
$\state{1}$ & 
$\state{1}$ & 
$\state{1}$ & 
$\state{0}$ & 
$0$ \\ \hline
16 & 
$\state{1}$ & 
$\state{1}$ & 
$\state{1}$ & 
$\state{1}$ & 
$x_{4}$ \\ \hline
\end{tabular}
\\ \\
(c) \quad
\begin{tabular}[c]{r|llll|l}\hline
\multicolumn{1}{l|}{\scshape {Index}} & 
\multicolumn{1}{l}{$t_{1}$} & 
\multicolumn{1}{l}{$t_{2}$} & 
\multicolumn{1}{l}{$t_{3}$} & 
\multicolumn{1}{l|}{$t_{4}$} & 
\multicolumn{1}{l}{$\condutil[0]{{\prob{R=\state{F},H=\state{T}}}}{{t_{1}, t_{2}, t_{3}, t_{4}}}$} \\ \hline\hline
1 & 
$\state{0}$ & 
$\state{0}$ & 
$\state{0}$ & 
$\state{0}$ & 
$x_{1}$ \\ \hline
2 & 
$\state{0}$ & 
$\state{0}$ & 
$\state{0}$ & 
$\state{1}$ & 
$x_{1}$ \\ \hline
3 & 
$\state{0}$ & 
$\state{0}$ & 
$\state{1}$ & 
$\state{0}$ & 
$x_{1}$ \\ \hline
4 & 
$\state{0}$ & 
$\state{0}$ & 
$\state{1}$ & 
$\state{1}$ & 
$x_{1}$ \\ \hline
5 & 
$\state{0}$ & 
$\state{1}$ & 
$\state{0}$ & 
$\state{0}$ & 
$x_{1} + x_{2}$ \\ \hline
6 & 
$\state{0}$ & 
$\state{1}$ & 
$\state{0}$ & 
$\state{1}$ & 
$x_{1} + x_{2}$ \\ \hline
7 & 
$\state{0}$ & 
$\state{1}$ & 
$\state{1}$ & 
$\state{0}$ & 
$x_{1} + x_{2}$ \\ \hline
8 & 
$\state{0}$ & 
$\state{1}$ & 
$\state{1}$ & 
$\state{1}$ & 
$x_{1} + x_{2}$ \\ \hline
9 & 
$\state{1}$ & 
$\state{0}$ & 
$\state{0}$ & 
$\state{0}$ & 
$0$ \\ \hline
10 & 
$\state{1}$ & 
$\state{0}$ & 
$\state{0}$ & 
$\state{1}$ & 
$0$ \\ \hline
11 & 
$\state{1}$ & 
$\state{0}$ & 
$\state{1}$ & 
$\state{0}$ & 
$0$ \\ \hline
12 & 
$\state{1}$ & 
$\state{0}$ & 
$\state{1}$ & 
$\state{1}$ & 
$0$ \\ \hline
13 & 
$\state{1}$ & 
$\state{1}$ & 
$\state{0}$ & 
$\state{0}$ & 
$x_{2}$ \\ \hline
14 & 
$\state{1}$ & 
$\state{1}$ & 
$\state{0}$ & 
$\state{1}$ & 
$x_{2}$ \\ \hline
15 & 
$\state{1}$ & 
$\state{1}$ & 
$\state{1}$ & 
$\state{0}$ & 
$x_{2}$ \\ \hline
16 & 
$\state{1}$ & 
$\state{1}$ & 
$\state{1}$ & 
$\state{1}$ & 
$x_{2}$ \\ \hline
\end{tabular}
&
(d) \quad
\begin{tabular}[c]{r|llll|l}\hline
\multicolumn{1}{l|}{\scshape {Index}} & 
\multicolumn{1}{l}{$t_{1}$} & 
\multicolumn{1}{l}{$t_{2}$} & 
\multicolumn{1}{l}{$t_{3}$} & 
\multicolumn{1}{l|}{$t_{4}$} & 
\multicolumn{1}{l}{$\condutil[0]{{\prob{R=\state{F},H=\state{F}}}}{{t_{1}, t_{2}, t_{3}, t_{4}}}$} \\ \hline\hline
1 & 
$\state{0}$ & 
$\state{0}$ & 
$\state{0}$ & 
$\state{0}$ & 
$x_{4}$ \\ \hline
2 & 
$\state{0}$ & 
$\state{0}$ & 
$\state{0}$ & 
$\state{1}$ & 
$0$ \\ \hline
3 & 
$\state{0}$ & 
$\state{0}$ & 
$\state{1}$ & 
$\state{0}$ & 
$x_{3} + x_{4}$ \\ \hline
4 & 
$\state{0}$ & 
$\state{0}$ & 
$\state{1}$ & 
$\state{1}$ & 
$x_{3}$ \\ \hline
5 & 
$\state{0}$ & 
$\state{1}$ & 
$\state{0}$ & 
$\state{0}$ & 
$x_{4}$ \\ \hline
6 & 
$\state{0}$ & 
$\state{1}$ & 
$\state{0}$ & 
$\state{1}$ & 
$0$ \\ \hline
7 & 
$\state{0}$ & 
$\state{1}$ & 
$\state{1}$ & 
$\state{0}$ & 
$x_{3} + x_{4}$ \\ \hline
8 & 
$\state{0}$ & 
$\state{1}$ & 
$\state{1}$ & 
$\state{1}$ & 
$x_{3}$ \\ \hline
9 & 
$\state{1}$ & 
$\state{0}$ & 
$\state{0}$ & 
$\state{0}$ & 
$x_{4}$ \\ \hline
10 & 
$\state{1}$ & 
$\state{0}$ & 
$\state{0}$ & 
$\state{1}$ & 
$0$ \\ \hline
11 & 
$\state{1}$ & 
$\state{0}$ & 
$\state{1}$ & 
$\state{0}$ & 
$x_{3} + x_{4}$ \\ \hline
12 & 
$\state{1}$ & 
$\state{0}$ & 
$\state{1}$ & 
$\state{1}$ & 
$x_{3}$ \\ \hline
13 & 
$\state{1}$ & 
$\state{1}$ & 
$\state{0}$ & 
$\state{0}$ & 
$x_{4}$ \\ \hline
14 & 
$\state{1}$ & 
$\state{1}$ & 
$\state{0}$ & 
$\state{1}$ & 
$0$ \\ \hline
15 & 
$\state{1}$ & 
$\state{1}$ & 
$\state{1}$ & 
$\state{0}$ & 
$x_{3} + x_{4}$ \\ \hline
16 & 
$\state{1}$ & 
$\state{1}$ & 
$\state{1}$ & 
$\state{1}$ & 
$x_{3}$ \\ \hline
\end{tabular}
\end{tabular}
\caption{Secondary analysis to distinguish humans from zombies, using
  the polynomials in the result table for \prob{R,H} shown in
  Equation~\ref{eq:trh} instantiated at different values of the
  parameters $t_1$ through $t_4$.  We search for values of
  $(t_1,t_2,t_3,t_4)$ such that exactly one of
  \prob{R=\state{T},H=\state{T}} or \prob{R=\state{T},H=\state{F}} is
  zero, and exactly one of \prob{R=\state{F},H=\state{T}} or
  \prob{R=\state{F},H=\state{F}} is zero.  The values at rows 6 and 11
  meet these criteria.}
\label{tbl:meta-zombie}
\end{table}

It takes two phases of analysis to find a question that will determine
whether the speaker is human or zombie.  For the primary analysis we
compute the joint probability \prob{R,H} of each identity and each
response using symbolic probability-network inference as in
Section~\ref{sec:symbolic}:
\begin{equation}
\begin{tabular}[c]{r|ll|l}\hline
\multicolumn{1}{l|}{\scshape {Index}} & 
\multicolumn{1}{l}{$R$} & 
\multicolumn{1}{l|}{$H$} & 
\multicolumn{1}{l}{$\prob{{R, H}}$} \\ \hline\hline
1 & 
$\state{T}$ & 
$\state{T}$ & 
$x_{2} + t_{1} x_{1} - t_{2} x_{2}$ \\ \hline
2 & 
$\state{T}$ & 
$\state{F}$ & 
$x_{3} - t_{3} x_{3} + t_{4} x_{4}$ \\ \hline
3 & 
$\state{F}$ & 
$\state{T}$ & 
$x_{1} - t_{1} x_{1} + t_{2} x_{2}$ \\ \hline
4 & 
$\state{F}$ & 
$\state{F}$ & 
$x_{4} + t_{3} x_{3} - t_{4} x_{4}$ \\ \hline
\end{tabular}
%
\label{eq:trh}
\end{equation}
For the secondary analysis we must find values of the parameters
$(t_1,t_2,t_3,t_4)$ such that the question $Q$ encoded by those values
successfully discriminates between humans and zombies.  With regard to
the result table for \prob{R,H} in Equation~\ref{eq:trh}, successful
discrimination requires that after substitution of the chosen values
of the $t_i$ parameters, exactly one of the first two elements of
\prob{R,H} is identically zero and that exactly one of the second two
elements is identically zero.  As in Section~\ref{sec:search} we can
set up a simple exhaustive search to find such values by substituting
each of the sixteen possible values of $(t_1,t_2,t_3,t_4)$ with each
$t_i \in \{0,1\}$ into each of the four polynomials in
Equation~\ref{eq:trh}.  Table~\ref{tbl:meta-zombie} shows these
substituted polynomial values.  There are two vectors of parameter
values that meet the search criteria: $(0,1,0,1)$ at row~6 and
$(1,0,1,0)$ at row~11.  With reference to Equation~\ref{eq:cpt-q},
here is the question $Q$ instantiated using the second solution
$(t_1,t_2,t_3,t_4)=(1,0,1,0)$:
\begin{equation}
\begin{tabular}[c]{ll|ll}\hline
\multicolumn{4}{l}{$\condprob[0]{{Q}}{{H, B}}$} \\ \hline\hline
\multicolumn{1}{l}{$H$} & 
\multicolumn{1}{l|}{$B$} & 
\multicolumn{1}{l}{${Q=\state{T}}$} & 
\multicolumn{1}{l}{${Q=\state{F}}$} \\ \hline\hline
$\state{T}$ & 
$\state{T}$ & 
$1$ & 
$0$ \\ \hline
$\state{T}$ & 
$\state{F}$ & 
$0$ & 
$1$ \\ \hline
$\state{F}$ & 
$\state{T}$ & 
$1$ & 
$0$ \\ \hline
$\state{F}$ & 
$\state{F}$ & 
$0$ & 
$1$ \\ \hline
\end{tabular}
\end{equation}
Using this solution the primary variable $Q$ expresses the logical
formula $B$ representing whether `Bal' means `yes'; in other words
$Q:=B$.  In this case the joint probability of identity and response,
obtained by substituting the selected parameter values into
Equation~\ref{eq:trh}, is given by:
\begin{equation}
\begin{tabular}[c]{r|ll|l}\hline
\multicolumn{1}{l|}{\scshape {Index}} & 
\multicolumn{1}{l}{$R$} & 
\multicolumn{1}{l|}{$H$} & 
\multicolumn{1}{l}{$\prob{{R, H}}$} \\ \hline\hline
1 & 
$\state{T}$ & 
$\state{T}$ & 
$x_{1} + x_{2}$ \\ \hline
2 & 
$\state{T}$ & 
$\state{F}$ & 
$0$ \\ \hline
3 & 
$\state{F}$ & 
$\state{T}$ & 
$0$ \\ \hline
4 & 
$\state{F}$ & 
$\state{F}$ & 
$x_{3} + x_{4}$ \\ \hline
\end{tabular}
\end{equation}
In other words, the question ``Does `Bal' mean `yes'?'' will reliably
distinguish a human speaker from a zombie: a human ($H=\state{T}$)
must answer `Bal' ($R=\state{T}$) but a zombie must answer `Da'.  The
other solution $(0,1,0,1)$ for the parameters $(t_1,t_2,t_3,t_4)$
gives the negation of this question which works just as well; in
response to the question ``Does `Da' mean `yes'?'' (corresponding to
the formula $\neg B$) a human must answer `Da' but a zombie must
answer `Bal'.  Notably, Kolany incorrectly matched the question from
the first solution with the responses from the second solution:
``\ldots we could ask him whether \emph{Bal} means \emph{Yes}.  If he
answers \emph{Bal}, he is a zombie.''  \cite{kolany}

Note that certain prior assumptions about whether the speaker is human
versus zombie will lead to exceptions.  Here are the computed
probability distributions \prob{H} and \prob{R}:
\begin{equation}
\begin{tabular}[c]{l|l}\hline
\multicolumn{1}{l|}{$H$} & 
\multicolumn{1}{l}{$\prob{{H}}$} \\ \hline\hline
$\state{T}$ & 
$x_{1} + x_{2}$ \\ \hline
$\state{F}$ & 
$x_{3} + x_{4}$ \\ \hline
\end{tabular}
\qquad
\begin{tabular}[c]{l|l}\hline
\multicolumn{1}{l|}{$R$} & 
\multicolumn{1}{l}{$\prob{{R}}$} \\ \hline\hline
$\state{T}$ & 
$x_{1} + x_{2}$ \\ \hline
$\state{F}$ & 
$x_{3} + x_{4}$ \\ \hline
\end{tabular}
\end{equation}
Had it been specified \emph{a priori} that there were no humans on the
island (with the constraint
\begin{math}
x_{1} + x_{2}
= 0
\end{math})
then it would be impossible under the rules of the island for the
speaker to answer `Bal'; in other words $\prob{R=\state{T}}$, which
has the polynomial value
\begin{math}
x_{1} + x_{2}
\end{math},
would also be constrained to equal zero.  Likewise had it been
specified that there were no zombies on the island (\begin{math}
x_{3} + x_{4}
= 0
\end{math})
then it would be impossible for the speaker to answer `Da'
($R=\state{F}$).

\section{Conclusion}

The method of \emph{parametric probability analysis}, introduced and
illustrated above, demonstrates that probability and classical logic
are not only compatible but also complementary.  To adopt a popular
turn of phrase, there is no daylight between logic and probability.
Many so-called `logic' problems are more specifically probability
problems, because they require reasoning about the probabilities of
formulas from the propositional calculus.  For such problems it is
useful to embed logical formulas inside parametric probability
networks.  Many other logic problems are better represented directly
as parametric probability networks, without use of the propositional
calculus or of first-order logic at all.  Parametric probability
analysis complements classical logic by providing a powerful set of
computational tools for modeling and reasoning about implication,
consequence, and quantification.

\begin{raggedright}
  \small
  \bibliographystyle{plain}
  \bibliography{../decisions}

\vspace*{\baselineskip}
\noindent
\emph{Implementation note}:
\texttt{pqlsh-8.0.7+cplex+openmp 2012-05-23}

\end{raggedright}

\clearpage
\appendix

\section{Probability Model Source Code}
\label{sec:source}

\subsection{Ace-King}

\begin{code}
\small
\begin{verbatim}
// ace-king.pql: Johnson-Laird Acta Psych 1996 via Bringsford J Applied Logic 2008

primary A { label = "There is an ace in the hand"; states = binary; }
primary K { label = "There is a king in the hand"; states = binary; }
clique _C; potential ( _C : A K ) { parametric(x); }

primary P { 
  label = "Value of $((K \leftrightarrow A) \leftrightarrow K)$"; 
  states = binary; 
}
probability ( P | A K ) { function = ( "P <-> ((K <-> A) <-> K) ? 1 : 0" ); }
\end{verbatim}
\end{code}

\subsection{Amphibian}

\begin{code}
\small
\begin{verbatim}
// amphibian.pql: adopted from Paris, Muino, Rosefield 2009

primary P { label = "Chicken killer"; states = binary; }
primary Q { label = "Japanese"; states = binary; }
primary R { label = "Salamander"; states = binary; }
// fully parametric prior distribution, no independence assumptions
clique _C; probability ( _C : P Q R ) { parametric(x); }

// beliefs (like axioms)

primary S_1 { label = "Value of $(P \wedge Q)$"; states = binary; }
probability ( S_1 | P Q ) { function = ( "S_1 <-> P && Q ? 1 : 0" ); }

primary S_2 { label = "Value of $(\neg (Q \wedge R) \wedge P)$"; states = binary; }
probability ( S_2 | P Q R ) { function = ( "S_2 <-> !(Q && R) && P ? 1 : 0" ); }
        
primary S_3 { label = "$S_3 :: R \wedge (\neg P \rightarrow (R \wedge Q))$"; states = binary; }
probability ( S_3 | P Q R ) { function = ( "S_3 <-> R && (!P -> (R && Q)) ? 1 : 0" ); }

// queries

primary S_4 { label = "Value of $(S_4 :: P \wedge R)$"; states = binary; }
probability ( S_4 | P R ) { function = ( "S_4 <-> P && R ? 1 : 0" ); }

primary S_5 { label = "Value of $(P \wedge (Q \vee R))$"; states = binary; }
probability ( S_5 | P Q R ) { function = ( "S_5 <-> P && (Q || R) ? 1 : 0" ); }

primary S_6 { label = "Value $(R)$"; states = binary; }
probability ( S_6 | R ) { function = ( "S_6 <-> R ? 1 : 0" ); }

primary S_7 { label = "Value of $(\neg R)$"; states = binary; }
probability ( S_7 | R ) { function = ( "S_7 <-> !R ? 1 : 0" ); }

primary S_8 { label = "Value of $(R \wedge \neg R)$"; states = binary; }
probability ( S_8 | R ) { function = ( "S_8 <-> R && !R ? 1 : 0" ); }

net { graph = "rankdir = TB;"; }
net { graph = 'subgraph { rank=same; "P"; "Q"; "R"}'; }
\end{verbatim}
\end{code}

\subsection{Counterfactual Conditions}

\begin{code}
\small
\begin{verbatim}
// butter.pql: Goodman's counterfactual from Fact, Fiction, and Forecast

parameter x { range = (0,1); }
parameter y { range = (0,1); }
parameter z { range = (0,1); }

primary H { label = "The butter was heated"; states = binary; }
probability ( H ) { data = (x, 1-x); }

primary M { label = "The butter melted"; states = binary; }
probability ( M | H ) { data = (y, 1-y, z, 1-z); }

primary C_1 { 
  label = "Value of $(H \rightarrow M)$"; tex = "(H \rightarrow M)"; 
  states = binary; 
}
probability ( C_1 | H M ) { function = "C_1 <-> H -> M ? 1 : 0"; }

primary C_2 { 
  label = "Value of $(H \rightarrow \neg M)$"; tex = "(H \rightarrow \neg M)";
  states = binary; 
}
probability ( C_2 | H M ) { function = "C_2 <-> H -> !M ? 1 : 0"; }

net { graph = 'subgraph { rank=same; "H"; "M"; }'; }
\end{verbatim}
\end{code}

\subsection{Knight or Knave}

\begin{code}
\small
\begin{verbatim}
// knight2.pql: Kolany's example 2, from Smullyan 1978 #36 p. 23
// I asked one of them, "Is either of you a knight?"

primary A { label = "A is a knight"; states = binary; }
primary B { label = "B is a knight"; states = binary; }
clique _C; probability ( _C : A B ) { parametric(x); }

primary Q { label = "Question: $A \vee B$"; states = binary; }
probability ( Q | A B ) { function = ( "(Q <-> A || B) ? 1 : 0" ); }

primary R { label = "A's response: $A \leftrightarrow Q$"; states = binary; }
probability ( R | Q A ) { function = ( "R <-> (Q <-> A) ? 1 : 0" ); }

net { graph = 'subgraph { rank=same; "Q"; "R"; }'; }
\end{verbatim}
\end{code}

\subsection{Human or Zombie: Primary Analysis}

\begin{code}
\small
\begin{verbatim}
// zombie1.pql: Kolany's example 3 from Smullyan 1978 #160 p. 150

primary H { label = "Speaker is human"; states = binary; }
primary B { label = "`Bal' means `yes'"; states = binary; }
primary Q { label = "Question (parametric in $t_i$)"; states = binary; }
probability ( Q | H B ) { parametric(t); }
clique _C; probability ( _C : H B ) { parametric(x); }

primary R { 
  label = "Response is `Bal': $((H \leftrightarrow Q) \leftrightarrow B)$"; 
  states = binary;
}
probability ( R | H Q B ) { function = ( "R <-> ((H <-> Q) <-> B) ? 1 : 0" ); }

net { graph = 'subgraph { rank=same; "Q"; "R"; }'; }
net { graph = 'subgraph { rank=max; "t1"; "t2"; "t3"; "t4"; }'; }
\end{verbatim}
\end{code}

\subsection{Human or Zombie: Secondary Analysis}

\begin{code}
\small
\begin{verbatim}
// zombie1-search.pql: Kolany's example 3 from Smullyan 1978 #160 p. 150

decision t[1] { states = values(0,1); }
decision t[2] { states = values(0,1); }
decision t[3] { states = values(0,1); }
decision t[4] { states = values(0,1); }

// sequence of decisions; not important here
probability ( t[1] ) {}
probability ( t[2] | t[1] ) {}
probability ( t[3] | t[2] ) {}
probability ( t[4] | t[3] ) {}

// to allow substitution of t[i] values; not really decision-theoretic utilities
utility U_1 { tex = "\prob{R=\state{T},H=\state{T}}"; range = (0,1); }
utility U_2 { tex = "\prob{R=\state{T},H=\state{F}}"; range = (0,1); }
utility U_3 { tex = "\prob{R=\state{F},H=\state{T}}"; range = (0,1); }
utility U_4 { tex = "\prob{R=\state{F},H=\state{F}}"; range = (0,1); }

// polynomials are from the result Pr( R, H ) using zombie1.pql
probability ( U_1 | t[1] t[2] t[3] t[4] ) { function = "x2 + t1*x1 - t2*x2"; }
probability ( U_2 | t[1] t[2] t[3] t[4] ) { function = "x3 - t3*x3 + t4*x4"; }
probability ( U_3 | t[1] t[2] t[3] t[4] ) { function = "x1 - t1*x1 + t2*x2"; }
probability ( U_4 | t[1] t[2] t[3] t[4] ) { function = "x4 + t3*x3 - t4*x4"; }

primary H { label = "Speaker is human"; states = binary; }
primary B { label = "`Bal' means `yes'"; states = binary; }

// Question 11: Does `Bal' mean `yes'?
set t[1] = 1; set t[2] = 0; set t[3] = 1; set t[4] = 0;

primary Q { label = "Question (parametric in $t_i$)"; states = binary; }
probability ( Q | H B ) { parametric(t); }
clique _C; probability ( _C : H B ) { parametric(x); }

primary R { 
  label = "Response is `Bal': value of $(Q \leftrightarrow (H \leftrightarrow B))$"; 
  states = binary;
}
probability ( R | Q H B ) { function = ( "R <-> (Q <-> (H <-> B)) ? 1 : 0" ); }

net { graph = 'subgraph { rank=same; "Q"; "R"; }'; }
net { graph = 'subgraph { rank=max; "t1"; "t2"; "t3"; "t4"; }'; }
\end{verbatim}
\end{code}

\clearpage
\tableofcontents

\end{document}